\documentclass[12pt,leqno]{amsart}
\usepackage{amssymb}
\usepackage{euscript}
\usepackage[hypertex]{hyperref}

\catcode`\@=11

\title[Operators and graph complexes]%
       {Natural differential operators and graph complexes}

\author[M.~Markl]{Martin~MARKL}
\thanks{Supported by the grant GA \v CR 201/05/2117 and by
   the Academy of Sciences of the Czech Republic,
   Institutional Research Plan No.~AV0Z10190503}

\address{Mathematical Institute of the Academy, {\v Z}itn{\'a} 25,
         115 67 Prague 1, The Czech Republic}
\email{markl@math.cas.cz}

\begin{document}

\bibliographystyle{plain}

\swapnumbers
\newtheorem{theorem}{Theorem}[section]
\newtheorem{corollary}[theorem]{Corollary}
\newtheorem{observation}[theorem]{Observation}
\newtheorem{lemma}[theorem]{Lemma}
\newtheorem{proposition}[theorem]{Proposition}
\newtheorem{problem}[theorem]{Problem}
\newtheorem{conjecture}[theorem]{Conjecture}
\newtheorem{odstavec}{\hskip -0.1mm}[section]
\newtheorem*{principle}{Principle}
\newtheorem*{itt}{Invariant Tensor Theorem}

\theoremstyle{definition}
\newtheorem{example}[theorem]{Example}
\newtheorem{remark}[theorem]{Remark}
\newtheorem{definition}[theorem]{Definition}

\def\Rada#1#2#3{#1_{#2},\dots,#1_{#3}} \def\GLV{{\GLname(V)}}
\def\aRda#1#2#3{#1^{#2},\dots,#1^{#3}} \def\bfk{{\mathbf k}}
\def\pa{\partial} \def\Con{{\it Con\/}} \def\card{{\rm card}}
\def\GLname{{\rm GL\/}}\def\GL#1#2{{\GLname\/}^{(#1)}_{#2}}
\def\NGLname{{\rm NGL\/}}\def\NGL#1#2{{\NGLname\/}^{(#1)}_{#2}}
\def\nglname{{\mathfrak {ngl}\/}}\def\ngl#1#2{{\nglname\/}^{(#1)}_{#2}}
\def\JGLname{{\rm J\GLname\/}}  \def\pg{{\mathbold p}}
\def\Frname{{\it Fr\/}}\def\IN{{\rm In}}\def\OUT{{\rm Ou}}
\def\fS{{\mathfrak S}} \def\fs{{\mathfrak s}}\def\hh{h}
\def\Ant{{\rm Ant\/}}\def\epi{ \twoheadrightarrow}\def\frH{{\mathfrak H}}
\def\fSout{\fS_{\it out}} \def\fSin{\fS_{\it in}}\def\Grbull{\Gr_\bullet}
\def\Fr#1{{\Frname}^{(#1)}} \def\wR{{\widehat R}} \def\wGr{{\widehat{\Gr}}}
\def\fin{{\mathfrak{{I}}}} \def\udelta{{\underline{\delta}}}
\def\Si{{\mathfrak{{I}}}} \def\So{{\mathfrak{{O}}}}
\def\wrR{\widehat{\mathrm R}}\def\uR{{\underline{R}}}
\def\adjust#1{{\raisebox{-#1em}{\rule{0pt}{0pt}}}} \def\adj{\adjust {.4}}
\def\fout{{\mathfrak{ou}}} \def\rR{{\mathrm R}}
\def\bbbR{{\mathbb R}} \def\bbR{\bbbR}
\def\dr#1#2{\frac{\partial #1}{\partial #2}}
\def\Man{{\tt Man\/}} \def\ctverecek{\raisebox{-0pt}{\smbbox}\hskip -4pt}
\def\gF{{\mathfrak F\/}}\def\gG{{\mathfrak G\/}}\def\gO{{\mathfrak O\/}}
\def\gB{{\mathfrak B\/}}\def\sB{{\EuScript B}}
\def\Fib{{\tt Fib\/}} \def\symexp#1#2{{#1^{\odot #2}}}
\def\Lin{{\mbox {\it Lin\/}}}\def\Sym{{\mbox {\it Sym\/}}}
\def\ext{\mbox{\Large$\land$}}\def\bp{{\mathbf p}}
\def\EXT{\mbox{\raisebox{-.6em}{\rule{1pt}{0pt}}%
         \raisebox{-.2em}{\Huge$\land$}}}
\def\LAND{\mbox{\raisebox{-.2em}{\rule{1pt}{0pt}}%
         \raisebox{-.0em}{\Large$\land$}}}
\def\Jo{{\overline{\JGLname}}}\def\uNat{\underline{\mathfrak {Nat}\/}} 
\def\jo{{\overline {\j}}} \def\CCE{C_{\it CE}}
\def\semidirect{\makebox{\hskip .3mm$\times \hskip -.8mm\raisebox{.2mm}%
                {\rule{.17mm}{1.9mm}}\hskip-.75mm$\hskip 2mm}}
\def\Rn{{\bbR^n}} \def\ot{\otimes} \def\zn#1{{(-1)^{#1}}}
\def\id{{1 \!\! 1}} \def\gh{{\mathfrak h}}
\def\sF{{\EuScript F}}\def\sG{{\EuScript G}}\def\sC{{\EuScript C}}
\def\Map{{\mbox {\it Map\/}}} \def\Nat{{\mathfrak {Nat}\/}} 
\def\WFG{W_{\gF,\gG}}  \def\sqot{\hskip -.3em \otimes \hskip -.3em}
\def\Gr{{\EuScript {G}\rm r}}\def\plGr{{\rm pl\widehat{\EuScript {G}\rm r}}}
\def\GrFG{{\Gr}_{\gF,\gG}}   \def\orGr{{\rm or\EuScript {G}\rm r}}
\def\otexp#1#2{#1^{\ot #2}} \def\semGH{G \semidirect H}
\def\sbbox{{\raisebox {.1em}{\rule{.4em}{.4em}} \hskip .1em}}
\def\bbox{{\raisebox {.1em}{\rule{.6em}{.6em}} \hskip .1em}}
\def\mbbox{{\raisebox {.1em}{\rule{.4em}{.4em}} \hskip .1em}}
\def\smbbox{{\raisebox {.0em}{\rule{.5em}{.5em}} \hskip .0em}}
\def\vmbbox{{\raisebox {.1em}{\rule{.2em}{.2em}} \hskip .1em}}
\def\ii{{(\infty)}} \def\phi{\varphi}\def\uGrFG{{\uGr}_{\gF,\gG}} 
\def\rada#1#2{{#1,\ldots,#2}}\def\uGr{{\EuScript {G}\underline{\rm r}}}
\def\dirlim{{{\mathop{{\rm lim}}\limits_{\longrightarrow}}\hskip 1mm}}
\def\dCE{\delta_{\it CE\/}} \def\plR{{\mathrm {pl}\widehat{\EuScript R}}}  
\def\orR{{\mathit orR}} \def\fA{{\widehat{\Sigma}}}
\def\Vert{{\it Vert\/}} \def\bfb{{\mathbf b}}
\def\Ker{{\it Ker}} \def\wc{\circlearrowright}
\def\Lie{{\mathcal L{\it ie\/}}}\def\Im{{\it Im}}
\def\Ext{\mathop{{\rm \EXT}}\displaylimits}
\def\Land{\mathop{{\LAND}}\displaylimits}
\def\skelet{
\begin{picture}(2,4)(0,0)
\put(1,1){\oval(2,2)[t]}
\put(1,-1){\oval(2,2)[b]}
\put(2,0.25){\makebox(0,0)[c]{\vdots}}
\end{picture}
}

\def\subb#1{\hskip -.2em {\raisebox{-.25em}{\scriptsize $#1$}}\hskip .2em}
\def\anchor{$\unitlength .25cm
\begin{picture}(1,1.4)(-1,-.7)
\put(-.45,.55){\makebox(0,0)[cc]{$\sbbox$}}
\put(-.5,-.8){\vector(0,1){1.2}}
\end{picture}$}

\def\uunit{$\unitlength .25cm
\begin{picture}(1,1.4)(-1,-.9)
\put(-.45,.55){\makebox(0,0)[cc]{$\sbbox$}}
\put(-.45,-1){\makebox(0,0)[cc]{$\bullet$}}
\put(-.5,-.8){\vector(0,1){1.2}}
\end{picture}$}

\def\black{$\unitlength .25cm
\begin{picture}(1,1.4)(-1,-0.4)
\put(-.45,-.3){\makebox(0,0)[cc]{\Large$\bullet$}}
\put(-.5,0){\vector(0,1){1.3}}
\end{picture}\hskip .2em $}

\def\white{$\unitlength .25cm
\begin{picture}(1,1.4)(-1,-.7)
\put(-.45,-.3){\makebox(0,0)[cc]{\Large$\circ$}}
\put(-.5,0.1){\vector(0,1){1.2}}
\end{picture}\hskip .2em $}

\def\osG{\overline{\Gr}}    \def\ti{{\times \infty}}      \def\td{{\times d}}
\def\Tr{{\it Tr\/}}         \def\oti{{\otimes \infty}}    
\def\Com{{\EuScript C}{\it om}} \def\calA{{\EuScript A}} 
\def\bfc{{\mathbf c}}       \def\calQ{{\EuScript Q}}
\def\calP{{\EuScript P}}     \def\pLie{p{\mathcal L{\it ie\/}}}
\def\Re{{\mathbb R}}        \def\od{{\otimes d}}
\def\crr{connected replacement rules}
\def\odrazka#1{{\raisebox{#1}{\rule{0pt}{0pt}}}}
\def\ccdot{\mbox {\scriptsize \hskip .3em $\bullet$\hskip .3em}}
\def\Grtr{\Gr_{\bullet\nabla\it Tr}}      \def\vt{\vartheta}
\def\Grd #1#2#3{\Gr^{#1}_{\bullet#2}[#3](d)} \def\Edg{{\it Edg}}
\def\Lab{{\it Lab}}

\def\borelioza#1#2{
\unitlength.7cm
\put(0,-1){
\put(0,-.1){
\put(0,2){\put(0.03,0){\makebox(0,0)[cc]{$\bbox$}}}
\put(0,1){\vector(0,1){.935}}
\put(0,1){\makebox(0,0)[cc]{\Large$\bullet$}}
\put(0,.7){\makebox(0,0)[t]{\scriptsize$#1$}}}
\put(.4,.2){
\put(2.09,.85){\makebox(0,0)[cc]{\oval(1.5,1.5)[b]}}
\put(2.09,1.15){\makebox(0,0)[cc]{\oval(1.5,1.5)[t]}}
\put(2.85,1.22){\line(0,1){.3}}
\put(.7,.7){\vector(1,1){.55}}
\put(.7,.7){\makebox(0,0)[cc]{\Large$\bullet$}}
\put(1.35,1.35){\makebox(0,0)[cc]{\Large$\bullet$}}
\put(1.32,1.25){\makebox(0,0)[tc]{\vector(0,1){0}}}
\put(1,1.45){\makebox(0,0)[r]{\scriptsize $F$}}
\put(.7,0.4){\makebox(0,0)[t]{\scriptsize $#2$}}
}}}
\def\cases#1#2#3#4{
                  \left\{
                         \begin{array}{ll}
                           #1,\ &\mbox{#2}
                           \\
                           #3,\ &\mbox{#4}
                          \end{array}
                   \right.
}
\def\boreliozaInv#1#2{
\unitlength.7cm
\put(0,-1){
\put(0,-.1){
\put(0,2){\put(0.03,0){\makebox(0,0)[cc]{$\bbox$}}}
\put(0,1){\vector(0,1){.935}}
\put(0,1){\makebox(0,0)[cc]{\Large$\bullet$}}
\put(0,.7){\makebox(0,0)[t]{\scriptsize$#1$}}
}
\put(3.8,.2){
\put(-2.09,.85){\makebox(0,0)[cc]{\oval(1.5,1.5)[b]}}
\put(-2.09,1.15){\makebox(0,0)[cc]{\oval(1.5,1.5)[t]}}
\put(-2.85,1.22){\line(0,1){.3}}
\put(-.7,.7){\vector(-1,1){.55}}
\put(-.7,.7){\makebox(0,0)[cc]{\Large$\bullet$}}
\put(-1.35,1.35){\makebox(0,0)[cc]{\Large$\bullet$}}
\put(-1.32,1.25){\makebox(0,0)[tc]{\vector(0,1){0}}}
\put(-1,1.45){\makebox(0,0)[l]{\scriptsize $F$}}
\put(-.7,0.4){\makebox(0,0)[t]{\scriptsize $#2$}}
}}}

\def\sigmadva{{\unitlength.5cm\thicklines
\put(0,-.5){\vector(1,1){1}}
\put(1,-.5){\vector(-1,1){1}}
\put(2,-.5){\vector(0,1){1}}
\put(0,-.7){\makebox(0,0)[t]{\scriptsize 1}}
\put(1,-.7){\makebox(0,0)[t]{\scriptsize 2}}
\put(2,-.7){\makebox(0,0)[t]{\scriptsize 3}}
\put(0,.7){\makebox(0,0)[b]{\scriptsize 1}}
\put(1,.7){\makebox(0,0)[b]{\scriptsize 2}}
\put(2,.7){\makebox(0,0)[b]{\scriptsize 3}}
}}

\def\sigmatri{{\unitlength.5cm\thicklines
\put(0,-.5){\vector(0,1){1}}
\put(1,-.5){\vector(1,1){1}}
\put(2,-.5){\vector(-1,1){1}}
\put(0,-.7){\makebox(0,0)[t]{\scriptsize 1}}
\put(1,-.7){\makebox(0,0)[t]{\scriptsize 2}}
\put(2,-.7){\makebox(0,0)[t]{\scriptsize 3}}
\put(0,.7){\makebox(0,0)[b]{\scriptsize 1}}
\put(1,.7){\makebox(0,0)[b]{\scriptsize 2}}
\put(2,.7){\makebox(0,0)[b]{\scriptsize 3}}
}}

\def\sigmactyri{{\unitlength.5cm\thicklines
\put(0,-.5){\vector(2,1){2}}
\put(1,-.5){\vector(0,1){1}}
\put(2,-.5){\vector(-2,1){2}}
\put(0,-.7){\makebox(0,0)[t]{\scriptsize 1}}
\put(1,-.7){\makebox(0,0)[t]{\scriptsize 2}}
\put(2,-.7){\makebox(0,0)[t]{\scriptsize 3}}
\put(0,.7){\makebox(0,0)[b]{\scriptsize 1}}
\put(1,.7){\makebox(0,0)[b]{\scriptsize 2}}
\put(2,.7){\makebox(0,0)[b]{\scriptsize 3}}
}}

\def\sigmapet{{\unitlength.5cm\thicklines
\put(0,-.5){\vector(1,1){1}}
\put(1,-.5){\vector(1,1){1}}
\put(2,-.5){\vector(-2,1){2}}
\put(0,-.7){\makebox(0,0)[t]{\scriptsize 1}}
\put(1,-.7){\makebox(0,0)[t]{\scriptsize 2}}
\put(2,-.7){\makebox(0,0)[t]{\scriptsize 3}}
\put(0,.7){\makebox(0,0)[b]{\scriptsize 1}}
\put(1,.7){\makebox(0,0)[b]{\scriptsize 2}}
\put(2,.7){\makebox(0,0)[b]{\scriptsize 3}}
}}

\def\sigmasest{{\unitlength.5cm\thicklines
\put(0,-.5){\vector(2,1){2}}
\put(1,-.5){\vector(-1,1){1}}
\put(2,-.5){\vector(-1,1){1}}
\put(0,-.7){\makebox(0,0)[t]{\scriptsize 1}}
\put(1,-.7){\makebox(0,0)[t]{\scriptsize 2}}
\put(2,-.7){\makebox(0,0)[t]{\scriptsize 3}}
\put(0,.7){\makebox(0,0)[b]{\scriptsize 1}}
\put(1,.7){\makebox(0,0)[b]{\scriptsize 2}}
\put(2,.7){\makebox(0,0)[b]{\scriptsize 3}}
}}

\def\brace{
\put(0,0){\line(0,-1){2}}
\put(0,0){\vector(-1,0){.2}}
\put(0,-2){\vector(-1,0){.2}}
}

\begin{abstract}
We show how the machine invented by
S.~Merkulov~\cite{merkulov-operads,merkulov:defquant,merkulov:PROP}
can be used to study and classify natural operators in differential
geometry. We also give an interpretation of graph complexes arising in
this context in terms of representation theory. As application, we
prove several results on classification of natural operators acting on
vector fields and connections.
\end{abstract}

\maketitle

\baselineskip15.8pt plus 0pt minus .5pt
\section*{Introduction}

This work started in an attempt to understand S. Merkulov's idea of 
``PROP profiles''~\cite{merkulov-operads,merkulov:PROP} and see 
if and how it may be used to investigate natural structures  in geometry.
It turned out that classifications of these geometric structures 
in many interesting cases boiled down to calculations of the cohomology of
certain graph complexes. More precisely, for a wide class of
natural operators, the following principle holds.

\begin{principle}
For a given type of natural differential operators, there exists a graph
cochain complex
$(\Gr^*,\delta) = (\Gr^0 \stackrel{\delta}{\to} \Gr^1
\stackrel{\delta}{\to} \Gr^2 \stackrel{\delta}{\to} \cdots )$
such that, in stable ranges,
\[
\left\{ \mbox {natural operators of a given type\/} \right\}
\cong
H^0(\Gr_*,\delta).
\]
\end{principle}

{\em Stability\/} means that the dimension of the underlying manifold
is bigger than some constant explicitly determined by the type of
operators. For example, for multilinear
natural operators $TM^{\times d} \to TM$ from the $d$-fold product of
the tangent bundle into itself the stability means that $\dim(M) \geq
d$. In smaller dimensions, ``exotic'' operations described in~\cite{dzh} occur.

In all cases we studied, the corresponding graph complex appeared to be
acyclic in positive dimensions, so
the cohomology describing natural operators was the only nontrivial
piece of the cohomology of $(\Gr_*,\delta)$.
Standard philosophy of strongly homotopy structures~\cite{markl:shalg} suggests
that the graph complex $(\Gr^*,\delta)$ describes {\em stable strongly homotopy
operators\/} of a given type.

Graph complexes arising in the Principle are in fact isomorphic to subspaces
of fixed elements in suitable Chevalley-Eilenberg complexes, so,
formally speaking, we claim  that a certain Chevalley-Eilenberg
cohomology is  the cohomology of some graph complex.  Instances of
this  phenomenon were 
systematically used by M.~Kontsevich in his
seminal paper~\cite{kontsevich:93}. The details of operadic graph
complexes were then written down by
J.~Conant~\cite{conant}, J.~Conant and
K.~Vogtmann~\cite{conant-vogtman:MA03,conant-vogtmann},  M.~Mulase and
M.~Penkava~\cite{mulase-penkava:AJM98}, 
M.~Penkava~\cite{penkava:infinity}, and M.~Penkava and
A.~Schwarz~\cite{penkava-schwarz:TAMS95}.  
What makes the Principle exciting is the miraculous fact 
that the corresponding graph complexes are of the type
studied during the ``renaissance of operads'' and powerful methods
developed in this period culminating in~\cite{mms,mv,merkulov:duke} apply.

Another way to view the proposed method is as a formalization of the
``abstract tensor calculus'' attributed to R.~Penrose.  When we studied
differential geometry in kindergarten, many of us, trying to avoid
dozens of indices, drew simple pictures consisting of nodes
representing tensors (which resembled little insects) and lines
joining legs of these insects symbolizing contraction of indices. We
attempt to put this kindergarten approach on a solid footing.

Thus the purpose of this paper is two-fold. 
The first one is to  set up principles of
abstract tensor calculus as a useful language for
`stable' geometric objects. This will be done in
\hbox{Sections~\ref{1}--\ref{sec4}}. The logical continuation should be 
translating textbooks on differential geometry into this language,
because all basic properties of fundamental objects 
(vector fields, forms, currents, connections and
their torsions and curvatures) are of stable nature.

We then show, in Sections~\ref{5}--\ref{7}, 
how results on graph complexes may give 
explicit classifications of natural operators in stable ranges.
As an example we derive from a rather deep result of~\cite{markl:JLT07} a
characterization of operators on vector fields (Theorem~\ref{main-A}
and its Corollary~\ref{B}). As another
application we prove that all natural
operators on linear connections and vector fields, with values in
vector fields, are freely generated by
compositions of covariant derivatives and Lie
brackets, and by traces of these compositions -- see
Theorems~\ref{P1} and~\ref{T3}, and their Corollaries~\ref{T2} and~\ref{C2},
in conjunction with Theorems~\ref{X} and~\ref{Y}. 

This article is supplemented by~\cite{markl:ig} in which we explain the
relation between invariant tensors and graphs. We believe
that~\cite{markl:ig}, which can be read independently, will help to
understand the constructions of Sections~\ref{3} and~\ref{sec4}.

The theory of invariant operators sketched out  
in this paper leads to directed, not
necessarily connected or simply-connected, graphs. A similar theory
can be formulated also for symplectic manifolds, where the
corresponding graph complexes would be those appearing in the context of
anti-modular operads (modular versions of anticyclic operads, 
see~\cite[Definition~5.20]{markl-shnider-stasheff:book}). 
Something close to a symplectic version of our theory has in
fact already been worked out in~\cite{weingart}.

\vskip .4em

\noindent
{\bf Acknowledgment.}  I would like to express my thanks to
S.~Merkulov for sharing his ideas with me, to A.~Alekseev for
suggesting an interpretation of the homological vector field in terms
of the Chevalley-Eilenberg differential, and to G.~Weingart who
pointed some flaws in my reasoning to me. Also conversations with J.~Jany\v
ska and J.~Slov\'ak were extremely useful. Suggestions of the referee
lead to a substantial improvement of the paper.

\section{Natural operators}
\label{1}

Informally, a natural differential operator is a recipe that constructs
from a geometric object another one, in a natural fashion, and which
is locally a function of coordinates and their derivatives.

\begin{example}
\label{Lie}
Let $M$ be a $n$-dimensional smooth manifold. 
The classical {\em Lie bracket\/} $X,Y
\mapsto [X,Y]$ is a natural operation that constructs from two
vector fields on $M$ a third one.
Given a local coordinate system $(\aRda x1n)$ on $M$, the vector fields
$X$ and $Y$ are locally expressions 
$X = \sum_{1 \leq i \leq n} X^i {\pa}/{\pa x^i},\
Y = \sum_{1 \leq i \leq n} Y^i {\pa}/{\pa x^i}$,
where $X^i, Y^i$ are smooth functions on $M$. If we define $X^i_j :=
{\pa X^i}/{\pa x^j}$ and $Y^i_j :=
{\pa Y^i}/{\pa x^j}$, $1 \leq i,j \leq n$,
then the Lie bracket is locally given by the formula
$[X,Y] = \sum_{1 \leq i,j \leq n}
\left(X^j Y^i_j - Y^jX^i_j\right){\pa}/{\pa x^i}$.
\end{example}

In the rest of the paper, we use Einstein's convention assuming
summations over repeated indices. In this context,
indices $i,j,k,\ldots$ will always be natural numbers 
between $1$ and the dimension of the underlying manifold, which will
typically be denoted $n$.

\begin{example}
\label{cov}
The {\em covariant derivative\/}
$(\Gamma, X, Y) \mapsto \nabla_X Y$ is a 
natural operator that constructs from a linear
connection $\Gamma$ and vector fields $X$ and $Y$, 
a vector field $\nabla_X Y$. 
In local coordinates,
\begin{equation}
\label{Sergej_je_magor.}
\nabla_X Y = \left(\Gamma^i_{jk} X^j Y^k + X^j {Y^i_j}
              \right) \dr{}{x^i},
\end{equation}
where $\Gamma^i_{jk}$ are Christoffel symbols.

Natural operations can be composed into more complicated
ones. Examples of  `composed' operations are the {\em torsion\/}
$ T(X,Y) := \nabla_XY - \nabla_Y X - [X,Y] $ and the {\em curvature\/}
$ R(X,Y)Z := \nabla_{[X,Y]} Z - [\nabla_X,\nabla_Y] Z $ of the linear
connection $\Gamma$.
\end{example}

\begin{example}
\label{Jituska_je_moje_pusinka.}
Let $X$ be a vector field and $\omega$ a 1-form on $M$. Denote by
$\omega(X)\in C^\infty(M)$ the evaluation of the form $\omega$ on $X$. Then 
$(X,\omega) \mapsto  \exp(\omega(X))$ defines a natural
differential operator with values in smooth functions. Clearly, the exponential
can be replaced by an arbitrary smooth function $\phi : \bbR \to \bbR$,
giving rise to a natural operator $\gO_\phi(X,\omega) : = \phi(\omega(X))$.
\end{example}

\begin{example}
\label{Cesti}
`Randomly' generated local formulas need not lead to
natural operators. As we will see later, neither
$O_1(X,Y) = X^1_{3} Y^4 \pa /{\pa x^2}$ nor
$O_2(X,Y) = X^j Y^i_j {\pa}/{\pa x^i}$
behaves properly under coordinate changes, so they do not
give rise to vector-field valued natural operators.
\end{example}

We may summarize the above examples by saying that a natural
differential operator is a~recipe given locally as a smooth function in
coordinates and their derivatives, such that the local formula is
invariant under coordinate changes.
After this motivation, we give precise definitions of geometric
objects and operators between them. Our
exposition follows~\cite{katsylo-timashev}, see
also~\cite{kolar-michor-slovak}. 

Denote by $\Man_n$ the category of $n$-dimensional manifolds and open
embeddings. Let $\Fib_n$ be the category of smooth fiber bundles over
$n$-dimensional manifolds with morphisms differentiable maps covering
morphisms of their bases in $\Man_n$.

\begin{definition}
\label{Bonn}
A {\em natural bundle\/} is a functor $\gB : \Man_n \to \Fib_n$ such
that for each $M \in \Man_n$, $\gB(M)$ is a bundle over $M$. Moreover,
$\gB(M')$ is the restriction of $\gB(M)$ for each open submanifold $M'
\subset M$, the map $\gB(M') \to \gB(M)$ induced by $M'
\hookrightarrow M$ being the inclusion $\gB(M') \hookrightarrow \gB(M)$.
\end{definition}

Let us recall a structure theorem for natural bundles
due to Krupka, Palais and Terng~\cite{krupka:AMB78,palais-terng,terng:AMJ78}.
For each $s \geq 1$ we denote by $\GL sn$ the group of $s$-jets of
local diffeomorphisms $\bbR^n \to \bbR^n$ at~$0$, 
so that $\GL 1n$ is the ordinary general
linear group $\GLname_n$ of linear invertible maps $A : \bbR^n \to \bbR^n$.
Let  $\Fr s(M)$ be the bundle of $s$-jets of frames on $M$
whose fiber over $z \in M$ consist of $s$-jets of local
diffeomorphisms of neighborhoods of $0 \in \bbR^n$ with  
neighborhoods of $z\in M$. It is clear that $\Fr s(M)$ is a principal
$\GL sn$-bundle and
$\Fr 1(M)$ the ordinary $\GLname_n$-bundle of frames~$\Frname(M)$.

\begin{theorem}[Krupka, Palais, Terng]
For each natural bundle $\gB$, there exists $l \geq 1$ and 
a~manifold $\sB$ with a smooth $\GL ln$-action such that there is a functorial
isomorphism 
\begin{equation}
\label{porad_mne_to_trochu_poboliva}
\gB(M) \cong \Fr l(M) \times_{\GL ln} \sB: =
(\Fr l(M) \times \sB)/\GL ln.
\end{equation}
\end{theorem}

Conversely, each smooth $\GL ln$-manifold $\sB$ induces, 
via~(\ref{porad_mne_to_trochu_poboliva}), a natural
bundle $\gB$. We will call $\sB$ the {\em fiber\/} of the natural 
bundle~$\gB$. If the action of $\GL ln$ on $\sB$ does not reduce  to
an action of the quotient $\GL {l-1}n$ we say that $\gB$ has {\em
order $l$\/}.

\begin{example}
Vector fields are sections of the {\em tangent bundle\/} $T(M)$. The
fiber of this bundle is $\bbbR^n$, with the standard action of
$\GLname_n$. The description $ T(M) \cong \Frname(M)
\times_{\GLname_n} \bbbR^n $ is classical.
\end{example}

\begin{example}
\label{Stava_se_ze_mne_celebrita.}
De~Rham $m$-forms are sections of the bundle $\Omega^m(M)$ whose fiber
is the space of anti-symmetric $m$-linear maps
$\Lin(\ext^m(\Rn),\bbR)$, with the obvious induced
$\GLname_n$-action. The~presentation $ \Omega^m(M) \cong \Frname(M)
\times_{\GLname_n}\Lin(\ext^m(\Rn),\bbR)$ is also classical.
A particular case is $\Omega^0(M) \cong \Frname(M) \times_{\GLname_n}
\bbR \cong M \times \bbR$, the bundle whose sections are smooth
functions. We will denote this natural bundle by $\Re$,
believing there will be no confusion with the symbol for the reals.
\end{example}

\begin{example}
\label{Co_se_deje_v_byte?}
Linear connections are sections of the 
{\em bundle of connections\/} 
$\Con(M)$~\cite[Section~17.7]{kolar-michor-slovak} 
which we recall below. Let us
first describe the group $\GL 2n$. Its elements are
expressions of the form $A = A_1 + A_2$, where $A_1 : \Rn \to \Rn$
is a linear invertible map and $A_2$ is a linear map from the
symmetric product  $\Rn \odot \Rn$ to  $\Rn$.
The multiplication in $\GL 2n$ is given by
\[
(A_1 + A_2)(B_1 + B_2) := A_1(B_1) + A_1(B_2) + A_2(B_1,B_1). 
\]
The unit of $\GL 2n$ is ${\it id}\/_\Rn + 0$ 
and the inverse is given by the formula
\[
(A_1 + A_2)^{-1} = A_1^{-1} - A^{-1}_1 (A_2(A^{-1}_1,A^{-1}_1)).
\]
   
Let $\sC$ be the space of linear maps $\Lin(\Rn \ot \Rn,\Rn)$,
with the left action of $\GL 2n$ given as
\begin{equation}
\label{Napise_mi_nekdo?}
(A f) (u \ot v) :=  A_1 f (A^{-1}_1(u),A^{-1}_1(v)) -
A_2(A^{-1}_1(u),A^{-1}_1(v)),  
\end{equation}
for $f \in \Lin(\Rn \ot \Rn,\Rn)$, $A = A_1 + A_2 \in \GL 2n$ and $u,v
\in \Rn$. The bundle of connections is then the order~$2$ natural
bundle represented as 
$\Con(M) := \Fr 2 (M) \times \subb{\GL 2n} \sC$.
Observe that, while the action of $\GL 2n$ on the vector space $\sC$
is not linear, the restricted action of $\GLname_n \subset \GL 2n$ on
$\sC$ is the standard action of the general linear group on the space
of bilinear maps.
\end{example}

For $k \geq 0$ we denote by $\gB^{(k)}$ the bundle of $k$-jets of
local sections of the natural bundle $\gB$ so that $\gB^{(0)}= \gB$.
If $\gB$ is represented as in~(\ref{porad_mne_to_trochu_poboliva}),
then $\gB^{(k)}(M) \cong \Fr {k +l}(M) \times \subb{\GL {k +l}n}
\sB^{(k)}$, where $\sB^{(k)}$ is the space of $k$-jets of local
diffeomorphisms $\Rn \to \sB$ defined in a neighborhood of $0 \in
\Rn$.

\begin{definition}
\label{uuiiu}
Let $\gF$ and $\gG$ be natural bundles.
A (finite order) {\em natural differential operator\/} $\gO : \gF \to \gG$ is a
natural transformation (denoted by the same symbol) 
$\gO : \gF^{(k)} \to \gG$, for some $k \geq 1$. We denote the 
space of all natural differential operators $\gF \to \gG$ by
$\Nat(\gF,\gG)$.
\end{definition}

If $\gF$ and $\gG$ are natural bundles of order
$\leq l$, with fibers $\sF$ and $\sG$, respectively, then each natural
operator in Definition~\ref{uuiiu} is 
induced by an $\GL{k+l}n$-equivariant map $O: \sF^{(k)} \to \sG$, for
some $k \geq 0$.
Conversely, such an equivariant map induces an
operator $\gO : \gF \to \gG$. This means that the study of natural
operators is reduced to the study of
equivariant maps. The procedure described above is therefore called
the {\em IT reduction\/} (from invariant-theoretic).

From this moment on, we impose the following assumptions on natural
bundles $\gF$, $\gG$  an operators $\gO : \gF \to \gG$ between them.

\label{page}
\begin{itemize}
\item[{\bf A1}] 
The fibers $\sF$ and $\sG$ of the  bundles $\gF$ and $\gG$ are vector
spaces and the restricted actions of $\GLname_n \subset \GL ln$ on
$\sF$ and $\sG$ are rational linear representations,

\item[{\bf A2}] 
the action of $\GL ln$ on the fiber $\sG$ of $\gG$ is linear, and

\item[{\bf A3}] 
we consider only {\em polynomial} differential
operators for which the induced map of the
fibers $O : \sF^{(k)} \to \sG$ is a polynomial map. 
\end{itemize}

Notice that we do not require the action of the {\em full\/} group $\GL ln$
on the fiber of $\gF$ to be linear. 
Assumption~A2 is needed for the cohomology in
Theorem~\ref{Zitra_jedu_do_Hamburku.} in Section~\ref{2} to be 
well-defined, assumptions~A1 and~A3 are necessary to relate this
cohomology to a graph complex.

Polynomiality A3 rules out operators as $\gO_\phi$ from
Example~\ref{Jituska_je_moje_pusinka.}. There is probably no
systematic way how to study operators of this type -- imagine that
$\phi$ is an arbitrary, not even real analytic, smooth function. Clearly most
if not all ``natural'' natural operators considered in differential
geometry are polynomial, so assumption~A3 seems to be
justified. As argued in~\cite[Section~24]{kolar-michor-slovak}
and as we will also see later in Remarks~\ref{pol1} and~\ref{pol2}, 
in some situations the operators possess a certain homogeneity 
which automatically implies 
polynomiality.

\begin{example}
Given natural bundles $\gB'$ and $\gB''$ with fibers $\sB'$ resp.\
$\sB''$, there is an obviously defined natural bundle $\gB' \times
\gB''$ with fiber $\sB' \times \sB''$. With this notation, the Lie
bracket is a natural operator $[-,-] : T \times T \to T$ and the
covariant derivative an operator $\nabla : \Con \times T \times T \to
T$, where $T$ is the tangent space functor and $\Con$ the bundle of
connections recalled in Example~\ref{Co_se_deje_v_byte?}. The
corresponding equivariant maps of fibers can be easily read off from
local formulas given in Examples~\ref{Lie} and~\ref{cov}.
\end{example}

\begin{example}
The operator $\gO_{\phi} : T \times   \Omega^1 \to C^\infty$ from
Example~\ref{Jituska_je_moje_pusinka.} is induced by
the $\GLname_n$-equivariant map $O_{\phi}: \Rn \times \Rn^* \to \bbR$
given by $o_\phi(v,\alpha) := \phi(\alpha(v))$. Clearly, $\gO_\phi$
satisfies A3 if and only if $\phi : \bbR \to \bbR$ is a polynomial.
\end{example}

\section{Natural operators and cohomology}
\label{2}

We start this section by a brief recollection of two classical
constructions.  For a Lie algebra $\gh$ and a $\gh$-module $W$, the
{\em Chevalley-Eilenberg\/} cohomology $H^*(\gh,W)$ of $\gh$ with
coefficients in $W$ is the cohomology of the cochain complex
$(C^*(\gh,W),\dCE)$ defined by $C^m (\gh,W) := \Lin(\ext^m \hskip
.2em \gh, W),\ m \geq 0$, with $\dCE$ the sum $\dCE = \delta_1 +
\delta_2$, where
\begin{eqnarray}
\label{d1}
(\delta_1 f)(\Rada h1{m+1}) &:=& \sum_{1 \leq i \leq m+1}
\zn {i+1} \cdot h_i f(h_1,\ldots,{\hat h}_i,\ldots,h_{m+1})\
\mbox { and }
\\
\label{d2}
(\delta_2 f)(\Rada h1{m+1})  &:=& \sum_{1 \leq i< j \leq m+1}
\zn {i+j} \cdot
f([h_i,h_j],h_1,\ldots,{\hat h}_i,\ldots,{\hat h}_j,\ldots,h_{m+1}),
\end{eqnarray}
for $f \in C^m (\gh,W)$, $\Rada h1{m+1} \in \gh$ and  $\ \hat{}\ $ denoting
the omission. If $m=0$, the summation in the right hand
side of~(\ref{d2}) runs over the empty set, so we put $(\delta_2 f)(h)
:= 0$ for $f \in C^0 (\gh,W)$.

The second notion we need to recall is the semidirect product of groups.
Assume that $G$ and $H$ are Lie groups, with $G$ acting on $H$ by
homomorphisms. 
One then defines the {\em semidirect product\/} $G \semidirect H$ as
the Cartesian product $G \times H$ with the multiplication
\[
(g_1,h_1)(g_2,h_2) := (g_1g_2,g_2^{-1}(h_1)h_2),\ g_1,g_2 \in G,\
h_1,h_2 \in H.
\]
Both $G$ and $H$ are subgroups of $G \semidirect H$
and their union $G \cup H$ generates $G \semidirect H$.
Let us close this introductory part by formulating a
proposition that ties the above two constructions together.

If $W$ is a left $\semGH$-module, the inclusion $H \subset G \semidirect H$
induces a left $H$-action on $W$ which in turn induces an
infinitesimal action of $\gh$ on $W$. 
One may therefore consider the
cochain complex $(C^*(\gh,W),\dCE)$. 
Since $G$ acts by homomorphisms, the unit
of $H$ is $G$-fixed, so there is an induced action of $G$ on the Lie
algebra $\gh$ of $H$. The group $G$ acts also on $W$, via the inclusion
$G \subset \semGH$. These two actions give rise, in the usual way, to
an action of $G$ on $C^*(\gh,W)$. Let us denote $C_G^*(\gh,W)$ the subspace of
$G$-fixed elements of $C^*(\gh,W)$. We have
the following:

\begin{proposition}
\label{Co_doma?}
The subspace of fixed elements  $C^*_G(\gh,W)\subset C^*(\gh,W)$ 
is $\dCE$-closed, so the cohomology
$H^*_G(\gh,W) := H^*(C^*_G(\gh,W),\dCE)$ is defined. For $H$ connected, there
is an  isomorphism
\begin{equation}
\label{Musim_zavolat_Jane_co_je_doma_pokud_se_neozve.}
H_G^0(\gh,W) \cong W^{\semGH},
\end{equation}
where $W^{\semGH}$ denotes, as usual, the space of $\semGH$-fixed 
elements in $W$.
\end{proposition}

\begin{proof}
We leave a direct verification of the $\dCE$-closeness of 
$C_G^*(\gh,W)$ as a simple exercise to the reader. It is equally
easy to see that
$H^0_G(\gh,W)$ consists of elements of $W$ which are simultaneously
$G$-fixed and $\gh$-invariant. If $H$ is connected, the exponential
map is an epimorphism, thus $\gh$-invariant elements in $W$ are
precisely those which are $H$-fixed. This, along with the fact that $G \cup H$
generates $\semGH$, 
gives~(\ref{Musim_zavolat_Jane_co_je_doma_pokud_se_neozve.}).  
\end{proof}

In Section~\ref{1} we recalled that natural
differential operators $\gO \in \Nat(\gF,\gG)$ 
between natural bundles of order $\leq l$ with fibers $\sF$ resp.\
$\sG$, correspond to $\GL
{k+l}n$-equivariant maps $O : \sF^{(k)} \to \sG$ with some $k \geq
0$. This can be expressed by the isomorphism:
\begin{equation}
\label{zub?}
\Nat(\gF,\gG) \cong \bigcup_{k \geq 0}
\Map_{\GL {k+l}n}(\sF^{(k)},\sG),
\end{equation}
where $\Map \hskip .2em \subb{\GL {k+l}n}(\sF^{(k)},\sG)$ 
is the space of {\em polynomial\/}
$\GL{k+l}n$-equivariant maps $\sF^{(k)} \to \sG$ -- see assumption~A3
on page \pageref{page}.
The space $\Map(\sF^{(k)},\sG)$ of all polynomial maps  
has the standard 
$\GL {k+l}n$-action induced from the actions on $\sF^{(k)}$ and
$\sG$. 

The space of equivariant maps is the fixed subspace $\Map \hskip .2em
\subb{\GL {k+l}n}(\sF^{(k)},\sG) = {\Map(\sF^{(k)},\sG)}^{\GL
{k+l}n}$.  Let us see how Proposition~\ref{Co_doma?} describes these
spaces.  The crucial observation is that $\GL sn$ is, for each $s \geq
1$, a semidirect product~\cite[Section~13]{kolar-michor-slovak}.  If
$\symexp {(\Rn)} r$ denotes the $r$th symmetric power of $\Rn$, $r
\geq 1$, then elements of $\GL sn$ are expressions $A = A_1 + A_2 +
A_3 + \cdots + A_{s}$, $A_i \in \Lin \left(\adj\symexp {(\Rn)}{i},
\Rn\right)$, $1 \leq i \leq s$, such that $A_1 : \Rn \to \Rn$ is
invertible.  It is a simple exercise to write formulas for the product
and inverse; for $s=2$ it was done in
Example~\ref{Co_se_deje_v_byte?}.

The space $\Lin \left(\symexp {(\Rn)}{i}, \Rn\right)$ is
canonically isomorphic to the space $\Sym\left(\otexp {(\Rn)}{i},
\Rn\right)$ of symmetric multilinear maps and we will identify these two
spaces in the sequel.  Denote by 
$
\NGL sn = \{ A = A_1 + A_2 + A_3 + \cdots + A_{s} \in \NGL sn;\
A_1 = {\it id\/}\}
$
the prounipotent radical of
$\GL sn$.
Then $\GL sn$ is the semidirect product
$\GL sn = \GLname_n \semidirect \NGL sn$,
with $\GLname_n$ acting on $\NGL sn$ by adjunction. Denote finally $\ngl sn$
the Lie algebra of $\NGL sn$,
\begin{equation}
\label{Andulka_by_taky_byla_pusinka_kdyby_vic_psala.}
\ngl sn = 
\{
a = a_2 + a_3 + \cdots + a_{s}; \hskip .5em 
a_i  \in \Sym\left(\otexp {(\Rn)}{i}, \Rn \right),\  2\leq i \leq s
\}.
\end{equation}

Assume that the action of $\GL ln$ on the fiber $\sG$ of $\gG$ is linear. 
Then $\Map(\sF^{(k)},\sG)$ is a linear representation of
$\GL {k+l}n$ and Proposition~\ref{Co_doma?} applied to
$G = \GLname_n$, $H = \NGL {k+l}n$ and $W = \Map(\sF^{(k)},\sG)$ gives 
\begin{equation}
\label{Boli_mne_pod_ramenem}
\Map_{\GL {k+l}n}(\sF^{(k)},\sG)
\cong H^0_{{\GLname_n}}\left(\adj\ngl {k+l}n,\Map(\sF^{(k)},\sG)\right).
\end{equation}

For each $k \geq 0$, the inclusion $\Map(\sF^{(k)},\sG)
\hookrightarrow
\Map(\sF^{(k+1)},\sG)$ together with the projection $\ngl {k+l+1}n \to \ngl
{k+l}n$ induces a $\GLname_n$-invariant inclusion
\[
C^*\left(\ngl{k+l}n, \Map(\sF^{(k)},\sG)\right) \
\hookrightarrow C^*\left(\adj\ngl{k+l +1}n, \Map(\sF^{(k+1)},\sG)\right)
\] 
which commutes with the differentials. Let us denote
\begin{equation}
\label{za_tyden_jedu_do_IHES}
C^*\left(\adj\nglname_n^\ii, \Map(\sF^\ii,\sG)\right) 
:= \bigcup_{k \geq 0}\bigcup_{l \geq 1}
C^*\left(\adj\ngl{k+l}n, \Map(\sF^{(k)},\sG)\right)
\end{equation}
and $C_{\GLname_n}^*(\nglname^\ii_n,
\Map(\sF^\ii,\sG))$ 
the $\GLname_n$-stable subspace of $C^*(\nglname_n^\ii,
\Map(\sF^\ii,\sG))$. Let
\[
H_{\GLname_n^\ii}^*\left(\adj\nglname_n^\ii, \Map(\sF^\ii,\sG)\right) 
:= H^*(C_{\GLname_n}^*\left(\adj\nglname_n^\ii, \Map(\sF^\ii,\sG)),\dCE\right).
\]
Then~(\ref{zub?}) together with~(\ref{Boli_mne_pod_ramenem}) and the
fact that cohomology commutes with direct limits implies:

\begin{theorem}
\label{Zitra_jedu_do_Hamburku.}
Let $\gF$ and $\gG$ be natural bundles with fibers $\sF$ resp.\
$\sG$ of orders $\leq l$. Suppose that the action of $\GL ln$ on $\sG$
is linear. Then, under the above notation
\begin{equation}
\label{Slavie_vyhrava!}
\Nat(\gF,\gG) \cong H_{\GLname_n}^0
\left(\adj\nglname_n^\ii,\Map(\sF^\ii,\sG)\right).
\end{equation}
\end{theorem}

In the following sections we show that, in many interesting cases, the 
cohomology in the right hand side of~(\ref{Slavie_vyhrava!}) is the
cohomology of a certain graph complex.

\section{Natural operators and graphs} 
\label{3}

We are going to describe natural differential operators by certain
spaces spanned by graphs. Roughly speaking, graphs, viewed as
contraction schemes for indices, will encode elementary
$\GLname_n$-invariant tensors in~(\ref{za_tyden_jedu_do_IHES}).
Our approach is based on a translation of the
Invariant Tensor Theorem into the graph language explained in~\cite{markl:ig}

Suppose that $\gB$ is a natural bundle satisfying A1 on
page~\pageref{page}, so that
the induced action of $\GLname_n
\subset \GL ln$ on the fiber $\sB$ is rational linear. 
According to standard facts of the representation theory of
$\GLname_n$ recalled, for instance,
in~\cite[\S~1.4]{katsylo-timashev}, 
an equivalent assumptions is
that, as a $\GLname_n$-module,
$\sB$ is the direct sum of $\GLname_n$-modules
\begin{equation}
\label{4}
\sB = \bigoplus_{1 \leq i \leq b} \sB_i,\
\end{equation}
where $\sB_i$ is, for each $1 \leq i \leq b$, either the space 
$\Lin(\otexp {\Rn}{q_i},\otexp {\Rn}{p_i})$ for some $p_i,q_i \geq 0$, with
the standard $\GLname_n$-action, or a subspace of this space
consisting of maps whose inputs and/or
outputs have a
specific symmetry, which can for example be expressed by a Young
diagram. 

In other words, $\sB_i$ are spaces of multilinear maps whose
coordinates are tensors $T^{\Rada a1{p_i}}_{\Rada b1{q_i}}$ with
$q_i$ input indices and $p_i$ output indices, which may or may not  
enjoy some kind of symmetry. We will graphically represent these 
tensors as corollas with $q_i$-inputs and $p_i$ outputs:
\begin{equation}
\raisebox{-3em}{\rule{0pt}{0pt}}
\label{88}
\unitlength 2.5mm
\linethickness{0.4pt}
\begin{picture}(20,3.7)(10.5,20.4)
\put(20,20){\vector(1,1){2}}
\put(20,20){\vector(-1,1){2}}
\put(20,20){\vector(-1,2){1}}
\put(18,18){\vector(1,1){1.9}}
\put(22,18){\vector(-1,1){1.9}}
\put(19,18){\vector(1,2){.935}}
\put(20,20){\makebox(0,0)[cc]{\Large$\bullet$}}
\put(20.5,18.1){\makebox(0,0)[cc]{$\ldots$}}
\put(20,16.5){\makebox(0,0)[cc]{%
   $\underbrace{\rule{10mm}{0mm}}_{\mbox{\scriptsize $q_i$ inputs}}$}}
\put(0,40){
\put(20.5,-18.25){\makebox(0,0)[cb]{$\ldots$}}
\put(20,-16.5){\makebox(0,0)[cc]{%
   $\overbrace{\rule{10mm}{0mm}}^{\mbox{\scriptsize $p_i$ outputs}}$}}}
\end{picture}\hskip -5em\raisebox{-.5em}{.}
\end{equation}
Instead of \raisebox{-.1em}{{\Large$\bullet$}} 
we may sometime use different
symbols for the node, such as $\nabla$, $\bbox$, {\Large$\circ$},~\&c.

\begin{example}
The fiber of the tangent bundle $T$ is $\Rn  = \Lin(\otexp \Rn0,
\otexp \Rn 1)$, so one has in~(\ref{4}) $b=1$, $p_1
= 1$, $q_0 = 0$. Elements of the fiber are tensors $X^a$
symbolized by \black.
The fiber $\sC$ of the connection bundle $\Con$ (see
Example~\ref{Co_se_deje_v_byte?}) is $\Lin(\otexp {\Rn}2, \otexp{\Rn}1)$,
therefore $b=1$, $p_1=1$ and $q_1 = 2$. Elements of $\sC$
are $\GLname_n$-tensors (Christoffel symbols) $\Gamma^a_{bc}$
represented by
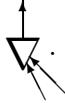
\[
\unitlength .5cm
\begin{picture}(.7,2)(-.5,-.5)
\put(-.5,0){\makebox(0,0)[cc]{\Large$\nabla$}}
\put(.3,0){\makebox(0,0)[cc]{.}}
\put(-.5,.4){\vector(0,1){1}}
\put(.1,-1.25){\vector(-1,2){.5}}
\put(.67,-1.1){\vector(-1,1){1}}
\end{picture}
\]

An example with a(n anti-)symmetry is the bundle $\Omega^m$ of de~Rham
$m$-forms, $m \geq 0$. Its fiber is the space $\Lin(\ext^m \Rn, \otexp
\Rn 0) = \Lin(\ext^m \Rn, \bbbR)$ of anti-symmetric tensors
$\omega_{)\Rada b1m(}$.

Sometimes we will need decorations of nodes. For example, the product bundle
$T \times T$ has fiber $\Rn \times \Rn$ generated by tensors $X^a,Y^a$
which will be denoted
\[
\unitlength .5cm
\begin{picture}(0,1)(0,.2)
\put(-.5,0){\makebox(0,0)[cc]{\Large$\bullet$}}
\put(0,0){\makebox(0,0)[bl]{\scriptsize $1$}}
\put(-.5,0){\vector(0,1){1.2}}
\end{picture}
\hskip 3em
\unitlength .5cm
\begin{picture}(0,1)(0,.2)
\put(-.5,0){\makebox(0,0)[cc]{\Large$\bullet$}}
\put(0,0){\makebox(0,0)[bl]{\scriptsize $2$}}
\put(-.5,0){\vector(0,1){1.2}}
\end{picture}
\hskip 2em 
\mbox { or }
\hskip 2em
\unitlength .5cm
\begin{picture}(0,1)(0,.2)
\put(-.5,0){\makebox(0,0)[cc]{\Large$\bullet$}}
\put(0,0){\makebox(0,0)[bl]{\scriptsize $X$}}
\put(-.5,0){\vector(0,1){1.2}}
\end{picture}
\hskip 3em
\unitlength .5cm
\begin{picture}(0,1)(0,.2)
\put(-.5,0){\makebox(0,0)[cc]{\Large$\bullet$}}
\put(0,0){\makebox(0,0)[bl]{\scriptsize $Y$}}
\put(-.5,0){\vector(0,1){1.2}} \hskip 1em .
\end{picture}
\]
\end{example}

Let $\gB$ be a natural bundle with fiber $\sB$ 
decomposed as in~(\ref{4}). It
is easy to see that the fiber $\sB^{(k)}$ of the $k$-jet bundle $\gB^{(k)}$
decomposes, as a $\GLname_n$-module, into
$\sB^{(k)} = \bigoplus_{1 \leq i \leq b}\sB^{(k)}_i$, where 
\begin{equation}
\label{preziju_Kryla?}
\sB^{(k)}_i = \bigoplus_{0 \leq v \leq k} 
\Sym(\otexp {\Rn} v,\bbR) \otimes \sB_i .
\end{equation}
This means that if elements of $\sB_i$ are tensors 
$T^{\Rada a1{p_i}}_{\Rada b1{q_i}}$, elements of $\sB^{(k)}_i$
are tensors
${}_{(\Rada s1{v})}T^{\Rada a1{p_i}}_{\Rada b1{q_i}}$,
$v \leq k$,
with braces indicating the symmetry in $(\Rada s1v)$. 
In terms of pictures this amounts to adding new symmetric inputs to
corollas~(\ref{88}), so elements
of $\sB^{(k)}_i$ will be symbolized by
\begin{equation}
\label{Zbynek_Hejda}
\raisebox{-3.5em}{\rule{0pt}{0pt}}
\unitlength 3mm
\linethickness{0.4pt}
\begin{picture}(20,3.5)(10.5,20.6)
\put(20,20){\makebox(0,0)[cc]{\Large$\bullet$}}
\put(26,20){\makebox(0,0)[cc]{,\ $v \leq k$.}}
\put(17.75,17.8){\makebox(0,0)[cc]{\scriptsize$( \hskip 12mm )$}}
\put(20,20){\vector(1,1){2}}
\put(20,20){\vector(-1,1){2}}
\put(20,20){\vector(-1,2){1}}
\put(20.5,18){\vector(-1,4){0.475}}
\put(21,18){\vector(-1,2){.95}}
\put(24,18){\vector(-2,1){3.8}}
\put(22.25,18){\makebox(0,0)[cc]{$\ldots$}}
\put(22.3,16.7){\makebox(0,0)[cc]{%
   $\underbrace{\rule{10mm}{0mm}}_{\mbox{\scriptsize $q_i$ inputs}}$}}
\put(19.5,18){\vector(1,4){0.475}}
\put(19,18){\vector(1,2){.95}}
\put(16,18){\vector(2,1){3.8}}
\put(17.75,18){\makebox(0,0)[cc]{$\ldots$}}
\put(17.7,16.7){\makebox(0,0)[cc]{%
   $\underbrace{\rule{10mm}{0mm}}_{\mbox{\scriptsize $v$ inputs}}$}}
\put(0,40){
\put(20.5,-18.25){\makebox(0,0)[cc]{$\ldots$}}
\put(20,-16.7){\makebox(0,0)[cc]{%
   $\overbrace{\rule{12mm}{0mm}}^{\mbox{\scriptsize $p_i$ outputs}}$}}}
\end{picture}
\end{equation}

\begin{example}
\label{To_jsem_zvedav_co_vsechno_mi_Sergej_provede.}
The fiber of the $k$th tangent bundle $T^{(k)}$ is the space of tensors
\begin{equation}
\label{zz}
X^a_{(s_1,\ldots,s_v)}:= \frac {\pa^u X^a}{\pa x^{s_1} \cdots \pa
x^{s_v}},\  v \leq k,
\end{equation}
which we draw as
\begin{equation}
\label{tyden_piti}
\raisebox{-2.2em}{\rule{0pt}{0pt}}
\unitlength 3mm
\linethickness{0.4pt}
\begin{picture}(20,2.5)(10.5,19.8)
\put(20,20){\vector(0,1){2}}
\put(18,18){\vector(1,1){1.9}}
\put(22,18){\vector(-1,1){1.9}}
\put(19,18){\vector(1,2){.935}}
\put(20,20){\makebox(0,0)[cc]{\Large$\bullet$}}
\put(25,20){\makebox(0,0)[cc]{,\ $v \leq k$.}}
\put(20,18){\makebox(0,0)[cc]{\scriptsize (\hskip 13.5mm)}}
\put(20.5,18){\makebox(0,0)[cc]{$\ldots$}}
\put(20,16.7){\makebox(0,0)[cc]{%
   $\underbrace{\rule{12mm}{0mm}}_{\mbox{\scriptsize $v$ inputs}}$}}
\end{picture}
\end{equation}
The fiber of the bundle $\Con^{(k)}$ is the space of tensors
${}_{(s_1,\ldots,s_v)}\Gamma^a_{bc}:= 
\frac {\pa^u \Gamma^a_{bc}}{\pa x^{s_1} \cdots \pa
x^{s_v}}$, $v \leq k$,
depicted as
\begin{equation}
\label{Uz_nevim_co.}
\unitlength .4cm
\begin{picture}(.7,2)(-.5,-.5)
\put(-.5,0){\makebox(0,0)[cc]{\Large$\nabla$}}
\put(2.5,0){\makebox(0,0)[cc]{, $v \leq k$.}}
\put(-.5,.4){\vector(0,1){1}}
\put(.1,-1.25){\vector(-1,2){.5}}
\put(.67,-1.1){\vector(-1,1){1}}
\unitlength 3mm
\put(-21,-20){
\put(19.4,18){\vector(1,4){0.475}}
\put(18.9,18){\vector(1,2){.95}}
\put(16,18){\vector(2,1){3.8}}
\put(17.75,18){\makebox(0,0)[cc]{$\ldots$}}
\put(17.7,16.6){\makebox(0,0)[cc]{%
   $\underbrace{\rule{10mm}{0mm}}_{\mbox{\scriptsize $v$ inputs}}$}}
\put(17.7,17.8){\makebox(0,0)[cc]{\scriptsize $( \hskip 12mm )$}}
}
\end{picture}
\raisebox{-2.5em}{\rule{0pt}{0pt}}
\end{equation}
\end{example}

As follows from~(\ref{Andulka_by_taky_byla_pusinka_kdyby_vic_psala.}),
$\ngl {k+l}n = \bigoplus_{2 \leq u \leq k+l}\Sym(\otexp {\Rn}u,\Rn)$.
Therefore $\ngl {k+l}n$ is the space of symmetric tensors
$\phi^b_{(\Rada s1u)}$, $2 \leq u \leq k+l$, or in pictures,
\begin{equation}
\label{O_byly_po_vsich_muziky/a_choraly_nam_hraly_temne}
\raisebox{-1.2em}{\rule{0pt}{0pt}}
\unitlength 3mm
\linethickness{0.4pt}
\begin{picture}(20,3)(10.5,18.3)
\put(20,20.3){\vector(0,1){1.7}}
\put(18,18){\vector(1,1){1.8}}
\put(22,18){\vector(-1,1){1.8}}
\put(19,18){\vector(1,2){.91}}
\put(20,20){\makebox(0,0)[cc]{\Large$\circ$}}
\put(27,20){\makebox(0,0)[cc]{,\ $2 \leq u \leq k+l$.}}
\put(20,18){\makebox(0,0)[cc]{\scriptsize (\hskip 13.5mm)}}
\put(20.5,18.25){\makebox(0,0)[cc]{$\ldots$}}
\put(20,17){\makebox(0,0)[cc]{%
   $\underbrace{\rule{12mm}{0mm}}_{\mbox{\scriptsize $u$ inputs}}$}}
\end{picture}
\end{equation}
In what follows, white
corollas~(\ref{O_byly_po_vsich_muziky/a_choraly_nam_hraly_temne}) 
will {\em always denote\/}
elements of $\ngl {k+l}n$ for some $k+l \geq 2$.

In the rest of this section we construct a graded space 
$\GrFG^*$ spanned by graphs representing $\GLname_n$-invariant
cochains  in $C^*_{\GLname_n}(\nglname_n^\ii,\Map(\sF^\ii,\sG))$. 
The differentials will be studied in the next section.

Suppose that the natural bundles $\gF$ and $\gG$ satisfy assumption A1
on page~\pageref{page}, and see what can be said about the space
$C^m_{\GLname_n}(\adj\ngl {k+l}n, \Map(\sF^{(k)},\sG))$ of
$\GLname_n$-equivariant polynomial maps from $\Map(\ngl {k+l}n \times
\sF^{(k)},\sG)$ that are $m$-homogeneous and antisymmetric in $\ngl
{k+l}n$.  By the polynomiality assumption~A3,
\begin{equation}
\label{Ten_clovek_bude_moje_smrt.}
C^m_{\GLname_n}\left(\adj\ngl {k+l}n, \Map(\sF^{(k)},\sG)\right) \cong
\bigoplus_{t \geq 0}\Lin_{\GLname_n}
\left(\ext^m \ngl {k+l}n \ot \otexp{{\sF^{(k)}}}t,\sG
\right),
\end{equation}
where $\Lin_{\GLname_n}(-,-)$ denotes the space of
${\GLname_n}$-equivariant linear
maps.

Let us decompose the fibers $\sF$ and $\sG$ of natural bundles 
$\gF$ and $\gG$ into the direct sum~(\ref{4}), 
$\sF = \bigoplus_{1 \leq i \leq f} \sF_i$ and
$\sG = \bigoplus_{1 \leq i \leq g} \sG_i$.
By~(\ref{preziju_Kryla?}), the components of the fiber 
$\sF^{(k)}$ of the $k$-jet bundle $\gF^{(k)}$, $k \geq 0$, are
the direct sums
$\sF^{(k)}_i = \bigoplus_{0 \leq v \leq k} \sF^{[v]}_i$, $1 \leq i \leq f$,
with $\sF^{[v]}_i := \Sym(\otexp {{\Rn}}v,\Rn) \ot \sF_i$.
Using the above decompositions and
description~(\ref{Andulka_by_taky_byla_pusinka_kdyby_vic_psala.}) of
$\ngl {k+l}n$, one can rewrite the right hand side
of~(\ref{Ten_clovek_bude_moje_smrt.}) into
\begin{equation}
\label{V_utery_se_mu_asi_nevyhnu.}
\bigoplus_{t \geq 0}\bigoplus_{S(k,l,t)}\Lin_{\GLname_n}
\left(
\Ext_{1 \leq i \leq m}  \Sym(\otexp {{\Rn}}{{u_i}},\Rn)
\ot
\bigotimes_{1 \leq s \leq t} \sF^{[v_s]}_{i_s}, \sG_i
\right),
\end{equation}
where $S(k,l,t)$ is the set of integers $\Rada u1m$, $\rada{i_1}{i_t}$,
$\rada{v_1}{v_t}$ and $i$ such that
\[
2 \leq \Rada u1m \leq k+l,\ 1 \leq \rada{i_1}{i_t} \leq f,\
0 \leq \rada{v_1}{v_t} \leq k\ \mbox { and } 1 \leq i \leq g.
\]

Let us fix a multiindex $\omega= (\Rada
u1m,\rada{i_1}{i_t},\rada{v_1}{v_t},i) \in S(k,l,t)$. By our
assumptions, the space $\sF^{[v_s]}_{i_s}$ is, for each $1 \leq s \leq
t$, isomorphic to the space 
\[
\Lin_{\Si_s}^{\So_s}(\otexp
{{\Rn}}{(v_s + q_{i_s})},\otexp {{\Rn}}{p_{i_s}})
:= 
\{ f \in  \Lin(\otexp{{\Rn}}{(v_s + q_{i_s})},\otexp {{\Rn}}{p_{i_s}})
;\ f {\mathfrak s} = 0 = {\mathfrak t} f \mbox {
for } {\mathfrak s} \in \Si_s,  {\mathfrak t} \in \So_s
\}.
\]
of linear maps having a symmetry
specified by subsets $\Si_s \subset \bfk[\Sigma_{v_s + q_{i_s}}]$, 
$\So_s \subset
\bfk[\Sigma_{p_{i_s}}]$, see also~\cite[Remark~4.4]{markl:ig}.  
Similarly, $\sG_i
\cong \Lin_{\Si}^{\So}(\otexp {{\Rn}}{c}, \otexp {{\Rn}}{d})$, for
some $c,d \geq 0$ and subsets $\Si \subset \bfk[\Sigma_{c}]$, $\So \subset
\bfk[\Sigma_{d}]$. The expression
\begin{equation}
\label{Jestli_preziju_utery+preziju_vsechno}
\Lin_{\GLname_n}\left(
\Ext_{1 \leq i \leq m}  \Sym(\otexp {{\Rn}}{{u_i}},\Rn)
\ot
\bigotimes_{1 \leq s \leq t} \sF^{[v_s]}_{i_s}, \sG_i
\right)
\end{equation}
in~(\ref{V_utery_se_mu_asi_nevyhnu.}) is therefore isomorphic to
\begin{equation}
\label{Piji_kafe}
\Lin_{\GLname_n}\left(
\Ext_{1 \leq i \leq m}  \Sym(\otexp {{\Rn}}{{u_i}},\Rn)
\ot
\bigotimes_{1 \leq s \leq t} 
\Lin_{\Si_s}^{\So_s}(\otexp {{\Rn}}{{(v_s + q_{i_s})}},
\otexp {{\Rn}}{p_{i_s}}),
\Lin_{\Si}^{\So}(\otexp{{\Rn}}{c}, \otexp{{\Rn}}{d})
\right).
\end{equation}
Let us remark that in all applications discussed in this paper, we
will always have $p_{i_s} = 1$ for $1 \leq s \leq t$, $c=0$ and $d=1$. 

Observe that~(\ref{Piji_kafe}) is the space
in~(24) of~\cite{markl:ig}, with an
appropriate choice of the parameters, which in this case is $r := t-m$,
and
\begin{eqnarray*}
h_i := u_i && \mbox { for } 1 \leq i \leq m, \mbox  { and }
\\
h_i := v_s + q_{i_s},\ \Si_i := \Si_s,\ \So_i := \So_s &&
\mbox { for } i = s+m,\ 1 \leq s \leq t,
\end{eqnarray*}
therefore the methods developed in~\cite{markl:ig} apply. We believe that
the reader can tolerate a certain incompatibility between the notation
used in this paper and the notation of~\cite{markl:ig} -- the alphabet
does not have enough letters to avoid notational conflicts.

By Proposition~4.8 and Remark~4.10 
of~\cite{markl:ig}, the space~(\ref{Piji_kafe})
is related to the space $\Gr^m_\omega$ spanned by graphs with vertices
of three types:
\[
\def\arraystretch{1.2}
\begin{array}{rl}
\mbox{\bf 1st type:} & \mbox{%
$t$ `black' vertices~(\ref{Zbynek_Hejda}) with $p_i := p_{i_s}$, 
$q_i:= q_{i_s}$ and $v := v_s$, representing
tensors}
\\
& \mbox  {in $\sF_{j_s}^{[v_s]}$, $1 \leq s \leq t$,}
\\
\mbox{\bf 2nd type:} & \mbox{%
one vertex~(\ref{88}) with $p_i : =c$ and $q_i : =d$ called the {\em
  anchor\/}, representing tensors}
\\
&\mbox  {in the dual $\sG_i^*$ of $\sG_i$, and}
\\
\mbox{\bf 3rd type:} & \mbox{%
$m$ `white'
  vertices~(\ref{O_byly_po_vsich_muziky/a_choraly_nam_hraly_temne})
  with $u := u_i$
representing generators of the Lie algebra}
\\
& \mbox{$\ngl {k+l}n$, $1 \leq i \leq m$.}     
\end{array}
\]
Our graphs are directed and oriented, where
an {\em orientation\/} is, by definition, 
an equivalence class of linear orders of the set of white vertices,
modulo the relation identifying orders that differ by an even
number of transpositions. If the orientations of two graphs $G'$ and
$G''$  differ by an odd number of
transpositions, we put $G' = -G''$ in $\Gr^m_\omega$. This
notion of 
orientation is not the traditional one 
but resembles orientations in various graph 
complexes~\cite[\S~II.5.5]{markl-shnider-stasheff:book}.

The graphs spanning $\Gr^m_\omega$ are not
required to be connected, and multiple edges and loops are allowed .
The vertices above are Merkulov's {\em genes\/}~\cite{merkulov:PROP}.
The unique vertex of the 2nd type marks the place where we evaluate
the composition along the graph at an element of $\sG^*$, which
explains the dualization in the definition of this vertex.

Proposition~4.8 of~\cite{markl:ig} (or its obvious extension mentioned 
in \cite[Remark~4.10]{markl:ig}), 
combined with the isomorphism
between~(\ref{Jestli_preziju_utery+preziju_vsechno})
and~(\ref{Piji_kafe}), gives an epimorphism
\begin{equation}
\label{Bude_tam_vecer_ten_misa?}
R_{n,\omega}^m : \Gr^m_\omega \epi
\Lin_{\GLname_n}\left(
\Ext_{1 \leq i \leq m}  \Sym(\otexp {{\Rn}}{{u_i}},\Rn)
\ot
\bigotimes_{1 \leq s \leq t} \sF^{[v_s]}_{i_s}, \sG_i
\right)
\end{equation}
which is, by \cite[Proposition~4.9]{markl:ig}, 
a monomorphism if $n+m \geq$ the number of edges of graphs in
$\Gr^m_\omega$.
The central result of this section,
Theorem~\ref{Zase_mne_boli_v_krku_ale_to_neni_muj_hlavni_problem.}
below, uses the limit
\begin{equation}
\label{stavka}
\Gr^m_{\gF,\gG} := 
\bigcup_{k \geq 0} \bigcup_{l \geq 1} \bigoplus_{t \geq 0}
\bigoplus_{\omega \in S(k,l,t)} \Gr^m_\omega.
\end{equation}
The space $\Gr^m_{\gF,\gG}$ is
spanned by graphs with an arbitrary number of the 1st
type vertices with an arbitrary $v \geq 0$ in~(\ref{Zbynek_Hejda}),
one 2nd type vertex representing tensors in $\sG^*_i$ for $1 \leq i
\leq g$, and $m$ 3rd type vertices with an arbitrary $u \geq 2$
in~(\ref{O_byly_po_vsich_muziky/a_choraly_nam_hraly_temne}).

\begin{theorem}
\label{Zase_mne_boli_v_krku_ale_to_neni_muj_hlavni_problem.}
The epimorphisms $R^m_{n,\omega}$ in~(\ref{Bude_tam_vecer_ten_misa?}) assemble,
for each $m \geq 0$, into a surjection
\begin{equation}
\label{Ticho_pred_bouri?}
R_n^m : \GrFG^m \epi
C^m_{\GLname_n}\left(\nglname_n^\ii,\Map(\sF^\ii,\sG)\right).
\end{equation}
The restriction
\[ 
R^m_n(e) : \GrFG^m(e)  
\to C^m_{\GLname_n}\left(\nglname_n^\ii,\Map(\sF^\ii,\sG)\right)
\]
of the map $R^m_n$ to the subspace $\GrFG^m(e) \subset \GrFG^m$ 
spanned by graphs with $\leq e$ edges, is a monomorphism whenever 
$n = \dim (M) \geq e -m$.
\end{theorem}

\begin{proof}
The maps $R^m_{n,\omega}$ of~(\ref{Bude_tam_vecer_ten_misa?})
assemble, for each $k \geq 0$ and $l \geq 1$, into an epimorphism 
\[
R^m_{n,k,l} := 
\bigoplus_{t \geq 0} \bigoplus_{\omega \in S(k,l,t)} R^m_{n,\omega} 
: \bigoplus_{t \geq 0} \bigoplus_{\omega \in S(k,l,t)}
\Gr^m_\omega \epi 
\bigoplus_{t \geq 0}
\Lin_{\GLname_n}
\left(\ext^m \ngl {k+l}n \ot \otexp{{\sF^{(k)}}}t,\sG
\right).
\]
Recalling~(\ref{za_tyden_jedu_do_IHES}),~(\ref{Ten_clovek_bude_moje_smrt.}),
and the definition~(\ref{stavka}) of the graph complex $\Gr^m_{\gF,\gG}$, we
conclude that $R_n^m := \bigoplus_{k \geq 0}\bigoplus_{l \geq 1} 
R^m_{n,k,l}$  is the desired surjection~(\ref{Ticho_pred_bouri?}).
The second part of the theorem follows from 
\cite[Proposition~4.9]{markl:ig} applied to 
the constituents $R_{n,\omega}^m$ of $R^m_n$.
\end{proof}
 
\begin{example}
\label{sileny_smutek_Zbynka_Hejdy}
Let us discuss the case $\gF =T \times T$  and $\gG = T$, where $T$ is the
tangent bundle functor. Graphs spanning the 
vector space $\Gr^m_{T\times T,T}$ have finite number of the 1st type
vertices~(\ref{tyden_piti})
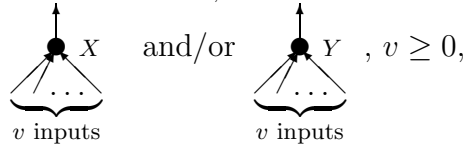
\[
\raisebox{-1.5em}{\rule{0pt}{0pt}}
\hskip -1.5cm
\unitlength 3mm
\linethickness{0.4pt}
\begin{picture}(20,3.3)(10.5,18)
\put(20,20){\vector(0,1){2}}
\put(18,18){\vector(1,1){1.9}}
\put(22,18){\vector(-1,1){1.9}}
\put(19,18){\vector(1,2){.935}}
\put(20,20){\makebox(0,0)[cc]{\Large$\bullet$}}
\put(21,20){\makebox(0,0)[lc]{\scriptsize$X$}}
\put(20.5,18){\makebox(0,0)[cc]{$\ldots$}}
\put(20,16.7){\makebox(0,0)[cc]{%
   $\underbrace{\rule{12mm}{0mm}}_{\mbox{\scriptsize $v$ inputs}}$}}
\end{picture}
\hskip -2cm \raisebox{.50cm}{\mbox {and/or}}\hskip -2cm
\unitlength 3mm
\linethickness{0.4pt}
\begin{picture}(20,3.3)(10.5,18)
\put(20,20){\vector(0,1){2}}
\put(18,18){\vector(1,1){1.9}}
\put(22,18){\vector(-1,1){1.9}}
\put(19,18){\vector(1,2){.935}}
\put(20,20){\makebox(0,0)[cc]{\Large$\bullet$}}
\put(21,20){\makebox(0,0)[lc]{\scriptsize$Y$}}
\put(20.5,18){\makebox(0,0)[cc]{$\ldots$}}
\put(20,16.7){\makebox(0,0)[cc]{%
   $\underbrace{\rule{12mm}{0mm}}_{\mbox{\scriptsize $v$ inputs}}$}}
\end{picture}
\hskip -2.3cm \raisebox{.5cm}{\mbox {, $v \geq 0,$}}
\]
marking the places where to insert tensors
$X^a_{(s_1,\ldots,s_v)}$  and $Y^a_{(s_1,\ldots,s_v)}$ of 
the fiber of $(T \times T)^{(\infty)}$. The unique vertex 
\anchor\ of the 2nd type is the
place to insert a tensor of
the fiber $\Rn^*$ of $T^*$.
There of course will also be 
$m$ vertices~(\ref{O_byly_po_vsich_muziky/a_choraly_nam_hraly_temne})
of the 3rd type for 
generators of $\ngl {\infty}n$.

Observe that we omitted braces
indicating the symmetry because inputs of all vertices are
symmetric and no confusion may occur.
Let us inspect how $\Gr_{T\times T,T}^0$ describes
$\GLname_n$-equivariant maps in
$Map_{\GLname_n}(\adj(\Rn \times \Rn)^{(\infty)},\Rn) = 
C^0_{\GLname_n}(\ngl \infty n,
\Map((\Rn \times \Rn)^{(\infty)},\Rn))$.
The graph
\[
\unitlength .5cm
\begin{picture}(0,1.7)(0,0.4)
\put(0,2){\makebox(0,0)[cc]{\hskip .5mm$\bbox$}}
\put(0,1){\vector(0,1){.935}}
\put(0,1){\makebox(0,0)[cc]{\Large$\bullet$}}
\put(.4,1){\makebox(0,0)[lc]{\scriptsize$Y$}}
\put(0,0){\vector(0,1){.935}}
\put(0,0){\makebox(0,0)[cc]{\Large$\bullet$}}
\put(.4,0){\makebox(0,0)[lc]{\scriptsize$X$}}
\end{picture}
\]
describes the equivariant map 
that sends an element $(X^a,X^a_b,Y^a,Y^a_b) \in (\Rn \times \Rn)^{(1)}$ 
into the element $(X^j Y^a_j) \in \Rn$. It is precisely the map
$O_2$ considered in Example~\ref{Cesti}. The linear combination
\[
\unitlength .5cm
\begin{picture}(0,2)(0,0.4)
\put(0,2){\makebox(0,0)[cc]{\hskip .5mm$\bbox$}}
\put(0,1){\vector(0,1){.935}}
\put(0,1){\makebox(0,0)[cc]{\Large$\bullet$}}
\put(.4,1){\makebox(0,0)[lc]{\scriptsize$Y$}}
\put(0,0){\vector(0,1){.935}}
\put(0,0){\makebox(0,0)[cc]{\Large$\bullet$}}
\put(.4,0){\makebox(0,0)[lc]{\scriptsize$X$}}
\end{picture}
\hskip 1cm \raisebox{.2cm}{-}\hskip .8cm
\unitlength .5cm
\begin{picture}(0,2)(0,0.4)
\put(0,2){\makebox(0,0)[cc]{\hskip .5mm$\bbox$}}
\put(0,1){\vector(0,1){.935}}
\put(0,1){\makebox(0,0)[cc]{\Large$\bullet$}}
\put(.4,1){\makebox(0,0)[lc]{\scriptsize$X$}}
\put(0,0){\vector(0,1){.935}}
\put(0,0){\makebox(0,0)[cc]{\Large$\bullet$}}
\put(.4,0){\makebox(0,0)[lc]{\scriptsize$Y$}}
\end{picture}
\raisebox{-1em}{\rule{0pt}{0pt}}
\]
represents the local formula 
$(X^a,X^a_b,Y^a,Y^a_b) \mapsto (X^j Y^a_j - Y^j X^a_j)$
for the Lie bracket $[X,Y]$ of two vector fields.
We allow also graphs as
\[
\unitlength .5cm
\begin{picture}(0,1.6)(1.1,1)
\put(0,2){\makebox(0,0)[cc]{\hskip .5mm$\bbox$}}
\put(0,1){\vector(0,1){.935}}
\put(0,1){\makebox(0,0)[cc]{\Large$\bullet$}}
\put(.4,1){\makebox(0,0)[lc]{\scriptsize$X$}}
\put(0,.2){
\put(2,1){\makebox(0,0)[cc]{\circle {1.5}}}
\put(1.3,1.35){\makebox(0,0)[cc]{\Large$\bullet$}}
\put(1.32,1.25){\makebox(0,0)[tc]{\vector(0,1){0}}}
\put(1.1,1.55){\makebox(0,0)[rb]{\scriptsize $Y$}}
\put(3.1,1.35){\makebox(0,0)[cc]{,}}
}
\end{picture}
\]
which represents the map $(X^a,X^a_b,Y^a,Y^a_b) \mapsto (X^a Y^i_i)$
involving the trace $Y^i_i$ of $Y$.
An example of a degree $1$ cochain in 
$C^1_{\GLname_n}(\ngl 2n,\Map(\Rn \times \Rn,\Rn))$ is provided by
\[
\unitlength .5cm
\begin{picture}(0,2.3)(0,0)
\put(0,2){\makebox(0,0)[cc]{\hskip .5mm$\bbox$}}
\put(0,1.22){\vector(0,1){.65}}
\put(0,1){\makebox(0,0)[cc]{\Large$\circ$}}
\put(-1,-1){
\put(0,1){{\vector(1,1){.92}}}
\put(0,1){\makebox(0,0)[cc]{\Large$\bullet$}}
\put(.4,1){\makebox(0,0)[lc]{\scriptsize$X$}}
}
\put(1,-1){
\put(0,1){{\vector(-1,1){.92}}}
\put(0,1){\makebox(0,0)[cc]{\Large$\bullet$}}
\put(.4,1){\makebox(0,0)[lc]{\scriptsize$Y$}}
}
\put(2,1){\makebox(0,0)[cc]{,}}
\end{picture}
\]
which defines the $\GLname_n$-equivariant 1-cochain
$(\phi^a_{bc},X^a,Y^a) \mapsto (\phi^a_{ij}X^iY^j)$.

As explained in~\cite[Remark~5.2]{markl:ig}, 
for degrees $\geq 2$ our interpretation of
graphs involves the antisymmetrization in white vertices. For instance,
the graph
\[
\raisebox{-2em}{\rule{0pt}{0pt}}
\unitlength .4cm
\begin{picture}(4,2.2)(-2,1)
\put(0,0){\oval(2,2)[b]}
\put(0,2){\oval(2,2)[t]}
\put(-1,0){\line(0,1){2}}
\put(1,0){\makebox(0,0)[b]{\Large$\circ$}}
\put(0.6,0){\makebox(0,0)[br]{\scriptsize $2$}}
\put(1,2){\makebox(0,0)[t]{\Large$\circ$}}
\put(0.6,1.6){\makebox(0,0)[br]{\scriptsize $1$}}
\put(2,-1){\makebox(0,0)[b]{\Large$\bullet$}}
\put(2.4,-1){\makebox(0,0)[lb]{\scriptsize $Y$}}
\put(0,-.5){
\put(2,1){\makebox(0,0)[b]{\Large$\bullet$}}
\put(2.4,1){\makebox(0,0)[lb]{\scriptsize $X$}}
\put(2.15,1.15){\vector(-1,1){.9}}
}
\put(1,.6){\vector(0,1){.85}}
\put(1,0.1){\vector(0,1){0}}
\put(2.15,-.8){\vector(-1,1){.95}}
\end{picture}
\]
represents the $2$-cochain $(\phi^a_{bc},\psi^a_{bc},X^a,Y^a) \mapsto
(\phi^i_{jk}\psi^j_{il} - \psi^i_{jk}\phi^j_{il})X^kY^l$.  The reason why
the expected traditional $\frac 1{2!}$\ -factor is missing
is explained in Remark~\ref{tlak}.
\end{example}

\begin{example}
In this example we express local formulas for the
covariant derivative, torsion and curvature in terms of
graphs. The covariant derivative is the 
operator $\nabla : \Con \times T^{\times 2} \to T$ locally given by the graph
\begin{equation}
\label{covar}
\raisebox{-2em}{\rule{0pt}{0pt}}
\unitlength .4cm
\begin{picture}(.7,2)(6,-.2)
\put(-.5,0){\makebox(0,0)[cc]{\Large$\nabla$}}
\put(2,0){\makebox(0,0)[c]{$+$}}
\put(-3,0){\makebox(0,0)[r]{$\nabla_X Y:$}}
\put(-.5,.4){\vector(0,1){1.1}}
\put(.1,-1.25){\vector(-1,2){.5}}
\put(.67,-1.1){\vector(-1,1){1}}
\put(-.45,1.4){\makebox(0,0)[bc]{$\bbox$}}
\put(0.2,-1.5){\makebox(0,0)[cc]{\Large$\bullet$}}
\put(0.2,-2){\makebox(0,0)[tc]{\scriptsize $X$}}
\put(1.5,-1.7){\makebox(0,0)[tc]{\scriptsize $Y$}}
\put(0.85,-1.3){\makebox(0,0)[cc]{\Large$\bullet$}}
\put(4,.3){
\unitlength .5cm
\begin{picture}(0,2.5)(0,1)
\put(0,2){\makebox(0,0)[cc]{\hskip .5mm$\bbox$}}
\put(0,1){\vector(0,1){.935}}
\put(0,1){\makebox(0,0)[cc]{\Large$\bullet$}}
\put(.4,1){\makebox(0,0)[lc]{\scriptsize$Y$}}
\put(0,0){\vector(0,1){.935}}
\put(0,0){\makebox(0,0)[cc]{\Large$\bullet$}}
\put(.4,0){\makebox(0,0)[lc]{\scriptsize$X$}}
\end{picture}
}
\end{picture}\hskip -.3cm,
\end{equation}
which is a graphical form of formula~(\ref{Sergej_je_magor.}).
The torsion $T: \Con \times T^{\times 2} \to T$ is given by
\[
\raisebox{-1.8em}{\rule{0pt}{0pt}}
\unitlength .4cm
\begin{picture}(.7,2.3)(6,-.5)
\put(-.5,0){\makebox(0,0)[cc]{\Large$\nabla$}}
\put(2,0.3){\makebox(0,0)[c]{$-$}}
\put(-3,0){\makebox(0,0)[r]{$T(X,Y):$}}
\put(-.5,.4){\vector(0,1){1.1}}
\put(.1,-1.25){\vector(-1,2){.5}}
\put(.67,-1.1){\vector(-1,1){1}}
\put(-.45,1.4){\makebox(0,0)[bc]{$\bbox$}}
\put(0.2,-1.5){\makebox(0,0)[cc]{\Large$\bullet$}}
\put(0.2,-2){\makebox(0,0)[tc]{\scriptsize $X$}}
\put(1.5,-1.7){\makebox(0,0)[tc]{\scriptsize $Y$}}
\put(0.85,-1.3){\makebox(0,0)[cc]{\Large$\bullet$}}
\put(5.3,-.5){
\unitlength .4cm
\begin{picture}(.7,2)(1,-.5)
\put(-.5,0){\makebox(0,0)[cc]{\Large$\nabla$}}
\put(-.5,.4){\vector(0,1){1.1}}
\put(.1,-1.25){\vector(-1,2){.5}}
\put(.67,-1.1){\vector(-1,1){1}}
\put(-.45,1.4){\makebox(0,0)[bc]{$\bbox$}}
\put(0.2,-1.5){\makebox(0,0)[cc]{\Large$\bullet$}}
\put(0.2,-2){\makebox(0,0)[tc]{\scriptsize $Y$}}
\put(1.5,-1.7){\makebox(0,0)[tc]{\scriptsize $X$}}
\put(0.85,-1.3){\makebox(0,0)[cc]{\Large$\bullet$}}
\end{picture}
}
\end{picture}
\]
and the curvature $R: \Con \times T^{\times 3} \to T$ as
\[
\raisebox{-3.2em}{\rule{0pt}{0pt}}
\unitlength .4cm
\begin{picture}(.7,3)(6,-.5)
\put(-.5,0){\makebox(0,0)[cc]{\Large$\nabla$}}
\put(2,0){\makebox(0,0)[c]{$-$}}
\put(-3,0){\makebox(0,0)[r]{$R(X,Y)Z:$}}
\put(-.5,.4){\vector(0,1){1.1}}
\put(.1,-1.25){\vector(-1,2){.5}}
\put(.67,-1.1){\vector(-1,1){1}}
\put(-1.3,-1.12){\vector(2,3){.65}}
\put(-.45,1.4){\makebox(0,0)[bc]{$\bbox$}}
\put(0.2,-1.5){\makebox(0,0)[cc]{\Large$\bullet$}}
\put(0.2,-2){\makebox(0,0)[tc]{\scriptsize $X$}}
\put(-1.7,-1.9){\makebox(0,0)[tc]{\scriptsize $Y$}}
\put(1.5,-1.7){\makebox(0,0)[tc]{\scriptsize $Z$}}
\put(0.85,-1.3){\makebox(0,0)[cc]{\Large$\bullet$}}
\put(-1.4,-1.3){\makebox(0,0)[cc]{\Large$\bullet$}}
\end{picture}
\put(-.5,0){
\unitlength .4cm
\begin{picture}(.7,3)(1,-.5)
\put(-.5,0){\makebox(0,0)[cc]{\Large$\nabla$}}
\put(2,0){\makebox(0,0)[c]{$+$}}
\put(6.8,0){\makebox(0,0)[c]{$-$}}
\put(10.8,0){\makebox(0,0)[c]{.}}
\put(-.5,.4){\vector(0,1){1.1}}
\put(.1,-1.25){\vector(-1,2){.5}}
\put(.67,-1.1){\vector(-1,1){1}}
\put(-1.3,-1.12){\vector(2,3){.65}}
\put(-.45,1.4){\makebox(0,0)[bc]{$\bbox$}}
\put(0.2,-1.5){\makebox(0,0)[cc]{\Large$\bullet$}}
\put(0.2,-2){\makebox(0,0)[tc]{\scriptsize $Y$}}
\put(-1.7,-1.9){\makebox(0,0)[tc]{\scriptsize $X$}}
\put(1.5,-1.7){\makebox(0,0)[tc]{\scriptsize $Z$}}
\put(0.85,-1.3){\makebox(0,0)[cc]{\Large$\bullet$}}
\put(-1.4,-1.3){\makebox(0,0)[cc]{\Large$\bullet$}}
\put(5.5,-.5){
\begin{picture}(.7,3)(1,-.5)
\put(-.5,0){\makebox(0,0)[cc]{\Large$\nabla$}}
\put(-.5,.4){\vector(0,1){1.1}}
\put(.1,-1.25){\vector(-1,2){.5}}
\put(1.1,-.9){\vector(-3,2){1.3}}
\put(-.45,1.4){\makebox(0,0)[bc]{$\bbox$}}
\put(0.2,-1.5){\makebox(0,0)[cc]{\Large$\bullet$}}
\put(0.2,-2){\makebox(0,0)[tc]{\scriptsize $Y$}}
\end{picture}
\put(.6,-1.4){
\unitlength .4cm
\begin{picture}(.7,3)(1,-.5)
\put(-.5,0){\makebox(0,0)[cc]{\Large$\nabla$}}
\put(.1,-1.25){\vector(-1,2){.5}}
\put(.67,-1.1){\vector(-1,1){1}}
\put(0.2,-1.5){\makebox(0,0)[cc]{\Large$\bullet$}}
\put(0.2,-2){\makebox(0,0)[tc]{\scriptsize $X$}}
\put(1.5,-1.7){\makebox(0,0)[tc]{\scriptsize $Z$}}
\put(0.85,-1.3){\makebox(0,0)[cc]{\Large$\bullet$}}
\end{picture}
}}
\put(10.2,-.5){
\begin{picture}(.7,3)(1,-.5)
\put(-.5,0){\makebox(0,0)[cc]{\Large$\nabla$}}
\put(-.5,.4){\vector(0,1){1.1}}
\put(.1,-1.25){\vector(-1,2){.5}}
\put(1.1,-.9){\vector(-3,2){1.3}}
\put(-.45,1.4){\makebox(0,0)[bc]{$\bbox$}}
\put(0.2,-1.5){\makebox(0,0)[cc]{\Large$\bullet$}}
\put(0.2,-2){\makebox(0,0)[tc]{\scriptsize $X$}}
\end{picture}
\put(0.6,-1.4){
\unitlength .4cm
\begin{picture}(.7,3)(1,-.5)
\put(-.5,0){\makebox(0,0)[cc]{\Large$\nabla$}}
\put(.1,-1.25){\vector(-1,2){.5}}
\put(.67,-1.1){\vector(-1,1){1}}
\put(0.2,-1.5){\makebox(0,0)[cc]{\Large$\bullet$}}
\put(0.2,-2){\makebox(0,0)[tc]{\scriptsize $Y$}}
\put(1.5,-1.7){\makebox(0,0)[tc]{\scriptsize $Z$}}
\put(0.85,-1.3){\makebox(0,0)[cc]{\Large$\bullet$}}
\end{picture}
}}
\end{picture}
}
\]
\end{example}

\begin{example}
This example shows that the map $R^m_n$ from 
Theorem~\ref{Zase_mne_boli_v_krku_ale_to_neni_muj_hlavni_problem.}
need not be a monomorphism below the `stable range.'
Consider again the two graphs from
Example~\ref{sileny_smutek_Zbynka_Hejdy}:
\[
\unitlength .5cm
\begin{picture}(0,2.5)(0,0)
\put(0,2){\makebox(0,0)[cc]{\hskip .5mm$\bbox$}}
\put(0,1){\vector(0,1){.935}}
\put(0,1){\makebox(0,0)[cc]{\Large$\bullet$}}
\put(.4,1){\makebox(0,0)[lc]{\scriptsize$Y$}}
\put(0,0){\vector(0,1){.935}}
\put(0,0){\makebox(0,0)[cc]{\Large$\bullet$}}
\put(.4,0){\makebox(0,0)[lc]{\scriptsize$X$}}
\put(-2,1){\makebox(0,0){$G_1 := $}}
\put(1.8,1.1){\makebox(0,0){and}}
\put(2.8,1){\makebox(0,0)[lc]{$G_2 := $}}
\put(8.7,1){\makebox(0,0)[lc]{.}}
\end{picture}
\hskip 3.5cm
\unitlength .5cm
\begin{picture}(0,2.5)(1.1,.5)
\put(0,2){\makebox(0,0)[cc]{\hskip .5mm$\bbox$}}
\put(0,1){\vector(0,1){.935}}
\put(0,1){\makebox(0,0)[cc]{\Large$\bullet$}}
\put(.4,1){\makebox(0,0)[lc]{\scriptsize$X$}}
\put(0,.2){
\put(2,1){\makebox(0,0)[cc]{\circle {1.5}}}
\put(1.3,1.35){\makebox(0,0)[cc]{\Large$\bullet$}}
\put(1.32,1.25){\makebox(0,0)[tc]{\vector(0,1){0}}}
\put(1.1,1.55){\makebox(0,0)[rb]{\scriptsize $Y$}}
}
\end{picture}
\]
The number of edges of both graphs is $2$.
As we already saw, $G_1$ represents the local formula
$\sum_{1 \leq i,j \leq n} X^j \, {\pa Y^i}/{\pa x^j} \,
\pa/\pa x^i$
and $G_2$ the formula
$\sum_{1 \leq i,j \leq n} {\pa Y^j}/{\pa x^j}\, X^i {\pa}/{\pa x^i}$.
For $n=1$ both formulas give the same result, namely 
$X\, {\pa Y}/{\pa x}\,  {\pa }/{\pa x}$,
therefore $R^0_1(G_1) = R^0_1(G_2)$.
For $n \geq 2$ one clearly has $R^0_n(G_1) \not= R^0_n(G_2)$.
\end{example}

\section{The differential}
\label{sec4}

In this section we express the restriction of the 
Chevalley-Eilenberg differential onto the subcomplex $
C^*_{\GLname_n}(\adj\nglname_n^\ii,\Map(\sF^\ii,\sG))$ of
$\GLname_n$-equivariant cochains in
terms of graph complexes.
Let us describe first the bracket in the limit $\nglname_n^\ii =
\bigcup_{s \geq 2}\ngl sn$ of Lie algebras $\ngl sn$ recalled
in~(\ref{Andulka_by_taky_byla_pusinka_kdyby_vic_psala.}).  If finite
sums $a = a_2 + a_3 + a_4 + \cdots$ and $b = b_2 + b_3 + b_4 + \cdots$
are elements of $\nglname_n^\ii$, $a_u,b_u \in \Sym(\otexp {(\Rn)}{u},
\Rn)$, $u \geq 2$, then $[a,b] = [a,b]_3 + [a,b]_4 + \cdots$ (no
quadratic term) with
\[
[a,b]_u = 
\sum_{s+t =u+1}\sum_{1 \leq i \leq s} 
\left(\adj S(a_s \circ_i b_t) -  S(b_s \circ_i a_t)\right),
\]
where $S(-)$ denotes the symmetrization (see Remark~\ref{tlak}) of a
linear map $\otexp {\Rn}{u} \to \Rn$, $a_s\circ_i b_t$ is the
insertion of $b_t$ into the $i$th slot of $a_s$ and $b_s\circ_i a_t$
has the similar obvious meaning. For $\Rada v1u \in \Rn$ we easily get
\begin{eqnarray}
\label{modified_bracket}
[a,b]_u(\Rada v1u) =
\sum_{s+t =u+1}&& \hskip -1.5em \frac{s!t!}{u!}
\sum_\sigma  \left\{ a_s(b_t(\Rada v{\sigma(1)}{\sigma(t)}),
\Rada v{\sigma(t+1)}{\sigma(u)})
- \right.
\\
\nonumber 
&&\hskip 3em -
\left.
b_s(a_t(\Rada v{\sigma(1)}{\sigma(t)}),
\Rada v{\sigma(t+1)}{\sigma(u)})\right\},
\end{eqnarray}
where $\sigma$ runs over all $(t,s-1)$-unshuffles $\sigma$,
i.e.~permutations $\sigma \in \Sigma_u$ such that 
$\sigma(1) < \ldots < \sigma(t)$, $\sigma(t+1) < \ldots < \sigma(u)$.

\begin{remark}
\label{maminka}
In the rest of the paper, we will consider $\nglname_n^\ii$ with the
{\em modified\/} Lie bracket, given by
formula~(\ref{modified_bracket}) {\em without\/} the
$\frac{s!t!}{u!}$\ -coefficients. Since this modified Lie
algebra is isomorphic to the original one, via the isomorphism $a_s
\mapsto s!  \cdot a_s$, for $a_s \in \Sym\left(\adj\otexp {(\Rn)}{s},
\Rn\right)$, $s \geq 2$, our modification is purely conventional.  The
advantage of this modified bracket is that the corresponding replacement
rule~(\ref{Co_Sergej_planuje?}) is a linear combination of graphs
without fractional coefficients.
\end{remark}

To help the reader to appreciate the idea of the differential, we
start with an informal definition. A precise formula including signs
and orientations is given in~(\ref{Porad_mi_nejak_lidi_nepisou.}).
At the beginning of Section~\ref{2}  we decomposed 
the CE-differential into the sum
$\dCE = \delta_1 + \delta_2$. Let us analyze the action of the
second piece $\delta_2$ first. 
A graph $G$ representing a
$\GLname_n$-invariant $m$-cochain has $m$ white vertices
that mark the places where to insert elements of
$\nglname_n^\ii$. Let us label, for $m \geq 1$, these white 
vertices by $\ell \in \{\rada 1m\}$ and denote the vertex
labelled $\ell$ by~$w_\ell$. If $m=0$, there are no white vertices and
no labelling is necessary.

The effect of the differential $\delta_2$ on the graph $G$ is, by the
definition recalled in~(\ref{d2}), the following.  For each $\ell \in
\{\rada 1m\}$ insert to the vertex $w_\ell$ the element
$[h_i,h_j]$ and to the remaining white vertices elements
$h_1,\ldots,{\hat h}_i,\ldots,{\hat h}_j,\ldots,h_{m+1}$, make the
summation over all $1 \leq i < j \leq m+1$ and antisymmetrize in
$\Rada h1{m+1}$. Denote the resulting $(m+1)$-cochain by
$G_\ell$. Then $\delta_2(G) = \varepsilon_1 \cdot G_1 + \cdots +
\varepsilon_1 \cdot G_m$, where $\Rada {\varepsilon}1m \in \{-1,+1\}$
are appropriate signs.  A moment's reflection reveals that $G_\ell$ is
obtained by replacing the vertex $w_\ell$ by:
\begin{equation}
\label{Co_Sergej_planuje?}
\raisebox{-3.7em}{\rule{0pt}{0pt}}
\unitlength 3mm
\linethickness{0.4pt}
\begin{picture}(20,3)(10.5,18.8)
\put(20,20.3){\vector(0,1){1.7}}
\put(18,18){\vector(1,1){1.8}}
\put(22,18){\vector(-1,1){1.8}}
\put(19,18){\vector(1,2){.91}}
\put(20,20){\makebox(0,0)[cc]{\Large$\circ$}}
\put(20.5,18){\makebox(0,0)[cc]{$\ldots$}}
\put(20,16.7){\makebox(0,0)[cc]{%
   $\underbrace{\rule{12mm}{0mm}}_{\mbox{\scriptsize $u$ inputs}}$}}
\end{picture}
\hskip -2.5cm \longmapsto \hskip .5cm \sum_{t+s = u + 1}
\hskip -1.2cm
\unitlength 3mm
\linethickness{0.4pt}
\begin{picture}(20,3)(10.5,18.8)
\put(20,20.3){\vector(0,1){1.7}}
\put(18,18){\vector(1,1){1.8}}
\put(22,18){\vector(-1,1){1.8}}
\put(19,18){\vector(1,2){.91}}
\put(20,20){\makebox(0,0)[cc]{\Large$\circ$}}
\put(20.5,18){\makebox(0,0)[cc]{$\ldots$}}
\put(20.27,16.8){\makebox(0,0)[cc]{%
   $\underbrace{\rule{10mm}{0mm}}_{\mbox{\scriptsize $s$}}$}}
\end{picture}
\put(-22.5,-2.2){
\begin{picture}(20,4)(10.5,18.8)
\put(18,18){\vector(1,1){1.8}}
\put(22,18){\vector(-1,1){1.8}}
\put(19,18){\vector(1,2){.91}}
\put(20,20){\makebox(0,0)[cc]{\Large$\circ$}}
\put(20.5,18){\makebox(0,0)[cc]{$\ldots$}}
\put(20,16.7){\makebox(0,0)[cc]{%
   $\underbrace{\rule{12mm}{0mm}}_{\mbox{\scriptsize $t$}}$}}
\put(21.7,19){\makebox(0,0)[cc]{$\left(\rule{0pt}{17pt} \hskip 2.1cm
    \right)_{\rm ush}$}}
\end{picture}
}\hskip -2cm, \hskip 1cm
\end{equation}
where the braces $(-)_{\rm ush}$ indicate that the
summation over all $(t,s-1)$-unshuffles of the inputs has been performed.
This is precisely the formula for the generators
of the {\em homological vector field\/} introduced by
Merkulov~\cite{merkulov:defquant,merkulov:PROP}.  One also
recognizes~(\ref{Co_Sergej_planuje?}) as the graphical representation 
of the axioms of 
$L_\infty$-algebras as given in~\cite[page~160]{markl:shalg}.

A similar analysis shows that $\delta_1$ acts by replacing each vertex
of type 1 or 2 by the pictorial representation of the action of $\ngl
\infty n$ on tensors corresponding to this vertex. We will show
instances of these `pictorial presentations' 
in the following two examples.

\begin{example}
\label{polar}
Consider a symmetric map $\xi : \otexp {\Rn}{v} \to \Rn$
representing an element in the fiber of the $k$-th tangent space $T^{(k)}$ 
with coordinates $X^a_{(s_1,\ldots,s_v)}$ (see~(\ref{zz}) of 
Example~\ref{To_jsem_zvedav_co_vsechno_mi_Sergej_provede.}). The
action of $a = a_2 + a_3 + a_4 + \cdots \in \ngl \infty n$
on $\xi$ is given by $a \xi = (a\xi)_{u+1} + (a\xi)_{u+2} + \cdots$, where
\[
(a\xi)_v = 
\sum_{s+u=v+1} 
\left(\sum_{1 \leq i \leq s} S(a_s \circ_i \xi) -  
\sum_{1 \leq i \leq v }S(\xi \circ_i a_s)\right).
\]
Removing fractional coefficients by modifying the $\ngl \infty
n$-action (compare Remark~\ref{maminka}), one can
graphically express the above rule by the following polarization
of~(\ref{Co_Sergej_planuje?}): 
\begin{equation}
\label{Nikdo}
\raisebox{-3.2em}{\rule{0pt}{0pt}}
\unitlength 3mm
\hskip -2cm
\linethickness{0.4pt}
\begin{picture}(20,3)(10.5,18.8)
\put(20,20.3){\vector(0,1){1.7}}
\put(18,18){\vector(1,1){1.8}}
\put(22,18){\vector(-1,1){1.8}}
\put(19,18){\vector(1,2){.91}}
\put(20,20){\makebox(0,0)[cc]{\Large$\bullet$}}
\put(19.5,20.5){\makebox(0,0)[rb]{\scriptsize $X$}}
\put(20.5,18){\makebox(0,0)[cc]{$\ldots$}}
\put(20,16.7){\makebox(0,0)[cc]{%
   $\underbrace{\rule{12mm}{0mm}}_{\mbox{\scriptsize $v$ inputs}}$}}
\end{picture}
\hskip -2.2cm \longmapsto \sum_{s+u=v+1}
\hskip .7cm
\unitlength 3mm
\linethickness{0.4pt}
\begin{picture}(15,3)(15.5,18.8)
\put(20,20.3){\vector(0,1){1.7}}
\put(18,18){\vector(1,1){1.8}}
\put(22,18){\vector(-1,1){1.8}}
\put(19.1,18){\vector(1,2){.91}}
\put(20,20){\makebox(0,0)[cc]{\Large$\circ$}}
\put(20.5,18){\makebox(0,0)[cc]{$\ldots$}}
\put(20.4,16.8){\makebox(0,0)[cc]{%
   $\underbrace{\rule{10mm}{0mm}}_{\mbox{\scriptsize $s$}}$}}
\end{picture}
\put(-22.5,-2.2){
\begin{picture}(10,4)(10.5,18.8)
\put(18,18){\vector(1,1){1.8}}
\put(22,18){\vector(-1,1){1.8}}
\put(19.1,18){\vector(1,2){.91}}
\put(20,20){\makebox(0,0)[cc]{\Large$\bullet$}}
\put(19.5,20.5){\makebox(0,0)[rb]{\scriptsize $X$}}
\put(20.5,18){\makebox(0,0)[cc]{$\ldots$}}
\put(20,16.8){\makebox(0,0)[cc]{%
   $\underbrace{\rule{11mm}{0mm}}_{\mbox{\scriptsize $u$}}$}}
\put(21.7,19){\makebox(0,0)[cc]{$\left(\rule{0pt}{17pt} \hskip 2.1cm
\right)_{\rm ush}$}}
\end{picture}
}
\hskip -2cm
-
\hskip -1cm
\unitlength 3mm
\linethickness{0.4pt}
\begin{picture}(20,3)(10.5,18.8)
\put(20,20.3){\vector(0,1){1.7}}
\put(18,18){\vector(1,1){1.8}}
\put(22,18){\vector(-1,1){1.8}}
\put(19.1,18){\vector(1,2){.91}}
\put(20,20){\makebox(0,0)[cc]{\Large$\bullet$}}
\put(19.5,20.5){\makebox(0,0)[rb]{\scriptsize $X$}}
\put(20.5,18){\makebox(0,0)[cc]{$\ldots$}}
\put(20.4,16.8){\makebox(0,0)[cc]{%
   $\underbrace{\rule{10.7mm}{0mm}}_{\mbox{\scriptsize $u$}}$}}
\end{picture}
\put(-22.5,-2.2){
\begin{picture}(20,3)(10.5,18.8)
\put(18,18){\vector(1,1){1.8}}
\put(22,18){\vector(-1,1){1.8}}
\put(19.1,18){\vector(1,2){.91}}
\put(20,20){\makebox(0,0)[cc]{\Large$\circ$}}
\put(20.5,18){\makebox(0,0)[cc]{$\ldots$}}
\put(20,16.7){\makebox(0,0)[cc]{%
   $\underbrace{\rule{12mm}{0mm}}_{\mbox{\scriptsize $s$}}$}}
\put(21.7,19){\makebox(0,0)[cc]{$\left(\rule{0pt}{20pt} \hskip 2.1cm 
\right)_{\rm ush}$}}
\end{picture}
}\hskip -2cm  .
\end{equation}
\end{example}

\begin{example}
Let us write two initial replacement rules for the connection
and its derivatives. The first one is the infinitesimal version 
of~(\ref{Napise_mi_nekdo?}):
\[
\unitlength .4cm
\begin{picture}(.7,1.5)(-.5,-.5)
\put(-.5,0){\makebox(0,0)[cc]{\Large$\nabla$}}
\put(-.5,.4){\vector(0,1){1}}
\put(.1,-1.25){\vector(-1,2){.5}}
\put(.67,-1.1){\vector(-1,1){1}}
\end{picture} \hskip .5cm
\longmapsto \hskip .3cm -
\unitlength .23cm
\begin{picture}(4,3)(-2,-1)
\put(0,0.3){\vector(0,1){1.7}}
\put(-1,-2){\vector(1,2){.89}}
\put(1,-2){\vector(-1,2){.89}}
\put(0,0){\makebox(0,0)[cc]{\Large$\circ$}}
\end{picture}.
\]
The next one is a graphical form of an equation that can be
found in~\cite[Section~17.7]{kolar-michor-slovak} (but notice a different 
convention for covariant derivatives 
used in~\cite{kolar-michor-slovak}):
\[
\unitlength .5cm
\begin{picture}(.7,1)(-.5,-.5)
\put(-.5,0){\makebox(0,0)[cc]{\Large$\nabla$}}
\put(-.5,.45){\vector(0,1){1}}
\put(.1,-1.25){\vector(-1,2){.5}}
\put(.67,-1.1){\vector(-1,1){1}}
\put(-1.7,-1.1){\vector(1,1){1}}
\end{picture} \hskip .5cm
\longmapsto \hskip .3cm
\unitlength 2.2mm
\begin{picture}(4,3)(-2,-2)
\put(0,0.5){\vector(0,1){1.7}}
\put(-2,-4){\vector(1,2){1.8}}
\put(1.9,-1.9){\vector(-1,1){1.6}}
\put(0,0){\makebox(0,0)[cc]{\Large$\circ$}}
\end{picture}
\unitlength .5cm
\begin{picture}(.7,1)(-.5,0.4)
\put(-.5,0){\makebox(0,0)[cc]{\Large$\nabla$}}
\put(.1,-.9){\vector(-1,2){.39}}
\put(.57,-.72){\vector(-1,1){.8}}
\end{picture} 
\hskip .5cm - \hskip 1cm
\unitlength .35cm
\begin{picture}(.7,1)(-.5,-1)
\put(-.5,0){\makebox(0,0)[cc]{\Large$\nabla$}}
\put(-.5,.6){\vector(0,1){1}}
\put(-.25,-.9){\vector(0,1){.7}}
\put(1.2,-2){\vector(-2,3){1.35}}
\end{picture}
\unitlength .3cm
\begin{picture}(4,3)(0.54,.34)
\put(-1.5,-1.5){\vector(1,1){1.3}}
\put(.9,-1.5){\vector(-1,2){.65}}
\put(0,0){\makebox(0,0)[cc]{\Large$\circ$}}
\end{picture}
\hskip -0.1cm - \hskip 1cm
\unitlength .4cm
\begin{picture}(.7,1)(-.5,-1)
\put(-.5,0){\makebox(0,0)[cc]{\Large$\nabla$}}
\put(-.5,.4){\vector(0,1){1}}
\put(.45,-1.9){\vector(-1,2){.83}}
\put(.78,-.9){\vector(-1,1){.9}}
\end{picture}
\unitlength .3cm
\begin{picture}(4,3)(-.9,.25)
\put(-2.5,-1.5){\vector(3,2){2.22}}
\put(1.5,-1.5){\vector(-1,1){1.3}}
\put(0,0){\makebox(0,0)[cc]{\Large$\circ$}}
\end{picture}
\hskip -.1cm - \hskip .5cm
\unitlength .3cm
\begin{picture}(4,3)(-2,-1)
\put(0,0.3){\vector(0,1){1.7}}
\put(-2,-2){\vector(1,1){1.8}}
\put(2,-2){\vector(-1,1){1.8}}
\put(0,-2){\vector(0,1){1.75}}
\put(0,0){\makebox(0,0)[cc]{\Large$\circ$}}
\end{picture}.
\raisebox{-1.3em}{\rule {0pt}{0pt}}
\]
We are not going to give a general formula. For our
purposes, it will be enough to know that it is of the form
\begin{equation}
\label{Chce_se_mi_curat.}
\unitlength .4cm
\begin{picture}(.7,1)(-.5,-.5)
\put(-.5,0){\makebox(0,0)[cc]{\Large$\nabla$}}
\put(-.5,.4){\vector(0,1){1}}
\put(.1,-1.25){\vector(-1,2){.48}}
\put(.67,-1.1){\vector(-1,1){.95}}
\unitlength 2.2mm
\put(-21,-20){
\put(19.4,18){\vector(1,4){0.42}}
\put(18.82,18){\vector(1,2){.85}}
\put(16,18){\vector(2,1){3.6}}
\put(17.55,18){\makebox(0,0)[cc]{$\ldots$}}
\put(17.7,16){\makebox(0,0)[cc]{%
   $\underbrace{\rule{7mm}{0mm}}_{\mbox{\scriptsize $w$ inputs}}$}}
}
\end{picture}
\hskip .5cm \longmapsto \hskip .5cm
G_w \hskip .2cm - \hskip -.2cm
\unitlength 2.2mm
\begin{picture}(20,3)(15.5,18.8)
\put(20,20.6){\vector(0,1){1.7}}
\put(18,18){\vector(1,1){1.8}}
\put(22,18){\vector(-1,1){1.8}}
\put(19,18){\vector(1,2){.88}}
\put(20,20){\makebox(0,0)[cc]{\Large$\circ$}}
\put(20.5,18){\makebox(0,0)[cc]{$\ldots$}}
\put(20,16){\makebox(0,0)[cc]{%
   $\underbrace{\rule{9 mm}{0mm}}_{\mbox{\scriptsize $w+2$}}$}}
\end{picture}\hskip -2.9cm,
\raisebox{-2.1em}{\rule{0pt}{0pt}}
\end{equation}
where $G_w$ is a linear combination of $2$-vertex trees with one
vertex~(\ref{Uz_nevim_co.}), with $v < w$, and one 
vertex~(\ref{O_byly_po_vsich_muziky/a_choraly_nam_hraly_temne}) with $u < w+2$.
\end{example}

Let us write a formal definition of the graph
differential. For each oriented graph $G \in \GrFG^m$ we
define $\delta(G) \in \GrFG^{m+1}$ as the sum over the set $\Vert(G)$
of vertices of $G$,
\begin{equation}
\label{Porad_mi_nejak_lidi_nepisou.}
\delta(G) = \sum_{v \in \Vert(G)}\varepsilon_v \cdot \delta_v(G),
\end{equation}
where $\delta_v$ is the replacement of the vertex $v$ determined by
the type of $v$ and geometric data  as explained above. The signs
$\varepsilon_v$ and the orientations of the graphs in $\delta_v(G)$ are
determined in the following way.

(i) The operation $\delta_v$ replaces a  1st or 2nd type vertex $v$ 
by a linear combination of graphs containing precisely one white vertex. The
orientation of the graphs in  
$\delta_v(G)$ is given by the unique linear order such
that this new white vertex is the minimal element and the relative order of the
remaining white vertices is unchanged. The sign $\varepsilon_v$ is
$+1$. Symbolically
\begin{equation}
\label{B1}
\delta_v(\circ < \cdots < \circ) =
+1 \cdot (\delta_v(\bullet) < \circ < \cdots < \circ).
\end{equation} 

(ii) Let $v$ be a white vertex. We may assume that, after changing the sign of
the graph $G$ if necessary, $v$ is the minimal element in an
order determining the orientation. 
The orientation of graphs in $\delta_v(G)$ is then given by requiring
that the lower left white vertex in the right hand side  
of~(\ref{Co_Sergej_planuje?}) is the
minimal one, the upper right white vertex of~(\ref{Co_Sergej_planuje?}) is
the next one, and that the relative order of the remaining white
vertices is unchanged. The sign $\varepsilon_v$ is
again $+1$. Symbolically, 
\[
\delta_v(\circ < \cdots < \circ) =
+1 \cdot (\delta_v(\circ) < \circ < \cdots < \circ).
\]

We leave as a simple exercise to derive from the rule~(ii) that, if
the white vertex $v$ is the $i$th element of a linear order
determining the orientation of $G$, for some $1 \leq i \leq m$, the
orientations of graphs in $\delta_v(G)$ are symbolically expressed as
\begin{equation}
\label{B2}
\delta_v(\circ < \cdots < \circ) =
(-1)^{i+1} \cdot (
\underbrace{\circ < \cdots < \circ}_{\mbox{\scriptsize $i-1$}} <
\delta_v(\circ) <    
\underbrace{\circ < \cdots < \circ}_{\mbox{\scriptsize $m-i$}}).
\end{equation}

Let us emphasize that the applications in this paper use only the initial
part $\delta : \GrFG^0 \to \GrFG^1$ of the differential.
Since the graphs spanning $\GrFG^0$ (resp.~$\GrFG^1$) 
have no white vertices (resp.~one white vertex), the
orientation issue is trivial and all $\varepsilon_v$'s
in~(\ref{Porad_mi_nejak_lidi_nepisou.}) are $+1$.

\begin{theorem}
\label{Nikdo_se_neozyva.}
The object $\GrFG^* = (\GrFG^*,\delta)$ is a cochain complex and the maps
$R_n^m$ in~(\ref{Ticho_pred_bouri?}) assemble into a cochain map
$R_n^* : (\GrFG^*,\delta) \to
(C^*_{\GLname_n}(\nglname_n^\ii,
\Map(\sF^\ii,\sG)),\dCE)$.
\end{theorem}

\begin{proof}
Using the antisymmetry of $f$, one can rewrite 
equations~(\ref{d1}) and~(\ref{d2}) into
\begin{eqnarray*}
(\delta_1 f)(\Rada h1{m+1}) &=& 
\frac1{m!}\ \Ant\left(\adj h_1 f(\Rada h2{m+1})\right)\
\mbox { and }
\\
(\delta_2 f)(\Rada h1{m+1})  &=& 
\frac 1{ 2(m-1)!}\ \Ant
\left(\adj f([h_2,h_1],\Rada h3{m+1})\right),
\end{eqnarray*}
where $\Ant(-)$ denotes the antisymmetrization, see Remark~\ref{tlak}
below. If the multilinear map $f$ itself is
an antisymmetrization $\Ant(F)$ of a map $F$, one can rewrite
the above displays into
\begin{eqnarray}
\label{A1}
(\delta_1 f)(\Rada h1{m+1}) &=& \Ant\left(\adj h_1 F(\Rada h2{m+1})\right)\
\mbox { and }
\\
\label{A2}
(\delta_2 f)(\Rada h1{m+1})  &=& 
\Ant \left(\adj\right.
\sum_{1 \leq i \leq m}\zn {i+1} F(\adj\Rada h1{i-1}, [h_{i+1},h_{i}],
\Rada h{i+2}{m+1})
\left. \adj\right).
\end{eqnarray}

After this preparation, we prove that $R^*_n$ is a chain map by
verifying that $(R^{m+1}_n \circ \delta) (G) = 
(\dCE \circ R^m_n)(G)$ for each graph $G$ generating  $\GrFG^m$. After
choosing a linear order of white vertices of $G$ compatible with its
orientation, an appropriate version of 
the `state sum'~(11) of~\cite{markl:ig} gives a
multilinear map $F$ such that $R^m_n(G) = \Ant (F)$, 
see~\cite[Remark~5.2]{markl:ig}.

It is not difficult to see that $R^*_n$ translates the part of the
differential $\delta(G)$ in~(\ref{Porad_mi_nejak_lidi_nepisou.}) given
by the summation over the 1st and 2nd type vertices into
formula~(\ref{A1}) for $\delta_1(f)$ and the part of $\delta(G)$ given by
the summation over the white vertices to formula~(\ref{A2}) for
$\delta_2(f)$.  This fact is also reflected by the obvious similarity
between formulas~(\ref{A1}) and~(\ref{A2}) for the Chevalley-Eilenberg
differential and symbolic formulas~(\ref{B1}) and~(\ref{B2}) for the
graph differential.

The condition $\delta^2=0$  
can be verified directly using the fact that the local
replacement rules used in~(\ref{Porad_mi_nejak_lidi_nepisou.})
are duals of Lie algebra actions and checking
that the orientations were defined in such a way that the signs
combine properly. One may, however, proceed also as follows.

Since both the domain and target of
the map $R^*_n$, as well as $R^*_n$ itself, are defined in terms of
``standard representations,''  $R^m_n$
makes sense for an arbitrary natural $n$.
Let $G \in \GrFG^m$. By the finitary nature of objects involved,
there exists $e \geq 0$ such that all graphs that constitute
$\delta^2 (G) \in \GrFG^{m+2}$ have $\leq e$ edges. 
Choose $n \geq e - m -2$. We already
know that $R^*_n$ commutes with the differentials, therefore
$R^{m+2}_n( \delta^2 (G) ) = \dCE^2 (R^m_n (G)) = 0$.
By the second part of 
Theorem~\ref{Zase_mne_boli_v_krku_ale_to_neni_muj_hlavni_problem.}
this implies that $\delta^2 (G) = 0$.  
\end{proof}

\begin{remark}
\label{tlak}
In this paper, the {\em antisymmetrization\/} of an element $x$ of some
(say) right $\Sigma_k$-module, $k \geq 1$, is given by the formula
$\Ant(x) := \sum_{\sigma \in \Sigma_k} {\rm sgn}(\sigma)\cdot x \sigma$,
{\em without\/} the traditional~$\frac 1{k!}$. 
This convention is forced by the standard definition
of the Lie algebra associated to an associative algebra $(A,\cdot)$ --
the bracket $[a',a''] := a'\cdot a'' - a'' \cdot a'$,
$a',a'' \in A$,
does not involve the  $\frac1{2!}$-factor. On the
other hand, we define the {\em symmetrization\/} of $x$ as above by
the expected formula 
$S(x) := \frac 1{k!} \sum_{\sigma \in \Sigma_k}  x \sigma$.
\end{remark}

\begin{remark}
\label{Snad_mam_hodnoceni_AAA}
Applications of our theory will often be based on a 
suitable choice of a subspace of
$\Nat(\gF,\gG)$, together with the corresponding subcomplex of the
graph complex $\GrFG^*$.  These subobjects, denoted  for the
purposes of this remark by $\uNat(\gF,\gG)$ and $\uGrFG^* =
(\uGrFG^*,\udelta)$, will be chosen so that the number of
edges of graphs spanning $\uGrFG^m$ will be, for each $m \geq
0$, bounded by $C +m$, where $C$ is a fixed constant.

An example is the subcomplex $\Gr_{\bullet}^*(d)$ of the graph complex
$\Gr^*_{T^{\times d},T}$, introduced in Section~\ref{5}, that
describes $d$-multilinear operators on vector fields. Graphs spanning
$\Gr_{\bullet}^*(d)$ have precisely $d+m$ edges, so $C=d$ for this
subcomplex.  Another example is the subcomplex
$\Gr_{\bullet\nabla}^*(d)$ of $\Gr^*_{\Con \times T^{\times d},T}$
describing `connected' $d$-multilinear operators used in
Section~\ref{7}. Each degree $m$ graph spanning this subcomplex has
at most $2d+m-1$ edges, i.e.~$C= 2d-1$ in this case. The third example
is the complex $\Gr^*_{\bullet\nabla\wc}(d)$ introduced on
page~\pageref{AAA} describing `connected' operators in $\Nat(\Con
\times T^{\ot d},\Re)$. For this complex, $C := 2d$.

Let $(\uGrFG^*,\udelta)$, $\uNat(\gF,\gG)$ and the constant $C$ be as
above. By
Theorem~\ref{Zitra_jedu_do_Hamburku.} combined with
Theorem~\ref{Nikdo_se_neozyva.}, the restriction $\uR^*_n$
of $R^*_n$ induces the map
\begin{equation}
\label{uuu}
H^0(\uR^*_n) : H^0(\uGrFG^*,\udelta) \to \uNat(\gF,\gG)
\end{equation}
which is an isomorphism in {\em stable dimensions\/}. By this we mean that
the dimension $n$ of the underlying manifold $M$ is $\geq C$. 
If this happens, then the map $\uR^*_n$ is, by
\cite[Proposition~4.9]{markl:ig}, a chain isomorphism, so 
$H^0(\uR^*_m)$ is an isomorphism, too. If the dimension of $M$ is
less than the stable dimension, one cannot say anything
about the induced map $H^0(\uR^*_m)$,  although 
the chain map $\uR^*_n$ is still a chain epimorphism.
\end{remark}

\begin{example}
\label{Ted_mam_kasel.}
In this example we prove a baby version of Theorem~\ref{main-A}.
Namely, we show that the only natural bilinear operations on
vector fields on manifolds of dimensions $\geq 2$ are scalar multiples
of the Lie bracket. 
It will be convenient to have ready some 
initial cases of formula~(\ref{Nikdo}) for 
the replacement rule of vertices representing vector fields and their
derivatives: 
\[
\raisebox{.4em}{\rule{0pt}{10pt}}
\raisebox{-.3em}{\rule{0pt}{10pt}}
\delta \left(\hskip 1em
\unitlength .4cm
\begin{picture}(0,0)(0,.2)
\put(-.5,0){\makebox(0,0)[cc]{\Large$\bullet$}}
\put(-.5,0){\vector(0,1){1.2}}
\end{picture} 
\right) =  0,\
\delta \left(\hskip 1em
\unitlength .4cm
\begin{picture}(0,0)(0,-.3)
\put(-.5,0){\makebox(0,0)[cc]{\Large$\bullet$}}
\put(-.5,0){\vector(0,1){1.2}}
\put(-.5,-1){\vector(0,1){.9}}
\end{picture} 
\right) =  
\unitlength 2.2 mm
\begin{picture}(4,1)(-2,-.6)
\put(0,0.6){\vector(0,1){1.7}}
\put(0,-2){\vector(0,1){1.78}}
\put(1.9,-1.9){\vector(-1,1){1.62}}
\put(1.7,-1.7){\makebox(0,0)[cc]{\Large$\bullet$}}
\put(0,0){\makebox(0,0)[cc]{\Large$\circ$}}
\end{picture},\
\delta \left(
\begin{picture}(4,1)(-2,-.6)
\put(0,0.3){\vector(0,1){1.7}}
\put(-1,-2){\vector(1,2){.91}}
\put(1,-2){\vector(-1,2){.91}}
\put(0,0){\makebox(0,0)[cc]{\Large$\bullet$}}
\end{picture}
\right) = -
\unitlength 2.2mm
\begin{picture}(4,1)(-2,0)
\put(0,0.6){\vector(0,1){2.5}}
\put(0,0.4){\vector(0,1){1}}
\put(-1,-2.3){\vector(1,2){.89}}
\put(1,-2.3){\vector(-1,2){.89}}
\put(0,1.5){\makebox(0,0)[cc]{\Large$\bullet$}}
\put(0,-0.2){\makebox(0,0)[cc]{\Large$\circ$}}
\end{picture} 
+
\unitlength 2.2mm
\begin{picture}(4,1)(-2,-.6)
\put(0,0.6){\vector(0,1){1.7}}
\put(-1.3,-2.7){\vector(1,2){1.2}}
\put(1.3,-2.7){\vector(-1,2){1.2}}
\put(1.3,-2.7){\vector(-1,2){.4}}
\put(.75,-1.7){\makebox(0,0)[cc]{\Large$\bullet$}}
\put(0,0){\makebox(0,0)[cc]{\Large$\circ$}}
\end{picture}
+
\unitlength 2.2mm
\begin{picture}(4,1)(-2,-.6)
\put(0,0.6){\vector(0,1){1.7}}
\put(-1,-2){\vector(1,2){.88}}
\put(1,-2){\vector(-1,2){.88}}
\put(2,-1){\vector(-2,1){1.7}}
\put(2,-1){\makebox(0,0)[cc]{\Large$\bullet$}}
\put(0,0){\makebox(0,0)[cc]{\Large$\circ$}}
\end{picture} \hskip .7em ,\ \ldots
\]
It is also clear that $\delta \left(\mbox{\anchor}\right) =  0$.

\vskip .4em
Let us denote by $\Gr^*_{T \ot T,T} \subset \Gr^*_{T \times T,T}$ the
subcomplex describing bilinear operators. Its degree~$0$ part
$\Gr^0_{T \ot T,T}$ is spanned by
\[
\unitlength .5cm
\begin{picture}(0,0)(0,0)
\put(0,2){\makebox(0,0)[cc]{\hskip .5mm$\bbox$}}
\put(0,1){\vector(0,1){.935}}
\put(0,1){\makebox(0,0)[cc]{\Large$\bullet$}}
\put(.4,1){\makebox(0,0)[lc]{\scriptsize$Y$}}
\put(0,0){\vector(0,1){.935}}
\put(0,0){\makebox(0,0)[cc]{\Large$\bullet$}}
\put(.4,0){\makebox(0,0)[lc]{\scriptsize$X$}}
\end{picture}
\hskip .7cm \raisebox{.6cm}{,}\hskip .8cm
\unitlength .5cm
\begin{picture}(0,0)(0,0)
\put(0,2){\makebox(0,0)[cc]{\hskip .5mm$\bbox$}}
\put(0,1){\vector(0,1){.935}}
\put(0,1){\makebox(0,0)[cc]{\Large$\bullet$}}
\put(.4,1){\makebox(0,0)[lc]{\scriptsize$X$}}
\put(0,0){\vector(0,1){.935}}
\put(0,0){\makebox(0,0)[cc]{\Large$\bullet$}}
\put(.4,0){\makebox(0,0)[lc]{\scriptsize$Y$}}
\end{picture}
\hskip .7cm \raisebox{.6cm}{,}\hskip 1.5cm
\unitlength .6cm
\begin{picture}(0,1.5)(1.1,.6)
\put(0,2){\makebox(0,0)[cc]{\hskip .5mm$\bbox$}}
\put(0,1){\vector(0,1){.935}}
\put(0,1){\makebox(0,0)[cc]{\Large$\bullet$}}
\put(.4,1){\makebox(0,0)[lc]{\scriptsize$X$}}
\put(0,.2){
\put(2,1){\makebox(0,0)[cc]{\circle {1.5}}}
\put(1.3,1.35){\makebox(0,0)[cc]{\Large$\bullet$}}
\put(1.32,1.25){\makebox(0,0)[tc]{\vector(0,1){0}}}
\put(1.1,1.55){\makebox(0,0)[rb]{\scriptsize $Y$}}
}
\end{picture}
\hskip 1.5 cm \raisebox{.6cm}{and}\hskip 1.5cm
\unitlength .6cm
\begin{picture}(0,1.5)(1.1,.6)
\put(0,2){\makebox(0,0)[cc]{\hskip .5mm$\bbox$}}
\put(0,1){\vector(0,1){.935}}
\put(0,1){\makebox(0,0)[cc]{\Large$\bullet$}}
\put(.4,1){\makebox(0,0)[lc]{\scriptsize$Y$}}
\put(0,.2){
\put(2,1){\makebox(0,0)[cc]{\circle {1.5}}}
\put(1.3,1.35){\makebox(0,0)[cc]{\Large$\bullet$}}
\put(1.32,1.25){\makebox(0,0)[tc]{\vector(0,1){0}}}
\put(1.1,1.55){\makebox(0,0)[rb]{\scriptsize $X$}}
}
\end{picture}
\hskip 1.2 cm \raisebox{.6cm}{.}
\]
\vglue -.2cm
\noindent 
One easily calculates the differential of the leftmost term:
\begin{equation}
\label{ref1}
\delta
\left(\hskip .2em
\unitlength .5cm
\begin{picture}(1,1.5)(-0.3,0.7)
\put(0,2){\makebox(0,0)[cc]{\hskip .5mm$\bbox$}}
\put(0,1){\vector(0,1){.935}}
\put(0,1){\makebox(0,0)[cc]{\Large$\bullet$}}
\put(.4,1){\makebox(0,0)[lc]{\scriptsize$Y$}}
\put(0,0){\vector(0,1){.935}}
\put(0,0){\makebox(0,0)[cc]{\Large$\bullet$}}
\put(.4,0){\makebox(0,0)[lc]{\scriptsize$X$}}
\end{picture} \hskip .2em
\right)= 
\hskip .5cm
\unitlength .5cm
\begin{picture}(1,1.5)(-0.3,0.7)
\put(0,2){\makebox(0,0)[cc]{\hskip .5mm$\bbox$}}
\put(0,1){\vector(0,1){.935}}
\put(0,1){\makebox(0,0)[cc]{\Large$\bullet$}}
\put(.4,1){\makebox(0,0)[lc]{\scriptsize$Y$}}
\put(0,0){\vector(0,1){.935}}
\put(0,0){\makebox(0,0)[cc]{\Large$\bullet$}}
\put(.4,0){\makebox(0,0)[lc]{\scriptsize$X$}}
\put(-1.1,2){\makebox(0,0)[lc]{$\delta( \hskip .4cm )$}}
\end{picture} \hskip .3cm
+
\hskip .5cm
\unitlength .5cm
\begin{picture}(1,1.5)(-0.3,0.7)
\put(0,2){\makebox(0,0)[cc]{\hskip .5mm$\bbox$}}
\put(0,1){\vector(0,1){.935}}
\put(0,1){\makebox(0,0)[cc]{\Large$\bullet$}}
\put(.4,1){\makebox(0,0)[lc]{\scriptsize$Y$}}
\put(0,0){\vector(0,1){.935}}
\put(0,0){\makebox(0,0)[cc]{\Large$\bullet$}}
\put(.4,0){\makebox(0,0)[lc]{\scriptsize$X$}}
\put(0,1){\makebox(0,0)[cc]{$\delta( \hskip .6cm )$}}
\end{picture} \hskip .4cm
+
\hskip .4cm
\unitlength .5cm
\begin{picture}(1,1.5)(-0.3,0.7)
\put(0,2){\makebox(0,0)[cc]{\hskip .5mm$\bbox$}}
\put(0,1){\vector(0,1){.935}}
\put(0,1){\makebox(0,0)[cc]{\Large$\bullet$}}
\put(.4,1){\makebox(0,0)[lc]{\scriptsize$Y$}}
\put(0,0){\vector(0,1){.935}}
\put(0,0){\makebox(0,0)[cc]{\Large$\bullet$}}
\put(.4,0){\makebox(0,0)[lc]{\scriptsize$X$}}
\put(0,0){\makebox(0,0)[cc]{$\delta( \hskip .6cm )$}}
\end{picture} 
\hskip .3cm
=
\hskip .4cm
\unitlength .5cm
\begin{picture}(1,1.5)(-0.3,0.7)
\put(0,2){\makebox(0,0)[cc]{\hskip .5mm$\bbox$}}
\put(0,1.15){\vector(0,1){.78}}
\put(0,1){\makebox(0,0)[cc]{\Large$\circ$}}
\put(1.3,.7){\makebox(0,0)[lb]{\scriptsize$Y$}}
\put(0,0){\vector(0,1){.9}}
\put(0,0){\makebox(0,0)[cc]{\Large$\bullet$}}
\put(.4,0){\makebox(0,0)[lc]{\scriptsize$X$}}
\put(1,0.5){\makebox(0,0)[cc]{\Large$\bullet$}}
\put(1,.5){\vector(-2,1){.88}}
\end{picture} 
\hskip 2em \in \Gr^1_{T \ot T,T}
\end{equation}
and similarly one gets
\[
\delta\left(
\unitlength .5cm
\begin{picture}(3.2,1.2)(-.4,1.3)
\put(0,2){\makebox(0,0)[cc]{\hskip .5mm$\bbox$}}
\put(0,1){\vector(0,1){.9}}
\put(0,1){\makebox(0,0)[cc]{\Large$\bullet$}}
\put(.4,1){\makebox(0,0)[lc]{\scriptsize$X$}}
\put(0,.2){
\put(2,1){\makebox(0,0)[cc]{\circle {1.5}}}
\put(1.3,1.35){\makebox(0,0)[cc]{\Large$\bullet$}}
\put(1.32,1.25){\makebox(0,0)[tc]{\vector(0,1){0}}}
\put(1.1,1.55){\makebox(0,0)[rb]{\scriptsize $Y$}}
}
\end{picture}
\right) =
\unitlength .5cm
\begin{picture}(3.2,1.2)(-.4,1.3)
\put(0,2){\makebox(0,0)[cc]{\hskip .5mm$\bbox$}}
\put(0,1){\vector(0,1){.9}}
\put(0,1){\makebox(0,0)[cc]{\Large$\bullet$}}
\put(.4,1){\makebox(0,0)[lc]{\scriptsize$X$}}
\put(1,.2){
\put(2.09,.85){\makebox(0,0)[cc]{\oval(1.5,1.5)[b]}}
\put(2.09,1.15){\makebox(0,0)[cc]{\oval(1.5,1.5)[t]}}
\put(2.85,1.22){\line(0,1){.3}}
\put(.5,.5){\vector(1,1){.75}}
\put(.6,.6){\makebox(0,0)[cc]{\Large$\bullet$}}
\put(1.35,1.35){\makebox(0,0)[cc]{\Large$\circ$}}
\put(1.32,1.25){\makebox(0,0)[tc]{\vector(0,1){0}}}
\put(.9,1.15){\makebox(0,0)[rb]{\scriptsize $Y$}}
}
\end{picture} \hskip 2.2em \in \Gr^1_{T \ot T,T}\ .
\]
The formula for the differential of the remaining two generators of
$\Gr_{T\ot T,T}^0$ is
obtained by interchanging $X \leftrightarrow Y$ in the previous two displays.
One clearly has
\[
\raisebox{.7em}{\rule{0pt}{12pt}}
\delta
\left(
\unitlength .5cm
\begin{picture}(1,1.5)(-0.3,0.7)
\put(0,2){\makebox(0,0)[cc]{\hskip .5mm$\bbox$}}
\put(0,1){\vector(0,1){.935}}
\put(0,1){\makebox(0,0)[cc]{\Large$\bullet$}}
\put(.4,1){\makebox(0,0)[lc]{\scriptsize$Y$}}
\put(0,0){\vector(0,1){.935}}
\put(0,0){\makebox(0,0)[cc]{\Large$\bullet$}}
\put(.4,0){\makebox(0,0)[lc]{\scriptsize$X$}}
\end{picture}
-
\unitlength .5cm
\begin{picture}(1,1.5)(-0.3,0.7)
\put(0,2){\makebox(0,0)[cc]{\hskip .5mm$\bbox$}}
\put(0,1){\vector(0,1){.935}}
\put(0,1){\makebox(0,0)[cc]{\Large$\bullet$}}
\put(.4,1){\makebox(0,0)[lc]{\scriptsize$X$}}
\put(0,0){\vector(0,1){.935}}
\put(0,0){\makebox(0,0)[cc]{\Large$\bullet$}}
\put(.4,0){\makebox(0,0)[lc]{\scriptsize$Y$}}
\end{picture}
\right)= 
\hskip .4cm
\unitlength .5cm
\begin{picture}(1,1.5)(-0.3,0.7)
\put(0,2){\makebox(0,0)[cc]{\hskip .5mm$\bbox$}}
\put(0,1.15){\vector(0,1){.78}}
\put(0,1){\makebox(0,0)[cc]{\Large$\circ$}}
\put(1.3,.7){\makebox(0,0)[lb]{\scriptsize$Y$}}
\put(0,0){\vector(0,1){.9}}
\put(0,0){\makebox(0,0)[cc]{\Large$\bullet$}}
\put(.4,0){\makebox(0,0)[lc]{\scriptsize$X$}}
\put(1,0.5){\makebox(0,0)[cc]{\Large$\bullet$}}
\put(1,.5){\vector(-2,1){.88}}
\end{picture} 
\hskip .8cm -\hskip .2cm
\unitlength .5cm
\begin{picture}(1,1.5)(-0.3,0.7)
\put(0,2){\makebox(0,0)[cc]{\hskip .5mm$\bbox$}}
\put(0,1.15){\vector(0,1){.78}}
\put(0,1){\makebox(0,0)[cc]{\Large$\circ$}}
\put(1.3,.7){\makebox(0,0)[lb]{\scriptsize$X$}}
\put(0,0){\vector(0,1){.9}}
\put(0,0){\makebox(0,0)[cc]{\Large$\bullet$}}
\put(.4,0){\makebox(0,0)[lc]{\scriptsize$Y$}}
\put(1,0.5){\makebox(0,0)[cc]{\Large$\bullet$}}
\put(1,.5){\vector(-2,1){.88}}
\end{picture} \hskip .8cm =0,
\]
because the inputs of white vertices are symmetric. 
It is easy to verify that the element
\begin{equation}
\label{Pozitri_snad_prijede_PP.}
\raisebox{-1em}{\rule{0pt}{0pt}}
\bfb :=  \hskip .3em
\unitlength .5cm
\begin{picture}(1,1.5)(-0.3,0.7)
\put(0,2){\makebox(0,0)[cc]{\hskip .5mm$\bbox$}}
\put(0,1){\vector(0,1){.935}}
\put(0,1){\makebox(0,0)[cc]{\Large$\bullet$}}
\put(.4,1){\makebox(0,0)[lc]{\scriptsize$Y$}}
\put(0,0){\vector(0,1){.935}}
\put(0,0){\makebox(0,0)[cc]{\Large$\bullet$}}
\put(.4,0){\makebox(0,0)[lc]{\scriptsize$X$}}
\end{picture}
-
\unitlength .5cm
\begin{picture}(1,1.5)(-0.3,0.7)
\put(0,2){\makebox(0,0)[cc]{\hskip .5mm$\bbox$}}
\put(0,1){\vector(0,1){.935}}
\put(0,1){\makebox(0,0)[cc]{\Large$\bullet$}}
\put(.4,1){\makebox(0,0)[lc]{\scriptsize$X$}}
\put(0,0){\vector(0,1){.935}}
\put(0,0){\makebox(0,0)[cc]{\Large$\bullet$}}
\put(.4,0){\makebox(0,0)[lc]{\scriptsize$Y$}}
\end{picture}
\hskip .5em \in \Gr^0_{T \ot T,T}
\end{equation}
representing the Lie bracket in fact spans all cochains in
$\Gr^0_{T\ot T,T}$. We conclude that $H^0(\Gr^*_{T\ot T,T}, \delta)$
is one-dimensional, generated by the bracket $[X,Y]$.  The
complex $\Gr^*_{T \ot T,T}$ clearly fits into the scheme discussed in
Remark~\ref{Snad_mam_hodnoceni_AAA} (with $C = 2$), which proves 
Theorem~\ref{main-A} for $d=2$.
\end{example}

\begin{example}
We close this section by an example suggested by the referee which
will further illuminate the meaning of the graph differential. The graph
\begin{equation}
\label{ref2}
\unitlength .5cm
\begin{picture}(0,1.2)(0,1)
\put(0,2){\makebox(0,0)[cc]{\hskip .5mm$\bbox$}}
\put(0,1){\vector(0,1){.935}}
\put(0,1){\makebox(0,0)[cc]{\Large$\bullet$}}
\put(.4,1){\makebox(0,0)[lc]{\scriptsize$Y$}}
\put(0,0){\vector(0,1){.935}}
\put(0,0){\makebox(0,0)[cc]{\Large$\bullet$}}
\put(.4,0){\makebox(0,0)[lc]{\scriptsize$X$}}
\end{picture}
\raisebox{-.8cm}{\rule{0em}{0em}}
\end{equation}
represents the local expression 
\begin{equation}
\label{zpatky_z_Polska}
\left(X^i\frac{\pa \hphantom{x^i}}{\pa x^i},
Y^i\frac{\pa \hphantom{x^i}}{\pa x^i}\right)
\longmapsto
X^i \frac{\pa Y^j}{\pa x^i}\frac{\pa \hphantom{x^j}}{\pa x^j}.
\end{equation}
If $\{y^i\}$  is a different set of coordinates, then $X$ and $Y$ transforms
to $X^i \frac{\pa y^s}{\pa x^i}\frac{\pa \hphantom{y^s}}{\pa {y^s}}$
and $Y^j \frac{\pa y^r}{\pa x^j}\frac{\pa \hphantom{y^r}}{\pa {y^r}}$,
respectively. Having this transformed $X$ act on the transformed $Y$
gives
\[
X^i \frac{\pa y^s}{\pa x^i}\frac{\pa \hphantom{y^s}}{\pa y^s}
\left(
Y^j \frac{\pa y^r}{\pa x^j}
\right)
\frac{\pa \hphantom{y^r}}{\pa y^r}
=
X^i \frac{\pa y^s}{\pa x^i}\frac{\pa Y^j}{\pa y^s}
\frac{\pa y^r}{\pa x^j}
\frac{\pa \hphantom{y^r}}{\pa y^r}
+
X^i \frac{\pa y^s}{\pa x^i} Y^j 
\frac{\pa \hphantom{y^s}}{\pa y^s}
\left(
\frac{\pa y^r}{\pa x^j}
\right)
\frac {\pa \hphantom{y^r}}{\pa y^r}.
\]

The first term in the right hand side is equal to the expression 
in~(\ref{zpatky_z_Polska}) under
change-of-coordinates, so the second term represents the extent to
which this expression is not invariant. It is equal to 
$X^iY^j \pa^2 y^r / \pa x^i \pa x^j \
\pa/{\pa y^r}$,
which translates directly to the formula~(\ref{ref1}) for the
differential of~(\ref{ref2}) in the graph complex.
\end{example}

\section{Operations on vector fields}
\label{5}

In this section we consider differential operators acting on a
finite number of vector fields $X,Y,Z,\ldots$ with values in vector
fields, that is, operators in $\Nat(T^{\times \infty},T) := \bigcup_{d \geq
  0} \Nat(T^{\times d},T)$
The first statement of this section is:

\begin{theorem}
\label{main-A}
Let $M$ be a smooth manifold and $d$ a natural number such that 
$\dim(M) \geq d$. Then each
$d$-multilinear natural operator from vector fields to
vector fields is a sum of iterations of the  Lie bracket
containing each of $d$ variables precisely once, and
all relations between these expressions  follow from
the Jacobi identity and antisymmetry. 
In particular, there are precisely $(d-1)!$ linearly independent
operators of the above type.
\end{theorem}

Theorem~\ref{main-A} is an obvious consequence of
Proposition~\ref{Aspon_Viktor_napsal.} below and the formula for the
dimension of the $k$th piece of the operad $\Lie$ for Lie algebras
that can be found for example 
in~\cite[Example~3.1.12]{ginzburg-kapranov:DMJ94}.
Theorem~\ref{main-A} 
describes {\em multilinear\/} operators and does not cover operators
as $\gO(X,Y,Z) := [X,Y] + [X,[X,Z]]$
but can easily be extended to cover also these cases. Since
all operators are assumed to be polynomial, they decompose into the sum of
their homogeneous parts. For instance, $\gO(X,Y,Z)$
is the sum of the homogeneity-2 part $[X,Y]$ and the homogeneity-3
part $[X,[X,Z]]$.

\begin{remark}
\label{pol1}
Let us explain the decomposition of operators $\gO
\in \Nat(T^{\times \infty},T)$ into homogeneous parts in more detail. 
The local formula $O$
for the operator $\gO$ is the sum $O = O_1 + \cdots + O_r$, where
$O_d$ is the part of $O$ consisting of terms with precisely $d$
occurrences of the vector field variables. The action of the structure
group $\GL {\infty} n$ on the typical fiber of the prolongation of
$T^{\times \infty}$ is linear, which is expressed by the manifest
linearity of the replacement rule~(\ref{Nikdo}) in the vector field
variable. This implies that the map $O$ is $\GL {\infty}
n$-equivariant if and only if each of its homogeneous components $O_d$
is $\GL {\infty} n$-equivariant, $1 \leq d \leq r$. Therefore $\gO =
\gO_1 + \cdots + \gO_r$, where $\gO_d$ is the operator defined by the
local formula $O_d$, $1 \leq d \leq r$.

We conclude that to classify operators of the above type, it suffices
to classify homogeneous operators. It is a standard fact that
each homogeneous operator of degree $d$ is either
$d$-multilinear or a sum of operators 
obtained from $d$-multilinear operators by
repeating one or more of their variables. We will call this procedure the
{\em depolarization\/} of multilinear operators.
Theorem~\ref{main-A} therefore implies the following corollary. 
\end{remark}

\begin{corollary}
\label{B}
Let $M$ be a smooth manifold. Each
natural differential operator from vector fields on $M$ to 
vector fields on $M$ whose all components are of homogeneity $\leq
\dim(M)$ is a sum of iterations of the Lie bracket. All
relations between these iterations  follow from the Jacobi identity
and antisymmetry. 
\end{corollary}

In Example~\ref{Ted_mam_kasel.} we studied the graph complex
$\Gr^*_{T \ot T,T}$ describing bilinear operators. Bearing this\label{AAA}
example in mind, we introduce 
$\Grbull^*(d)= \Gr^*_{T^{\otimes d},T} \subset \Gr^*_{T^{\times d},T}$, the
subcomplex describing $d$-multilinear operators. Its degree $m$
component is 
spanned by graphs with $d$ vertices of the first type labelled by
$\Rada X1d$, $m$ white vertices of the third type and one 2nd type
vertex \anchor\ which we call the {\em anchor\/}.
Observe that $\Grbull^*(d)$ is precisely the graph complex
$\Gr^*_{\bullet(b)\nabla(c)}$ of \cite[Corollary~5.1]{markl:ig}
with $b:= d$ and $c := 0$. 
The collection $\Grbull^0 = \{\Grbull^0(d)\}_{d \geq 1}$ of degree~$0$
subspaces 
admits two types of operations.

(i) For graphs $G' \in \Grbull^0(u)$, $G'' \in \Grbull^0(v)$ and $1
\leq i \leq u$, one has the $\circ_i$-product $G' \circ_i G'' \in
\Grbull^0(u+v-1)$ given by the following straightforward extension of
the Chapoton-Livernet vertex
insertion~\cite[\S~1.5]{chapoton-livernet:pre-lie} to non-simply
connected graphs.
Assume that $\Rada {X'}1u$ are the black vertices of
$G'$, $\Rada {X''}1v$ the black vertices of $G''$ and ${\it In}
(X'_i)$ the set of inputs of $X'_i$ in $G'$. Then
\[
G' \circ_i G'' := 
\sum_{
f:{\it In} (X'_i)\to \{\Rada {X''}1v \}
} G' \circ_i^f G'' \in \Grbull^0(u+v-1),
\]
where $G' \circ^f_i G'' \in \Grbull^0(u+v-1)$ 
is the graph obtained by replacing the vertex $X'_i$ of $G'$ by $G''$ and
grafting the inputs of  $X'_i$ on black vertices of $G''$ following
$f$.

In more detail, one starts by cutting off the anchor \anchor\ of $G''$ 
and grafts the resulting free edge on the 
vertex of $G'$ immediately above $X'_i$. Then one grafts 
each input edge $e$ of $X'_i$ on the vertex $f(e)$ of $G''$. 
Finally, one changes the labels $\Rada {X'}1{i-1},\Rada {X''}1v,\Rada
{X'}{i+1}u$ of the black vertices of the graph obtained in this way 
into $\Rada X1{u+v-1}$. 

(ii) One has the right action of the symmetric
group: for each $G \in \Grbull^0(d)$ and a~permutation $\sigma \in
\Sigma_d$, one has $G\sigma \in \Grbull^0(d)$ given by permuting the
labels $\Rada X1d$ of the black vertices of $G$ according to $\sigma$.

\begin{proposition}
\label{zpatky_v_Koline}
The collection $\Grbull^0 = \{\Grbull^0(d)\}_{d \geq 1}$ with the above
operations is an operad
with unit \hskip .2em 
$\uunit \hskip .3em \in \Grbull^0(1)$~\cite{markl-shnider-stasheff:book}. 
The operad structure of $\Grbull^0$ restricts to
$H^0(\Grbull^*,\delta) = \Ker(\delta: \Grbull^0 \to \Grbull^1)$.
\end{proposition}

\begin{proof}

The operad axioms for the operations in (i) and (ii) above are
verified directly, compare
also~\cite[\S~1.5]{chapoton-livernet:pre-lie}. The simplest way to see
that the operad structure of $\Grbull^0$ restricts to the kernel of
$\delta$ is to extend the operations (i) and (ii), in the obvious
manner, to the graded collection $\Grbull^*$, making
$(\Grbull^*,\delta)$ a dg-operad. This, in particular, would mean that
$\delta$ is a derivation with respect to these extended $\circ_i$-operations,
which implies the second part of the proposition.
\end{proof}

\begin{example}
\label{psanovIHES}
An instructive example of the vertex insertion can be found
in~\cite[\S~1.5]{chapoton-livernet:pre-lie}. We
present here a simpler one, taken from the proof 
of~\cite[Theorem~1.9]{chapoton-livernet:pre-lie}.
Let $\bp$ be the graph
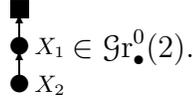
\[
\unitlength .5cm
\begin{picture}(0,1)(1,1)
\put(0,2){\makebox(0,0)[cc]{\hskip .5mm$\bbox$}}
\put(0,1){\vector(0,1){.935}}
\put(0,1){\makebox(0,0)[cc]{\Large$\bullet$}}
\put(.4,1){\makebox(0,0)[lc]{\scriptsize$X_1$}}
\put(1.4,1){\makebox(0,0)[lc]{$\in \Grbull^0(2)$.}}
\put(0,0){\vector(0,1){.935}}
\put(0,0){\makebox(0,0)[cc]{\Large$\bullet$}}
\put(.4,0){\makebox(0,0)[lc]{\scriptsize$X_2$}}
\end{picture}
\]
Then one has
\[
\adjust {1.6}
\unitlength .5cm
\begin{picture}(3,2.5)(0,0)
\put(0,2){\makebox(0,0)[cc]{\hskip .5mm$\bbox$}}
\put(0,1){\vector(0,1){.935}}
\put(0,1){\makebox(0,0)[cc]{\Large$\bullet$}}
\put(-.4,.5){\makebox(0,0)[rc]{$\bp \circ_1 \bp=$}}
\put(.4,1){\makebox(0,0)[lc]{\scriptsize$X_1$}}
\put(0,0){\vector(0,1){.935}}
\put(0,0){\makebox(0,0)[cc]{\Large$\bullet$}}
\put(.4,0){\makebox(0,0)[lc]{\scriptsize$X_2$}}
\put(0,-1){\vector(0,1){.935}}
\put(0,-1){\makebox(0,0)[cc]{\Large$\bullet$}}
\put(.4,-1){\makebox(0,0)[lc]{\scriptsize$X_3$}}
\put(3.3,-.5){
\put(0,2){\makebox(0,0)[cc]{\hskip .5mm$\bbox$}}
\put(0,1.15){\vector(0,1){.78}}
\put(-.4,1){\makebox(0,0)[rc]{\scriptsize$X_1$}}    
\put(-1.3,1){\makebox(0,0)[rc]{$+$}}    
\put(0,1){\makebox(0,0)[cc]{\Large$\bullet$}}
\put(2.3,1){\makebox(0,0)[cl]{$\in \Grbull^0(3)$  \hskip .5em and}}  
\put(1.3,.7){\makebox(0,0)[lb]{\scriptsize$X_3$}}
\put(0,0){\vector(0,1){.9}}
\put(0,0){\makebox(0,0)[cc]{\Large$\bullet$}}
\put(.4,0){\makebox(0,0)[lc]{\scriptsize$X_2$}}
\put(1,0.5){\makebox(0,0)[cc]{\Large$\bullet$}}
\put(1,.5){\vector(-2,1){.88}}
}
\end{picture}
\hskip 14em 
\unitlength .5cm
\begin{picture}(4,2.5)(0,.1)
\put(0,2){\makebox(0,0)[cc]{\hskip .5mm$\bbox$}}
\put(0,1){\vector(0,1){.935}}
\put(0,1){\makebox(0,0)[cc]{\Large$\bullet$}}
\put(-.4,.5){\makebox(0,0)[rc]{$\bp \circ_2 \bp=$}}
\put(3,.6){\makebox(0,0)[cc]{$\in \Grbull^0(3)$.}} 
\put(.4,1){\makebox(0,0)[lc]{\scriptsize$X_1$}}
\put(0,0){\vector(0,1){.935}}
\put(0,0){\makebox(0,0)[cc]{\Large$\bullet$}}
\put(.4,0){\makebox(0,0)[lc]{\scriptsize$X_2$}}
\put(0,-1){\vector(0,1){.935}}
\put(0,-1){\makebox(0,0)[cc]{\Large$\bullet$}}
\put(.4,-1){\makebox(0,0)[lc]{\scriptsize$X_3$}}
\end{picture}
\]
The above display implies that the associator 
${\it Ass}(\bp) :=  \bp \circ_1 \bp - \bp \circ_2 \bp$ equals
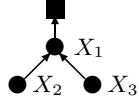
\[
\adjust{.8}
\unitlength .5cm
\begin{picture}(0,1.45)(0,.55)
\put(0,2){\makebox(0,0)[cc]{\hskip .5mm$\bbox$}}
\put(0,1.08){\vector(0,1){.85}}
\put(0,1){\makebox(0,0)[cc]{\Large$\bullet$}}
\put(0.5,1){\makebox(0,0)[lc]{\scriptsize$X_1$}}
\put(-1,-1){
\put(0,1){{\vector(1,1){.92}}}
\put(0,1){\makebox(0,0)[cc]{\Large$\bullet$}}
\put(.4,1){\makebox(0,0)[lc]{\scriptsize$X_2$}}
}
\put(1,-1){
\put(0,1){{\vector(-1,1){.92}}}
\put(0,1){\makebox(0,0)[cc]{\Large$\bullet$}}
\put(.4,1){\makebox(0,0)[lc]{\scriptsize$X_3$}}
}
\end{picture}
\]
and is therefore symmetric in $X_2$ and $X_3$. 
This, by definition, means that $\bp$ represents a pre-Lie
multiplication~\cite[\S~1.1]{chapoton-livernet:pre-lie}. We will see
that $\Grbull^0$ is indeed closely related to the pre-Lie operad $\pLie$.
\end{example}

Let $\tau \in \Sigma_2$ be the generator.  By standard properties of
pre-Lie algebras~\cite[Proposition~1.2]{chapoton-livernet:pre-lie},
the antisymmetrization \hbox{$\bp(\tau - \id)$} of the element
$\bp$ from Example~\ref{psanovIHES} is a Lie bracket. Observe that
\hbox{$\bp(\tau - \id)$}
equals the element $\bfb$ introduced
in~(\ref{Pozitri_snad_prijede_PP.}).

\begin{proposition}
\label{Aspon_Viktor_napsal.}
The $0$th cohomology $H^0(\Grbull^*(d),\delta)$ is, for each $d \geq
2$, generated by the Lie bracket $\bfb=\bp(\tau - \id) 
\in H^0(\Grbull^*(2),\delta)$, by
iterating operations~(i) and~(ii) above. There
are no relations between these iterations other than those following
from the Jacobi identity and antisymmetry.
\end{proposition}

A compact formulation of Proposition~\ref{Aspon_Viktor_napsal.} is that  
the operad $H^0(\Grbull^*,\delta) = \{H^0(\Grbull^*(d),\delta)\}_{d \geq 1}$ 
is isomorphic to the operad $\Lie = \{\Lie(d)\}_{d \geq 1}$ for Lie
algebras~\cite[Example~II.3.34]{markl-shnider-stasheff:book}, via an
isomorphism that sends the generator 
$\beta \in \Lie(2)$  of $\Lie$ into  $\bfb \in \Grbull^0(2)$. 
Graphs spanning $\Grbull^0(d)$ have $d$
edges which explains the stability condition $\dim(M) \geq d$ in
Theorem~\ref{main-A}.
The rest of this section is devoted to a proof of its main result.

\begin{proof}[Proof of Proposition~\ref{Aspon_Viktor_napsal.}] 
It is clear from formulas~(\ref{Co_Sergej_planuje?}),~(\ref{Nikdo})
and $\delta \left(\mbox{\anchor}\right) = 0$ that the differential
preserves connected components of underlying graphs. Therefore, for
each $d \geq 1$, $\Grbull^*(d)$ is the direct sum $\Grbull^*(d) =
\bigoplus_{c \geq 1}\Gr_{\bullet c}^*(d)$, where $\Gr_{\bullet
c}^*(d)$ denotes the subcomplex spanned by graphs with $c$ connected
components. In particular, $\Gr_{\bullet 1}^*(d)$ is the subcomplex of
connected graphs. It is easy to see that $\Gr^0_{\bullet1}$ is a
suboperad of~$\Grbull^*$.

As the Lie bracket represented by $\bfb \in \Gr^0_{\bullet1}(2)$ is
antisymmetric and satisfies the Jacobi
identity, the rule $F(\beta) := \bfb$, where $\beta \in
\Lie(2)$ is the generator, defines an operad homomorphism 
$F : \Lie \to \Gr^0_{\bullet1}$. Since the
Lie bracket and its iterations are natural operators, $\Im(F) \subset
\Ker(\delta :  \Gr^0_{\bullet1} \to  \Gr^1_{\bullet1})$.
Proposition~\ref{Aspon_Viktor_napsal.} will clearly be established
if we prove that

(i) the operad map $F : \Lie \to \Gr^0_{\bullet1}$ induces an isomorphism
$\Lie \cong H^0(\Gr_{\bullet 1}^*,\delta)$, and

(ii) $H^0(\Gr_{\bullet c}^*(d),\delta) = 0$, for each $c \geq 2$, $d \geq 1$.

Part~(i) is highly nontrivial, but it in fact has already been proved
in~\cite{markl:JLT07}. Indeed, the operad $\Gr^0_{\bullet1}$ 
is precisely the operad
$\pLie$ describing pre-Lie algebras~\cite{chapoton-livernet:pre-lie} 
and $F : \Lie \to \Gr^0_{\bullet1}$ corresponds, under the identification
$\Gr^0_{\bullet1} \cong \pLie$, to 
the inclusion $\iota : \Lie \hookrightarrow \pLie$ 
induced by the antisymmetrization of the
pre-Lie product. The dg operad ${{\rm r}{\sf pL}}^*$
of~\cite{markl:JLT07}  coincides, in
degrees $0$ and $1$, with the complex $\Gr_{\bullet 1}^*$ and the
isomorphism in (i) is isomorphism~(2) of~\cite{markl:JLT07}.

Let us prove~(ii). For each $m \geq 0$, $d \geq 1$, consider 
the span $\Gr_{\bullet\wc}^m(d)$ of
connected graphs with $d$ vertices $\Rada X1d$ of type~1, $m$ `white'
vertices of type~3 and no vertex of
type~2. The direct sum
$\Gr_{\bullet\wc}^*(d) = \bigoplus_{m \geq 0} \Gr_{\bullet\wc}^m(d)$
is a cochain complex, with the differential
defined in the same way as the differential in $\Grbull^*(d)$ and denoted
again by $\delta$.
We claim that, for each $c \geq 2$ and $d \geq 1$, there is an
isomorphism of cochain complexes
\begin{equation}
\label{rozklad}
\Gr_{\bullet c}^*(d) \cong \bigoplus_{i_1 + \cdots + i_c = d} 
\Gr_{\bullet 1}^*(i_1) \otimes  \left(\Gr_{\bullet\wc}^*(i_2) 
\odot \cdots \odot \Gr_{\bullet\wc}^*(i_c)\right)
\end{equation}
where $\odot$ as usual denotes the symmetric product.
To prove this isomorphism, observe that each graph 
$G \in \Gr_{\bullet c}^*(d)$
decomposes into the disjoint union
\begin{equation}
\label{g}
G = G_1 \sqcup G_2 \sqcup \cdots \sqcup G_c,
\end{equation}
of its connected components. Precisely one of
these components contains the unique type 2 vertex
$\unitlength .25cm
\begin{picture}(1,1.4)(-1,-.7)
\put(-.45,.55){\makebox(0,0)[cc]{$\sbbox$}}
\put(-.5,-.8){\vector(0,1){1.2}}
\end{picture}$, 
assume it is $G_1$. Then $G_1 \in
\Gr_{\bullet 1}^*(i_1)$ 
and $G_s \in \Gr_{\bullet\wc}^*(i_s)$ for $2 \leq s \leq c$, with
some $i_1 + \cdots + i_c = d$. Decomposition~(\ref{g}) is clearly
unique up to the order of $\Rada G2c$ and is preserved by the differential.
This proves~(\ref{rozklad}). By K\"unneth and Mashke's theorems, (ii)
follows from $H^0(\Gr_{\bullet\wc}^*(d),\delta) = 0$, $d \geq 1$, 
which is the same as showing that
\begin{equation}
\label{central}
\mbox {the map $\delta : \Gr_{\bullet\wc}^0(d) \to
\Gr_{\bullet\wc}^1(d)$ is a
monomorphism for each $d \geq 1$.}
\end{equation}

Let us inspect the structure of~$\Gr_{\bullet\wc}^*(d)$. It is
clear from simple graph combinatorics that each graph in
$\Gr_{\bullet\wc}^m(d)$ has genus $1$, therefore it contains a unique
wheel. Denote $\Gr_{\bullet\wc}^m(d,w) \subset \Gr_{\bullet\wc}^m(d)$ 
the subspace
spanned by graphs that have precisely $w$ vertices (of either type) on
the wheel, $w \geq 0$.
It is obvious from~(\ref{Co_Sergej_planuje?}) and~(\ref{Nikdo}) that
$\delta(\Gr_{\bullet\wc}^m(d,w)) \subset \Gr_{\bullet\wc}^{m+1}(d,w) \oplus 
\Gr_{\bullet\wc}^{m+1}(d,w+1)$, for  $d \geq 1$, $w \geq 0$;
see also Figure~\ref{fig1}.
Let us denote by $\delta^0$ the component of $\delta$ that preserves
the number of vertices on the wheel and $\delta^1$ the component that
raises it by one.
\begin{figure}[t]
\[
\unitlength .37cm
\begin{picture}(20,12.5)(9,9)
\put(1,20){
\put(-1,0){\makebox(0,0)[r]{$\delta \left(\rule{0pt}{25pt} \right.$}}
\put(1,2){\makebox(0,0){\skelet}}
\put(0,0){\makebox(0,0){\Large$\bullet$}}
\put(0,-1){\vector(0,1){.81}}
\put(0,0.25){\line(0,1){.81}}
\put(3,0){\makebox(0,0)[l]{$\left.\rule{0pt}{25pt} \right)=$}}
\put(7,0){
\put(1,2){\makebox(0,0){\skelet}}
\put(0,0){\makebox(0,0){\Large$\circ$}}
\put(-1,-1){\makebox(0,0){\Large$\bullet$}}
\put(-1.2,-1.2){\vector(1,1){1}}
\put(0,-1){\vector(0,1){.81}}
\put(0,0.25){\line(0,1){.81}}
\put(3.5,0){$+$ \hskip .2cm replacements of remaining vertices of the graph,}
}
}
\put(1,15){
\put(-1,0){\makebox(0,0)[r]{$\delta \left(\rule{0pt}{25pt} \right.$}}
\put(1,2){\makebox(0,0){\skelet}}
\put(0,0){\makebox(0,0){\Large$\bullet$}}
\put(-1.2,-1.2){\vector(1,1){1}}
\put(0,-1){\vector(0,1){.81}}
\put(0,0.25){\line(0,1){.81}}
\put(3,0){\makebox(0,0)[l]{$\left.\rule{0pt}{25pt} \right)=$}}
\put(11,0){\makebox(0,0){$+$}}
\put(17,0){\makebox(0,0){$+$}}
\put(23,0){\makebox(0,0){$-$}}
\put(29,0){\makebox(0,0){$+ \cdots$}}
\put(7,0){
\put(1,2){\makebox(0,0){\skelet}}
\put(0,0){\makebox(0,0){\Large$\circ$}}
\put(-1,-1){\makebox(0,0){\Large$\bullet$}}
\put(-1,-2){\vector(0,1){.9}}
\put(-1.2,-1.2){\vector(1,1){1}}
\put(0,-1){\vector(0,1){.81}}
\put(0,0.25){\line(0,1){.81}}
}
\put(13,0){
\put(1,2){\makebox(0,0){\skelet}}
\put(0,0){\makebox(0,0){\Large$\circ$}}
\put(1.1,-1.1){\makebox(0,0){\Large$\bullet$}}
\put(-1.2,-1.2){\vector(1,1){1}}
\put(1.2,-1.2){\vector(-1,1){1}}
\put(0,-1){\vector(0,1){.81}}
\put(0,0.25){\line(0,1){.81}}
}
\put(19,0){
\put(1,2){\makebox(0,0){\skelet}}
\put(0,.78){\makebox(0,0){\Large$\circ$}}
\put(0,-.8){\makebox(0,0){\Large$\bullet$}}
\put(0,-.58){\vector(0,1){1.2}}
\put(0,-1){\vector(0,1){0}}
\put(-1.2,-.4){\vector(1,1){1}}
}
\put(25,0){
\put(1,2){\makebox(0,0){\skelet}}
\put(0,.78){\makebox(0,0){\Large$\bullet$}}
\put(0,-.8){\makebox(0,0){\Large$\circ$}}
\put(0,-.58){\vector(0,1){1.2}}
\put(0,-1){\vector(0,1){0}}
\put(-1.2,-2){\vector(1,1){1}}
}
}
\put(1,10){
\put(-1,0){\makebox(0,0)[r]{$\delta \left(\rule{0pt}{25pt} \right.$}}
\put(1,2){\makebox(0,0){\skelet}}
\put(0,0){\makebox(0,0){\Large$\bullet$}}
\put(-1.2,-1.2){\vector(1,1){1}}
\put(-.6,-1.2){\vector(1,2){.52}}
\put(0,-1){\vector(0,1){.81}}
\put(0,0.25){\line(0,1){.81}}
\put(3,0){\makebox(0,0)[l]{$\left.\rule{0pt}{25pt} \right)=$}}
\put(11,0){\makebox(0,0){$+$}}
\put(17,0){\makebox(0,0){$+$}}
\put(23,0){\makebox(0,0){$+$}}
\put(7,0){
\put(1,2){\makebox(0,0){\skelet}}
\put(0,0){\makebox(0,0){\Large$\circ$}}
\put(-1,-1){\makebox(0,0){\Large$\bullet$}}
\put(-1.6,-2.2){\vector(1,2){.51}}
\put(-.4,-2.2){\vector(-1,2){.51}}
\put(-1.2,-1.2){\vector(1,1){1}}
\put(0,-1){\vector(0,1){.81}}
\put(0,0.25){\line(0,1){.81}}
\put(6,0){
\put(1,2){\makebox(0,0){\skelet}}
\put(0,0){\makebox(0,0){\Large$\circ$}}
\put(-1,-1){\makebox(0,0){\Large$\bullet$}}
\put(-1.6,-2.2){\vector(1,2){.51}}
\put(-.69,-2){\vector(1,3){0.61}}
\put(-1.2,-1.2){\vector(1,1){1}}
\put(0,-1){\vector(0,1){.81}}
\put(0,0.25){\line(0,1){.81}}
}
\put(12,0){
\put(1,2){\makebox(0,0){\skelet}}
\put(0,0){\makebox(0,0){\Large$\circ$}}
\put(-.39,-1.1){\makebox(0,0){\Large$\bullet$}}
\put(-2,-2){\line(1,1){1.41}}
\put(-.39,-1.1){\vector(1,3){0.31}}
\put(-1.2,-1.2){\vector(1,1){1}}
\put(-.39,-2){\vector(0,1){.75}}
\put(0,-1){\vector(0,1){.81}}
\put(0,0.25){\line(0,1){.81}}
}
\put(18,0){
\put(1,2){\makebox(0,0){\skelet}}
\put(0,0){\makebox(0,0){\Large$\circ$}}
\put(-1.2,-1.2){\vector(1,1){1}}
\put(-.6,-1.2){\vector(1,2){.52}}
\put(0,-1){\vector(0,1){.81}}
\put(0,0.25){\line(0,1){.81}}
\put(1,-1){\makebox(0,0){\Large$\bullet$}}
\put(1.2,-1.2){\vector(-1,1){1}}
}
\put(24,0){
\put(1,2){\makebox(0,0){\skelet}}
\put(0,0){\makebox(0,0){\Large$\bullet$}}
\put(-2.5,0){\makebox(0,0){$-$}}
\put(-.2,-.2){
\put(-1,-1){\makebox(0,0){\Large$\circ$}}
\put(-1.6,-2.2){\vector(1,2){.51}}
\put(-.4,-2.2){\vector(-1,2){.51}}
\put(-.85,-.85){\vector(1,1){.95}}
}
\put(0,-1){\vector(0,1){.81}}
\put(0,0.25){\line(0,1){.81}}
\put(4,0){\makebox(0,0)[l]{$+ \cdots$}}
}
}
}
\end{picture}
\]
\caption{\label{fig1}%
Action of $\delta$ on $\Gr_{\bullet\wc}^0$ -- the replacement
rule for a type $1$ vertex on the wheel.}
\end{figure} 
We claim that in order to
prove~(\ref{central}), it is enough to verify that
\begin{equation}
\label{centr}
\mbox {the map $\delta^0 : \Gr_{\bullet\wc}^0(d) \to 
\Gr_{\bullet\wc}^1(d)$ is a monomorphism for each $d \geq 1$.} 
\end{equation}

The spaces $\Gr_{\bullet\wc}^{m}(d,p)$ form a 
bicomplex $(\Gr_{\bullet\wc}^{*,*}(d),\delta)$ with
$\Gr_{\bullet\wc}^{p,q}(d) := \Gr_{\bullet\wc}^{p+q}(d,p)$ 
and $\delta$ the sum
$\delta^0 + \delta^1$, where $\delta^0 : \Gr_{\bullet\wc}^{*,*}(d) \to
\Gr_{\bullet\wc}^{*,*+1}(d)$ and $\delta^0 : \Gr_{\bullet\wc}^{*,*}(d) \to
\Gr_{\bullet\wc}^{*+1,*}(d)$ are defined above. 
Condition~(\ref{centr}) then implies~(\ref{central}) via a
standard spectral sequence argument. The only subtlety is that our
bicomplex is not a first quadrant one, thus the convergence of the
related spectral sequence has to be checked. 
We therefore decided to prove the implication \hbox{(\ref{centr})
$\Longrightarrow$~(\ref{central})} by  the following
elementary calculation.

Suppose that~(\ref{central}) does not hold
and let $x \in \Gr_{\bullet\wc}^0(d)$ be such that
$\delta(x) = 0$ while $x \not= 0$. There exists a
decomposition $x = x_a + x_{a+1} + \cdots + x_{a+s}$ with $x_w \in
\Gr_{\bullet\wc}^0(d,w)$ for $a \leq w \leq a+s$ in which 
$x_a \not=0$. Since $\delta^0(x_a)$
is the component of $\delta(x)$ in 
$\Gr_{\bullet\wc}^1(d,a)$, $\delta^0(x_a)
= 0$. Then~(\ref{centr}) implies $x_a = 0$, a contradiction.

Denote by $\osG_{\bullet\wc}^1(d,w) \subset \Gr_{\bullet\wc}^1(d,w)$ the
subspace spanned by graphs with one binary white vertex on the
wheel, as in the left graph in Figure~\ref{fig2}. Both
$\osG_{\bullet\wc}^1(d,w)$ and $\Gr_{\bullet\wc}^1(d,w)$ 
have canonical bases provided by
isomorphism classes of graphs, therefore one has a canonical
projection $\pi : \Gr_{\bullet\wc}^1(d,w) \to \osG_{\bullet\wc}^1(d,w)$. In addition to
the projection, there is a second map $r : \osG_{\bullet\wc}^1(d,w) \to
\Gr_{\bullet\wc}^0(d,w)$ whose definition is clear from Figure~\ref{fig2}.
\begin{figure}[t]
\[
\unitlength .4cm
\begin{picture}(20,2)(1,-1)
\put(7,0){
\put(1,2){\makebox(0,0){\skelet}}
\put(0,0){\makebox(0,0){\Large$\circ$}}
\put(-1,-1){\makebox(0,0){\Large$\bullet$}}
\put(-1.6,-2.2){\vector(1,2){.51}}
\put(-.4,-2.2){\vector(-1,2){.51}}
\put(-1.2,-1.2){\vector(1,1){1}}
\put(0,-1){\vector(0,1){.81}}
\put(0,0.25){\line(0,1){.81}}
\put(-1,-2.1){\makebox(0,0){\scriptsize $...$}}
\put(5,0){\makebox(0,0){$\stackrel{\mbox {\it r}}\longmapsto$}}
}
\put(15,0){
\put(1,2){\makebox(0,0){\skelet}}
\put(0,0){\makebox(0,0){\Large$\bullet$}}
\put(-1.2,-1.2){\vector(1,1){1}}
\put(0,-1){\vector(0,1){.81}}
\put(0,0.25){\line(0,1){.81}}
\put(-.5,-1.1){\makebox(0,0){\scriptsize $...$}}
}
\end{picture}
\]
\caption{\label{fig2}%
The map $r: \osG_{\bullet\wc}^1(d,w) \to
\Gr_{\bullet\wc}^0(d,w)$ contracts the unique edge connecting the
binary white vertex on the wheel with a black vertex outside the
wheel.}
\end{figure}

Let $G \in \Gr_{\bullet\wc}^0(d,w)$ be a graph. Observe that
$\Gr_{\bullet\wc}^0(d,0) = 0$, we may therefore assume $w \geq
1$. Recall that the differential $\delta(G)$ is the
sum~(\ref{Porad_mi_nejak_lidi_nepisou.})  of local replacements
$\delta_v(G)$ over $v \in \Vert(G)$.  Let $\Vert_\wc(G) \subset
\Vert(G)$ be the subset of vertices on the wheel.  For $v \in
\Vert_\wc(G)$, the contribution $\delta_v(G)$ contains precisely one
graph in $\Gr_{\bullet\wc}^1(d,w)$ with the binary white vertex -- see
again Figure~\ref{fig1}.  Denote this graph 
$\overline{\delta}\hskip .1em {}^0_v(G)$
and define $\overline{\delta}\hskip .1em {}^0(G) := \sum_{v \in \Vert_\wc(G)}
\overline{\delta}\hskip .1em {}^0_v(G)$.  
It is clear that $\Im(\overline{\delta}\hskip .1em {}^0)
\subset \osG\hskip .1em {}_{\bullet\wc}^1(d,w)$, 
$\overline{\delta}\hskip .1em {}^0 
= \pi \circ
\delta^0$ and $r \circ \overline{\delta}\hskip .1em {}^0 = w \cdot {\it id}$.
Combining these facts, we obtain $r \circ \pi \circ \delta\hskip .1em
{}^0 = w \cdot
{\it id}$, which implies~(\ref{centr}) and finishes the proof.
\end{proof}

We believe that one can even show  that the complex 
$(\Gr_{\bullet\wc}^*(d),\delta)$
used in the above proof is
acyclic in all dimensions. Let us close this section by formulating
the following interesting consequence of the proof of
Proposition~\ref{Aspon_Viktor_napsal.}.

\begin{corollary}
\label{C}
In stable dimensions, there are no nontrivial 
differential operators from vector fields to functions.
\end{corollary}

\begin{proof}
It is clear that $d$-multilinear operators from vector fields to functions are
described by the graph complex $\Gr_{\bullet\wc}^*(d)$ introduced in
our proof of
Proposition~\ref{Aspon_Viktor_napsal.}. Condition~(\ref{central})
implies that there are no nontrivial $d$-multilinear operators of this
type. The corollary
then follows from the standard (de)polarization trick.    
\end{proof}

\section{Structure of the space of natural operators}

In Example~\ref{Stava_se_ze_mne_celebrita.} we considered the 
trivial natural bundle  $\Re$ whose 
sections are smooth functions. Let $\gF$ be another natural
bundle. The space $\Nat(\gF,\Re)$ of natural operators $\gO : \gF \to
\Re$ with the `pointwise' multiplication is a commutative
algebra, with unit ${\mathfrak 1}$ the operator that sends all
sections of $\gF$ into the constant section $1 \in \bbR$. 
This indicates that spaces of natural operators may sometimes have a
rich algebraic structure that can be used to simplify their classification. 

\begin{definition}
\label{con}
We say that $\gF$ is a bundle with {\em \crr\/}
if the replacement rules send a connected graph to a linear
combination of connected graphs.
\end{definition}

All natural bundles considered in this paper have \crr, 
and the author does not know any `natural' natural
operator that has not. We will see that
the space of natural operators between bundles with \crr\ exhibits some
freeness property. Before we formulate the first statement of this
type, we introduce the following convention.

The graph complex $\Gr^*_{\gF,\Re}$ for operators
in $\Nat(\gF,\Re)$ is spanned by graphs with
vertices of the 1st type 
representing tensors in a prolongation of the fiber of $\gF$,
vertices~(\ref{Uz_nevim_co.}) of the third type and one 2nd type vertex
\raisebox{-.1em}{\ $\bbox$} which in this case has no inputs and no
outputs. Therefore \raisebox{-.1em}{\ $\bbox$} is an isolated vertex
bearing no information and we discard it from the picture. With this
convention, graphs spanning $\Gr^*_{\gF,\Re}$ have vertices of the 1st
and 3rd type only. The disjoint union of graphs spanning
$\Gr^*_{\gF,\Re}$ translates into the
pointwise multiplication of the corresponding operators and the 
unit ${\mathfrak 1} \in \Nat(\gF,\Re)$ is
represented by the `exceptional' empty graph.

\begin{theorem}
\label{X}
Let $\gF$ be a natural bundle with \crr. Then, in stable dimensions,
the commutative unital algebra $\Nat(\gF,\Re)$ is free, generated by the
subspace $\Nat_1(\gF,\Re)$ of natural operators represented by
connected graphs. In other words,
$\Nat(\gF,\Re) \cong \bbR[\Nat_1(\gF,\Re)]$,
where $\bbR[-]$ denotes the polynomial algebra functor.
\end{theorem}

\begin{proof}
Each graph spanning $\Gr^*_{\gF,\Re}$ decomposes
into the disjoint union of its connected
components. The differential
$\delta$,  by assumption, preserves this decomposition which is
clearly unique up to the order of components.
The proof is finished by recalling
that the disjoint union of graphs expresses the
pointwise multiplication of operators.
\end{proof}

Let  $\gF,\gG$ be natural bundles. 
The pointwise multiplication makes the space $\Nat(\gF,\gG)$ a
unital module over the unital algebra $\Nat(\gF,\Re)$.
We prove a structure theorem also for this space.

\begin{theorem}
\label{Y}
Suppose that both $\gF$ and $\gG$ are bundles with \crr. Then, in
stable dimensions,
$\Nat(\gF,\gG)$ is the free $\Nat(\gF,\Re)$-module generated by the
subspace $\Nat_1(\gF,\gG)$ of operators represented by connected graphs,
\begin{equation}
\label{pro}
\Nat(\gF,\gG) \cong \Nat_1(\gF,\gG) \otimes \Nat(\gF,\Re).
\end{equation}
\end{theorem}

\begin{proof}
The proof is similar to the proof of Theorem~\ref{X}.
The graph complex $\Gr^*_{\gF,\gG}$ describing operators
in  $\Nat(\gF,\gG)$ is spanned by graphs with
vertices of the first and third types,
and one vertex of the second type.
Each such a graph is the disjoint union of its connected
components as in~(\ref{g}) and the differential
preserves this  decomposition. Precisely one of these components
contains the vertex of the third type thus representing an operator in
$\Nat_1(\gF,\gG)$. The remaining components
describe operators from $\Nat_1(\gF,\Re)$ and assemble, via the
pointwise multiplication, into an operator in $\Nat(\gF,\bbR)$.  
\end{proof}

Theorems~\ref{X} and~\ref{Y} imply that in order to classify operators
in $\Nat(\gF,\gG)$, it is enough to understand the `connected'
subspaces $\Nat_1(\gF,\Re)$ and $\Nat_1(\gF,\gG)$. We will use this fact
in the next section.

\begin{example}
In Section~\ref{5} we studied natural operators on vector fields with
values in vector fields, that is, operators in
$\Nat(T^\ti,T) := \bigcup_{d \geq 0}\Nat(T^\td,T)$. We also considered
operators with values in functions and proved, 
in Corollary~\ref{C}, that there are no
nontrivial operators of this type in stable dimensions.

This means that $\Nat(T^\ti,\Re)$ is the trivial commutative algebra
$\bbR$ and~(\ref{pro}) reduces to the  isomorphism 
$\Nat(T^\ti,T) \cong \Nat_1(T^\ti,T)$ which says that all operators
from vector fields to vector fields live, in stable dimensions, on
connected graphs.
\end{example}

\section{Operators on connections and vector fields}
\label{7}

We will consider operators acting on a linear connection $\Gamma$
and a finite number of vector fields $X,Y,Z,\ldots$, with values
in vector fields, such as the covariant derivative 
$\nabla_XY$, torsion $T(X,Y)$
and curvature $R(X,Y)Z$ recalled in Example~\ref{cov}.
By Theorems~\ref{X} and~\ref{Y}, the structure of the space 
$\Nat(\Con \times T^{\times \infty},T) := \bigcup_{d \geq 0}
\Nat(\Con \times T^{\times d},T)$
of these operators is determined by the
`connected' subspaces $\Nat_1(\Con \times T^{\times \infty},T)$ and
$\Nat_1(\Con \times T^{\times \infty},\Re)$. 
In this section we describe these spaces.
The following remark should be compared to Remark~\ref{pol1} in 
Section~\ref{5}. 

\begin{remark}
\label{pol2}
 The local formula $O$ for a natural differential operator $\gO$ in
$\Nat(\Con \times T^{\times \infty},T)$ or in $\Nat(\Con \times
T^{\times \infty},\Re)$ decomposes into 
$O = \sum_{a,b \geq 0} O_{a,b}$ (finite sum),
where $O_{a,b}$ is the part of $O$ containing precisely $a$
$\nabla$-variables and $b$ vector field variables. For example, the
local formula~(\ref{Sergej_je_magor.}) for the covariant derivative
represented by the graph in~(\ref{covar}) is the sum $O_{1,2} +
O_{0,2}$, where $O_{1,2}(X,Y,\Gamma) := \Gamma^i_{jk}X^j Y^k 
\pa/\pa{x^i}$ and 
$O_{0,2}(X,Y,\Gamma) := X^j {Y^i_j}\pa/\pa{x^i}$. 

In contrast to Section~\ref{5}, here the action
of the structure group $\GL {\infty} n$ on the typical fiber is linear
only in the vector-field variables -- the non-linearity in the
$\nabla$-variables is manifested in the presence of the `isolated'
white vertex in the replacement
rule~(\ref{Chce_se_mi_curat.}). Nevertheless, one may still decompose 
$\gO =  \gO_1 + \cdots + \gO_r$,
with $\gO_k$ the operator represented by the local formula $O_d :
= \sum_{a\geq 0} O_{a,d}$, $1 \leq d \leq r$. 
Therefore {\em homogeneity\/} and {\em
multilinearity\/} in this section always refer to the vector fields
variables. 
\end{remark}

The first half of this section will be devoted to the study of the
space $\Nat_1(\Con \times T^{\times \infty},T)$, the
space $\Nat_1(\Con \times T^{\times \infty},\Re)$ will be addressed
in the second half of this section.  
As in Section~\ref{5}, we start with multilinear operators.

\begin{theorem}
\label{P1}
Let $d \geq 0$. On smooth manifolds of dimension $\geq 2d-1$, 
each $d$-multilinear operator in
$\Nat_1(\Con \times T^\od,T)$  is a linear
combination of iterations of the covariant derivative and the Lie
bracket which contains each of the vector fields $\Rada X1d$ exactly
once. All relations follow from the
anticommutativity and the Jacobi identity of the Lie bracket.

If $g_d$ denotes the number of linearly independent operators of
this type, the generating function
$g(t) = \sum_{d \geq 1} \frac 1{d!}\ g_dt^d$
is determined by the functional equation
\begin{equation}
\label{func}
e^{g(t)} \left(1-t-g^2(t)\right) = 1.
\end{equation}
\end{theorem}

Equation~(\ref{func}) can be expanded into inductive 
formula~(\ref{expansion}) from which one can calculate some initial
values of $g_k$ as  $g_1 = 1$, $g_2 = 3$, $g_3 = 26$, \&c.
Theorem~\ref{P1} will follow from Proposition~\ref{P11} below. The
depolarization of Theorem~\ref{P1} is:

\begin{corollary}
\label{T2}
On a smooth manifold $M$, 
each operator from $\Nat_1(\Con \times T^\ti,T)$ whose all components
are of homogeneity $\leq  \frac 12 (\dim(M) + 1)$ is a linear
combination of compositions of the covariant derivative and the Lie
bracket. All relations between these compositions follow from the
anticommutativity and the Jacobi identity of the Lie bracket.
\end{corollary}

The central object will be the 
subcomplex  $\Gr^*_{\bullet\nabla1}(d)$ of the graph complex 
$\Gr^*_{\Con \times T^{\times d},T}$
describing `connected' $d$-multilinear operators. Its degree $m$ piece
$\Gr^m_{\bullet\nabla1}(d)$ is
spanned by connected graphs with $d$ vertices~(\ref{tyden_piti}) labelled by
$\Rada X1d$, some number of vertices~(\ref{Uz_nevim_co.}) labelled
$\nabla$, $m$ white
vertices~(\ref{O_byly_po_vsich_muziky/a_choraly_nam_hraly_temne}) 
and one vertex
$\unitlength .25cm
\begin{picture}(1,1.4)(-1,-.7)
\put(-.45,.55){\makebox(0,0)[cc]{$\sbbox$}}
\put(-.5,-.8){\vector(0,1){1.2}}
\end{picture}$. 
It is clear that $\Gr^*_{\bullet\nabla1}(d)$ is precisely the
subcomplex spanned by connected graphs, of the direct sum 
$\Gr^*_{\bullet\nabla}(d) 
:= \bigoplus_{c \geq 0} \Gr^*_{\bullet(d)\nabla(c)}$,
where $\Gr^*_{\bullet(d)\nabla(c)}$ is
the graph complex of \cite[Corollary~5.1]{markl:ig}.
As in Proposition~\ref{zpatky_v_Koline}, 
one easily sees that the collection $\Gr^0_{\bullet\nabla1} =
\{\Gr^0_{\bullet\nabla1}(d)\}_{d \geq 1}$ forms an operad. It is also
not difficult to verify that each graph spanning
$\Gr^m_{\bullet\nabla1}(d)$ has at most $2d+m-1$ edges, which explains
the stability condition in Theorem~\ref{P1}.

Let $\calP = \{\calP(d)\}_{d\geq 1}$\label{psano_v_IHES} 
be the operad describing
algebras with two independent operations --  a~bilinear product $\star$
satisfying no other conditions and a Lie bracket. Of course, $\calP$ is the
free product (= the coproduct in the category of operads, 
see~\cite[p.~137]{markl:ha})
of the free operad  $\Gamma(\star)$
generated by the bilinear operation $\star$ and the operad $\Lie$ for
Lie algebras, $\calP =  \Gamma(\star)* \Lie$. 
Recall that we denoted by $\beta \in \Lie(2)$ the generator.

Define the operad homomorphism $F: \calP \to \Gr^0_{\bullet\nabla1}$
by $F(\beta) := \bfb$ and $F(\star) : = \bfc$, where $\bfb \in
\Gr^0_{\bullet\nabla1}(2)$ is the
graph~(\ref{Pozitri_snad_prijede_PP.}) representing the Lie bracket
and $\bfc \in \Gr^0_{\bullet\nabla1}(2)$ the graph~(\ref{covar}) for
the covariant derivative. As in Section~\ref{5} we easily see that $F$
is well-defined and that $\Im(F) \subset \Ker(\delta :
\Gr^0_{\bullet\nabla1} \to \Gr_{\bullet\nabla1}^1)$.  Theorem~\ref{P1} clearly
follows from

\begin{proposition}
\label{P11}
The map $F: \calP \to \Gr^0_{\bullet\nabla1}$ induces an isomorphism
$\calP \cong H^0(\Gr^0_{\bullet\nabla1},\delta)$. The generating
function $p(t) := \sum_{d \geq 1} \frac 1{d!}
{\dim(\calP(d))}\cdot t^d$ for the operad $\calP$
satisfies~(\ref{func}).
\end{proposition}

\begin{proof}
The map $F$ embeds into the following diagram of 
operads and their homomorphisms:
\begin{equation}
\label{diagram}
\odrazka{-2.3em}
\unitlength1.5cm
\begin{picture}(3,.7)(-1.5,.5)
\put(-2,1){\makebox(0,0)[cc]{$\calP = \Gamma(\star) * \Lie$}}
\put(-0,1){\makebox(0,0)[cc]{$\Gr_{\bullet\nabla1}^0$}}
\put(2,1){\makebox(0,0)[cc]{$\Gamma(\star) * \pLie$}}
\put(2,0){\makebox(0,0)[cc]{\hphantom{.}$\Gamma(\star) * \pLie$.}}
\put(-2,0){\vector(1,0){3.2}}
\put(-1,1){\vector(1,0){.56}}
\put(.35,1){\vector(1,0){.8}}
\put(-2,.80){\line(0,-1){.8}}
\put(2,.80){\vector(0,-1){.6}}
\put(-.7,1.1){\makebox(0,0)[bc]{\scriptsize $F$}}
\put(.7,1.1){\makebox(0,0)[bc]{\scriptsize $A$}}
\put(2.1,.5){\makebox(0,0)[lc]{\scriptsize $T$}}
\put(-.2,.1){\makebox(0,0)[bc]{\scriptsize ${\it id} * \iota$}}
\end{picture}
\end{equation}
Let us define the remaining maps in~(\ref{diagram}).
As in~\cite{chapoton-livernet:pre-lie}, one can
show that the operad 
$\Gr^0_{\bullet\nabla1}$ is isomorphic to the operad $\Gamma(\star) *
\pLie$ governing structures consisting of a bilinear multiplication
$\star$ and an independent 
pre-Lie product $\circ$. The map $A :
\Gr^0_{\bullet\nabla1} \to \Gamma(\star) * \pLie$ in~(\ref{diagram}) is 
the isomorphism that sends the graph
\[
\raisebox{-2em}{\rule{0pt}{0pt}}
\hskip 4.8cm
\unitlength .34cm
\begin{picture}(.7,1.9)(6,-.2)
\put(2,0){\makebox(0,0)[cl]{$\in \Gr_{\bullet\nabla1}^0(2)$}}
\put(-.5,0){\makebox(0,0)[cc]{\Large$\nabla$}}
\put(-.5,.5){\vector(0,1){1}}
\put(.3,-1.65){\vector(-1,2){.6}}
\put(1.07,-1.3){\vector(-1,1){1.2}}
\put(-.45,1.4){\makebox(0,0)[bc]{$\bbox$}}
\put(0.2,-1.5){\makebox(0,0)[cc]{\Large$\bullet$}}
\put(0.2,-2){\makebox(0,0)[tc]{\scriptsize $X$}}
\put(1.5,-1.7){\makebox(0,0)[tc]{\scriptsize $Y$}}
\put(1,-1.3){\makebox(0,0)[cc]{\Large$\bullet$}}
\end{picture}
\]
into  $X\star Y \in \Gamma(\star)(2)$ and the graph
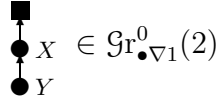
\[
\hskip 1em
\raisebox{-1.5em}{\rule{0pt}{0pt}}
\unitlength .5cm
\begin{picture}(0,1)(0,.9)
\put(0,2){\makebox(0,0)[cc]{\hskip .5mm$\bbox$}}
\put(0,1){\vector(0,1){.935}}
\put(0,1){\makebox(0,0)[cc]{\Large$\bullet$}}
\put(.4,1){\makebox(0,0)[lc]{\scriptsize$X$}}
\put(0,0){\vector(0,1){.935}}
\put(0,0){\makebox(0,0)[cc]{\Large$\bullet$}}
\put(.4,0){\makebox(0,0)[lc]{\scriptsize$Y$}}
\end{picture}
\hskip 1.5em \in \Gr^0_{\bullet\nabla1}(2)
\]
into $X \circ Y \in \pLie(2)$. The
map $T : \Gamma(\star) * \pLie \to \Gamma(\star) * \pLie$ is the `twist'
$T(X\star Y ) := X\star Y - Y \circ X$ and
$T(X\circ Y) := X \circ Y$.
It is evident that the composition $\it T A F$ coincides with the
coproduct ${\it id} * \iota$ of the identity ${\it id} : \Gamma(\star)
\to \Gamma(\star)$ and the map $\iota : \Lie
\to \pLie$ given by the antisymmetrization of the pre-Lie product
$\iota ([X,Y]) :=  Y \circ X - X \circ Y$,
which is an inclusion by~\cite[Proposition~3.1]{markl:JLT07}.
This implies that ${\it id}* \iota$ is a monomorphism, therefore  
$F$ is a monomorphism,~too. 

Now, to prove that $F$ induces an isomorphism $\calP \cong
H^*(\Gr^0_{\bullet\nabla1},\delta)$, it suffices to show that the
dimensions of the spaces $H^0(\Gr_{\bullet\nabla1}^0(d),\delta)$
and $\calP(d)$ are the same, for each $d \geq 1$.
Our calculation of the dimension of  
$H^0(\Gr_{\bullet\nabla1}^0(d),\delta)$ will be
based on the fact that $(\Gr_{\bullet\nabla1}^*(d),\delta)$ forms a
bicomplex. For integers $p,q$ denote by
$\Gr_{\bullet\nabla1}^{p,q}(d)$ the subspace of 
$\Gr_{\bullet\nabla1}^{p+q}(d)$ spanned by graphs
with precisely $-p$ $\nabla$-vertices. It immediately follows from
the replacement rules~(\ref{Co_Sergej_planuje?}),~(\ref{Nikdo}) 
and~(\ref{Chce_se_mi_curat.}) that
$\delta = \delta' + \delta''$, where
$\delta'(\Gr_{\bullet\nabla1}^{p,q}(d)) \subset  
\Gr_{\bullet\nabla1}^{p+1,q}(d)$ and
$\delta''(\Gr_{\bullet\nabla1}^{p,q}(d)) \subset  
\Gr_{\bullet\nabla1}^{p,q+1}(d)$.
It is also clear from simple graph combinatorics that the bicomplex
$(\Gr_{\bullet\nabla1}^{*,*}(d),\delta)$ is bounded by the triangle
$p=0$, $p+q = 0$ and $q=d-1$, see Figure~\ref{fig3}.
\begin{figure}[t]
\[
\unitlength 1.3cm
\begin{picture}(3,1.8)(-3.5,.2)
\put(0,0){\makebox(0,0)[cc]{$\Gr^{0,0}_{\bullet\nabla1}(3)$}}
\put(0,1){\makebox(0,0)[cc]{$\Gr^{0,1}_{\bullet\nabla1}(3)$}}
\put(0,2){\makebox(0,0)[cc]{$\Gr^{0,2}_{\bullet\nabla1}(3)$}}
\put(-2,1){\makebox(0,0)[cc]{$\Gr^{-1,1}_{\bullet\nabla1}(3)$}}
\put(-4,2){\makebox(0,0)[cc]{$\Gr^{-2,2}_{\bullet\nabla1}(3)$}}
\put(-2,2){\makebox(0,0)[cc]{$\Gr^{-1,2}_{\bullet\nabla1}(3)$}}
\put(0,0.28){\vector(0,1){.5}}
\put(0,1.28){\vector(0,1){.5}}
\put(-2,1.28){\vector(0,1){.5}}
\put(-1.3,1){\vector(1,0){.7}}
\put(-1.3,2){\vector(1,0){.7}}
\put(-3.3,2){\vector(1,0){.6}}
\put(0.1,.5){\makebox(0,0)[lc]{\scriptsize $\delta''$}}
\put(0.1,1.5){\makebox(0,0)[lc]{\scriptsize $\delta''$}}
\put(-1.9,1.5){\makebox(0,0)[lc]{\scriptsize $\delta''$}}
\put(-.9,2.1){\makebox(0,0)[bc]{\scriptsize $\delta'$}}
\put(-.9,1.1){\makebox(0,0)[bc]{\scriptsize $\delta'$}}
\put(-2.9,2.1){\makebox(0,0)[bc]{\scriptsize $\delta'$}}
\end{picture}
\]
\caption{\label{fig3}%
The bicomplex $(\Gr_{\bullet\nabla1}^{*,*}(3),\delta' + \delta'')$.}
\end{figure}
The horizontal differential $\delta'$ in $\Gr_{\bullet\nabla1}^{*}(d)$ 
is easy to describe -- it
replaces $\nabla$-vertices according the rule
\begin{equation}
\label{rule}
\unitlength .4cm
\begin{picture}(.7,1.4)(-.5,-.5)
\put(-.5,0){\makebox(0,0)[cc]{\Large$\nabla$}}
\put(-.5,.4){\vector(0,1){1}}
\put(.1,-1.25){\vector(-1,2){.5}}
\put(.67,-1.1){\vector(-1,1){1}}
\unitlength 3mm
\put(-21,-20){
\put(19.4,18){\vector(1,4){0.475}}
\put(18.9,18){\vector(1,2){.95}}
\put(16,18){\vector(2,1){3.8}}
\put(17.75,18){\makebox(0,0)[cc]{$\ldots$}}
\put(17.7,16.6){\makebox(0,0)[cc]{%
   $\underbrace{\rule{10mm}{0mm}}_{\mbox{\scriptsize $v$ inputs}}$}}
}
\end{picture}
\hskip .5cm \longmapsto \hskip .5cm - \hskip -.2cm
\unitlength 3mm
\begin{picture}(20,2)(15.5,18.8)
\put(20,20.3){\vector(0,1){1.4}}
\put(18,18){\vector(1,1){1.8}}
\put(22,18){\vector(-1,1){1.8}}
\put(19,18){\vector(1,2){.88}}
\put(20,20){\makebox(0,0)[cc]{\Large$\circ$}}
\put(20.5,18){\makebox(0,0)[cc]{$\ldots$}}
\put(20,16.6){\makebox(0,0)[cc]{%
   $\underbrace{\rule{12 mm}{0mm}}_{\mbox{\scriptsize $v+2$}}$}}
\end{picture}
\raisebox{-2.3em}{\rule{0pt}{0pt}} \hskip -3.8cm   , \ v \geq 0,
\end{equation}
and leaves other vertices unchanged.

\begin{remark}
\label{R1}
At this point we need to make a digression and observe that
$(\Gr^*_{\bullet\nabla1}(d),\delta')$ is a particular case of 
the following construction.
{}For each collection $(U^*,\vt_U) = \{(U^*(s),\vt_U)\}_{s
\geq 2}$ of right dg-$\Sigma_s$-modules $(U^*(s),\vt_U)$, one may
consider the complex
$\Grd *1U = (\Grd *1{U^*},\vt)$ spanned by connected graphs with 
$d$ vertices~(\ref{tyden_piti}) labelled $\Rada X1d$, one vertex
$\unitlength .25cm
\begin{picture}(1,1.4)(-1,-.7)
\put(-.45,.55){\makebox(0,0)[cc]{$\sbbox$}}
\put(-.5,-.8){\vector(0,1){1.2}}
\end{picture}$ and a finite number of vertices decorated by elements
of $U$. The grading of $\Grd *1{U^*}$ 
is induced by the grading of $U^*$ and the differential $\vt$ replaces
$U$-decorated vertices, one at a time, by their $\vt_U$-images and leaves other
vertices unchanged.
It is a standard fact~\cite{mv} (see also~\cite[Theorem~21]{markl:ha}) 
that the assignment 
$(U^*,\vt_U) \mapsto (\Grd *1{U^*},\vt)$ 
is a polynomial, hence exact, functor, so
\begin{equation}
\label{21}
H^*(\Grd *1{U^*},\vt) \cong \Grd *1{H^*(U,\vt_U)}.
\end{equation}

Let now $(E^*,\vt_E) = \{(E^*(s),\vt_E)\}_{s \geq 2}$ 
be such that $E^0(s)$ is spanned by 
symbols~(\ref{Uz_nevim_co.}),  with $v+2 = s$, $E^1(s)$ by
symbols~(\ref{O_byly_po_vsich_muziky/a_choraly_nam_hraly_temne}) with
$u = s$, and $E^m(s) = 0$ for $m \geq 2$. The
differential $\vt_E$ is defined by replacement rule~(\ref{rule}). 
More formally, $E^0(s) = {\rm Ind}^{\Sigma_s}_{\Sigma_{s-2}} 
({{\mathbf 1}}_{s-2})$ and  $E^1(s) = {{\mathbf 1}}_s$, where
${{\mathbf 1}}_{s-2}$ (resp.~${{\mathbf 1}}_{s}$) denotes the trivial
representation of the symmetric group $\Sigma_{s-2}$
(resp.~$\Sigma_s$). The differential $\vt_E$ then sends the generator $1
\in {{\mathbf 1}}_{s-2}$ into   $-1 \in {{\mathbf 1}}_{s}$.
It is clear that, with this particular choice of the collection
$(E^*,\vt_E)$,
\begin{equation}
\label{22}
(\Gr_{\bullet\nabla1}^*(d),\delta') \cong (\Grd *1{E^*},\vt).
\end{equation}
\end{remark}

Let us continue with the proof of
Proposition~\ref{P11}. Equations~(\ref{21}) and~(\ref{22}) in
Remark~\ref{R1} imply that
\begin{equation}
\label{Teplo_jako_na_jare.}
H^*(\Gr_{\bullet\nabla1}^*(d),\delta') = \Grd *1{H^*(E,\vt_E)}.
\end{equation}
Since $\vt_E : E^0(s) \to E^1(s)$ is an epimorphism, the collection
$H^*(E,\vt_E) = \{H^*(E(s),\vt_E)\}_{s \geq 2}$ is concentrated in degree 
$0$ and $H^0(E(s),\vt_E)$ 
is the kernel of the map $\vt_E : E^0(s) \to E^1(s)$. 
We conclude that $\Grd *1{H^*(E,\vt_E)}$ is spanned
by graphs with $d$ vertices~(\ref{tyden_piti}) 
labelled $\Rada X1d$, one vertex 
$\unitlength .25cm
\begin{picture}(1,1.4)(-1,-.7)
\put(-.45,.55){\makebox(0,0)[cc]{$\sbbox$}}
\put(-.5,-.8){\vector(0,1){1.2}}
\end{picture}$
and some number of vertices decorated by the collection 
$H^0(E,\vt_E) = \{H^0(E(s),\vt_E)\}_{s \geq 2}$. 

In particular, the graded space $\Grd *1{H^*(E,\vt_E)}$ and hence,
by~(\ref{Teplo_jako_na_jare.}), also the horizontal cohomology
$H^*(\Gr_{\bullet\nabla1}^*(d),\delta')$, is concentrated in degree $0$. This implies
that the first term $(E^{p,q}_1,d_1) = (H^p(\Gr_{\bullet\nabla1}^{*,q},\delta'),d_1)$ 
of the corresponding spectral sequence is
supported by the diagonal $p + q =0$, so this spectral
sequence degenerates at this level and 
\[
\dim(H^0(\Gr_{\bullet\nabla1}^*(d),\delta)) =
\dim(H^0(\Gr_{\bullet\nabla1}^*(d),\delta')) = \dim(\Grd 01{H^0(E,\vt_E)}.
\]

Denote the common value of the dimensions in the above display
$g_d$. We claim that the sequence $\{g_d\}_{d \geq 1}$ satisfies 
the recursion:
\begin{eqnarray}
\label{expansion}
\frac{g_{n+1}}{(n+1)!}
&=&
\frac{g_n}{n!} + 
{\frac 1{2!}} \sum_{i+j = n} \frac{g_ig_j}{i!j!} 
+\frac 1{3!} \sum_{i+j +k = n} \frac{g_ig_jg_k}{i!j!k!}
+{\frac {1}{4!}} \sum_{i+j +k + l =n} \frac{g_ig_jg_kg_l}{i!j!k!l!}+ \cdots
\\
\nonumber 
&& + \frac {2(2-1) - 1}{2!} \sum_{i+j = n+1} \frac{g_ig_j}{i!j!} 
+\frac {3(3-1)-1}{3!} 
\sum_{i+j +k= n +1} \frac{g_ig_jg_k}{i!j!k!} + \cdots.
\end{eqnarray}
This can be seen as follows.
Graphs $G$ spanning $\Grd 01{H^0(E,\vt_E)}$ are rooted trees with a
distinguished vertex (= root) $\unitlength .25cm
\begin{picture}(1,1.4)(-1,-.7)
\put(-.45,.55){\makebox(0,0)[cc]{$\sbbox$}}
\put(-.5,-.8){\vector(0,1){1.2}}
\end{picture}$. The vertex of $G$ adjacent to the root might either be
a vertex~(\ref{tyden_piti}) or a vertex decorated by $H^0(E,\vt_E)$. The
contribution from trees of the first type is reflected by the
first line of~(\ref{expansion}), in which the coefficients
$1,1/{2!},1/{3!},\ldots$ equal
${\dim({\mathfrak 1}_s)}/{s!}$, $s \geq 1$, 
where ${\mathfrak 1}_s$ is the trivial representation of the 
symmetric group $\Sigma_s$ spanned by the vertex~(\ref{tyden_piti}) with $u=s$.
The second line of~(\ref{expansion}) counts contributions from
trees of the second type. The coefficients are
${\dim(H^0(E(s),\vt_E))}/{s!}$, $s \geq 2$. It is simple to
assemble~(\ref{expansion}) into equation~(\ref{func}).

Let us show that the generating function $p(t) := \sum_{d \geq 1}
\frac 1{d!} {\dim(\calP(d))}\cdot t^d$ for the operad $\calP$ also
satisfies~(\ref{func}).  Since $\calP$ is, as the coproduct of
quadratic Koszul operads, itself quadratic Koszul, one has the
functional equation~\cite[Theorem~3.3.2]{ginzburg-kapranov:DMJ94}:
\begin{equation}
\label{related}
q(- p(t)) = -t.
\end{equation}
relating $p$ with the generating function $q(t) := \sum_{d \geq
1}\frac 1{d!}  {\dim(\calQ(d))}\cdot t^d$ of its quadratic dual
$\calQ$.

For convenience of the reader, we make a digression and briefly recall
the definition of quadratic operads and their quadratic duals. Details can
be found in~\cite[II.3.2]{markl-shnider-stasheff:book} or in the
original source~\cite{ginzburg-kapranov:DMJ94}.  An operad $\calA$ is
{\em quadratic} if it is the quotient $\Gamma(E)/(R)$ of the free operad
$\Gamma(E)$ on the right $\Sigma_2$-module $E := \calA(2)$ of 
arity-two operations of $\calA$, modulo the operadic ideal $(R)$ generated
by some subspace $R \subset \Gamma(E)(3)$.

Each quadratic operad $\calA = \Gamma(E)/(R)$ 
as above has its {\em quadratic dual\/}
$\calA^!$~\cite[Definition~II.3.37]{markl-shnider-stasheff:book}
defined as follows. Let us denote $E^\vee := E^* \otimes {\rm
sgn\/}_2$ the linear dual of the right $\Sigma_2$-module $E$ twisted
by the signum representation. One then has a natural isomorphism
$\Gamma(E^\vee)(3) \cong \Gamma(E)(3)^*$ of right $\Sigma_3$-modules. 
Let $R^{\perp }\subset
\Gamma(E^{\vee })(3)$ denote the annihilator of $R$ in
$\Gamma(E^\vee)(3) \cong \Gamma(E)(3)^*$. The quadratic dual of
$\calA$ is the quotient $\calA^! := \Gamma(E^\vee)/(R^\perp)$.

To describe the quadratic dual $\calQ$ of the operad $\calP$
introduced on page~\pageref{psano_v_IHES} is an easy task.  The operad
$\calQ$ governs algebras $V$ with two bilinear operations, $\ccdot$
and $\ast$, such that $\ccdot$ is commutative associative, $\ast$ is
`nilpotent' $(a \ast b) \ast c = a \ast (b \ast c) = 0$, $a,b,c \in
V$, and these two operations annihilate each other: $(a \ccdot b) \ast
c = a \ast (b \ccdot c) = a \ccdot (b \ast c) = 0$, $a,b,c \in V$.  It
is immediately obvious that
\[
\dim(\calQ(1)) = 1,\ \dim(\calQ(2)) =  3 \ \mbox { and }
\dim(\calQ(d)) = \dim(\Com(d)) = 1 \mbox { for } d \geq 3, 
\]
where $\Com$ denotes the operad for commutative associative algebras. 
The generating function for $\calQ$ therefore equals $q(t) = e^t - 1
+ t^2$ and equation~(\ref{related}) gives
$e^{-p(t)} - 1 + p(t)^2 = -t$, 
which is equivalent to~(\ref{func}). We proved that the
generating functions $g(t)$ and $p(t)$ satisfy the same functional equation
and, by definition, 
the same initial condition $p(0) = g(0) = 0$, therefore they
coincide and $\dim (H^0(\Gr_{\bullet\nabla1}^0(d),\delta)) 
= \dim(\calP(d))$ for each
$d \geq 1$.
\end{proof}

In the rest of this section we study operators in $\Nat_1(\Con \times
T^{\times \infty},\Re)$. Roughly speaking, we prove that all operators
in this space are traces in the following sense. 
Let $\gO \in \Nat(\Con \times T^\ti,T)$
 be an operator acting on vector
fields $X_0,X_1,X_2,\ldots$ and a connection $\Gamma$. 
Suppose that $\gO$ is a linear order $0$ differential operator in
 $X_0$. This means that the local formula $O(X_0,X_1,X_2,\ldots,\Gamma) \in
 \bbR$ for $\gO$ is a linear function of $X_0$ and does not contain
 derivatives of $X_0$.
For such an operator we define 
$\Tr_{X_0}(\gO) \in \Nat(\Con \times T^\ti,\Re)$ by the local formula
\[
\Tr_{X_0}(O)(X_1,X_2,\ldots,\Gamma) := \mbox{Trace}( 
O (-,X_1,X_2,\ldots,\Gamma) : \Rn \to \Rn) \in \bbR.
\] 
It is easy to see that $\Tr_{X_0}(\gO)$  is well defined.
Let us formulate a structure theorem for multilinear operators from 
$\Nat_1(\Con \times T^{\times \infty},\Re)$.

\begin{theorem}
\label{T3}
Let $d \geq 0$. On smooth manifolds of dimension $\geq 2d$, 
each $d$-multilinear 
operator in $\Nat_1(\Con \times T^{\ot d},\Re)$ is the trace of
a $(d+1)$-multilinear operator from  $\Nat_1(\Con \times T^{\ot (d+1)},T)$.
\end{theorem}

Theorem~\ref{T3} will follow from Proposition~\ref{P12} below. A
depolarized version of Theorem~\ref{T3} is:

\begin{corollary}
\label{C2}
On a smooth manifold $M$, 
each operator from $\Nat_1(\Con \times T^\ti,\Re)$ whose all components
are of homogeneity $\leq  \textstyle\frac 12 \dim(M)$ is a trace of an
operator from $\Nat_1(\Con \times T^\ti,T)$.
\end{corollary}

Denote by $\Gr^*_{\bullet\nabla\wc}(d)$ 
the graph complex describing operators in
$\Nat_1(\Con \times T^{\ot d},\Re)$. The
degree $m$-component of this complex is
spanned by connected graphs with $d$ vertices~(\ref{tyden_piti}) 
labelled $\Rada X1d$,  some number of vertices~(\ref{Uz_nevim_co.}) labelled
$\nabla$ and $m$ white
vertices~(\ref{O_byly_po_vsich_muziky/a_choraly_nam_hraly_temne}). It
is not difficult to see that the number of edges of
graphs spanning $\Gr^0_{\bullet\nabla\wc}(d)$ is $\leq 2d$, 
which explains the stability
assumption in Theorem~\ref{T3}.

We will also consider the subcomplex $\Grtr^* (d) \subset
\Gr_{\bullet\nabla1}^*(d+1)$ of graphs describing operators in
$\Nat_1(\Con \times T^{\otimes (d+1)},T)$ for which the trace is
defined. Clearly, the degree $m$ component $\Grtr^m (d)$ of this
subcomplex is spanned by connected graphs with one vertex $\unitlength
.25cm
\begin{picture}(1,1.4)(-1,-.7)
\put(-.45,-.55){\makebox(0,0)[cc]{$\bullet$}}
\put(-.5,1){\vector(0,-1){1.4}}
\end{picture}$
labelled $X_0$, one vertex  $\unitlength .25cm
\begin{picture}(1,1.4)(-1,-.7)
\put(-.45,.55){\makebox(0,0)[cc]{$\sbbox$}}
\put(-.5,-.8){\vector(0,1){1.2}}
\end{picture}$,
$d$ vertices~(\ref{tyden_piti}) labelled $\Rada X1d$,  
a finite number of vertices~(\ref{Uz_nevim_co.}) labelled
$\nabla$ and $m$ white
vertices~(\ref{O_byly_po_vsich_muziky/a_choraly_nam_hraly_temne}). The
trace is represented by the map
$\Tr : \Grtr^* (d) \to  \Gr^*_{\bullet\nabla\wc}(d)$ 
that removes the vertices 
$\unitlength .25cm
\begin{picture}(1,1.4)(-1,-.7)
\put(-.45,-.55){\makebox(0,0)[cc]{$\bullet$}}
\put(-.5,1){\vector(0,-1){1.4}}
\end{picture}$
and $\unitlength .25cm
\begin{picture}(1,1.4)(-1,-.7)
\put(-.45,.55){\makebox(0,0)[cc]{$\sbbox$}}
\put(-.5,-.8){\vector(0,1){1.2}}
\end{picture}$
and connects the two loose edges created in this way by
a directed wheel. 
It is clear that this map commutes with the differentials.
We now establish Theorem~\ref{T3} by proving the following.

\begin{proposition}
\label{P12}
The map $\Tr : (\Grtr^* (d),\delta) \to
(\Gr^*_{\bullet\nabla\wc}(d),\delta)$ induces an epimorphism of cohomology 
$H^0(\Grtr^* (d),\delta) \to H^0(\Gr^*_{\bullet\nabla\wc}(d),\delta)$.
\end{proposition}

\begin{proof}
As in the proof of Proposition~\ref{P11} we observe that both
$(\Grtr^* (d),\delta)$ and $(\Gr^*_{\bullet\nabla\wc}(d),\delta)$ are
bicomplexes, with $\Grtr^{p,q} (d)$
(resp.~$(\Gr^{p,q}_{\bullet\nabla\wc}(d)$) spanned by graphs in
$\Grtr^{p+q} (d)$ (resp.~$(\Gr^{p+q}_{\bullet\nabla\wc}(d)$) with
precisely $-p$ $\nabla$-vertices. The differential in both complexes
decomposes as $\delta = \delta' + \delta''$ where $\delta'$ (the
`horizontal part') raises the $p$-degree by one and preserves the
$q$-degree, and $\delta''$ (the `vertical part') preserves the
$q$-degree and raises the $p$-degree by one.

The map $\Tr : (\Grtr^* (d),\delta) \to
(\Gr^*_{\bullet\nabla\wc}(d),\delta)$ obviously preserves the bigradings,
therefore it induces the map 
\begin{equation}
\label{24}
H^*(\Tr,\delta') : H^*(\Grtr^* (d),\delta') 
\to H^*(\Gr^*_{\bullet\nabla\wc}(d),\delta')
\end{equation}
of the horizontal cohomology. Using the same considerations as in the
proof of Proposition~\ref{P11}, we identify this map with
\begin{equation}
\label{E1}
\Tr : \Grd *{\Tr}{H^*(E,\vt_E)} \to \Grd *{\wc}{H^*(E,\vt_E)},
\end{equation}
where $(E^*,\vt_E)$ is the dg-collection introduced in Remark~\ref{R1}
and the graph complexes in~(\ref{E1}) are defined analogously as 
the graph complex $\Grd *{1}{H^*(E,\vt_E)}$ used in the proof of
Proposition~\ref{P11}.

Let us show that 
the map in~(\ref{E1}) is an epimorphism. Consider a graph $G$
in $\Grd *{\wc}{H^*(E,\vt_E)}$ and choose a directed edge $e$ in the
(unique) wheel of $G$. Let $\widehat G$ be the graph in
$\Grd *{\Tr}{H^*(E,\vt_E)}$ obtained by cutting $e$ in the middle
and decorating the loose ends thus created by vertices $\unitlength .25cm
\begin{picture}(1,1.4)(-1,-.7)
\put(-.45,-.55){\makebox(0,0)[cc]{$\bullet$}}
\put(-.5,1){\vector(0,-1){1.4}}
\end{picture}$
and $\unitlength .25cm
\begin{picture}(1,1.4)(-1,-.7)
\put(-.45,.55){\makebox(0,0)[cc]{$\sbbox$}}
\put(-.5,-.8){\vector(0,1){1.2}}
\end{picture}$ as in the following display:
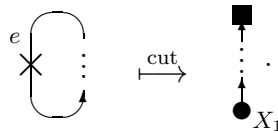
\[
\odrazka{-1em}
\unitlength .35cm
\begin{picture}(20,3.3)(1,-1.5)
\put(7,0){
\put(1,2){\makebox(0,0){\skelet}}
\put(0,-1){\line(0,1){2}}
\put(1.975,-1){\vector(0,1){0}}
\put(5,0){\makebox(0,0){$\stackrel{\rm cut}\longmapsto$}}
\put(-.4,1){\makebox(0,0)[r]{\scriptsize $e$}}
\put(0,0){\makebox(0,0)[]{\Large $\times$}}
}
\put(15,-.5){
\put(0,1.4){\vector(0,1){.8}}
\put(0,-1){\vector(0,1){1}}
\put(0,0.2){\makebox(0,0)[b]{$\vdots$}}
\put(0.05,1.9){\makebox(0,0)[b]{$\bbox$}}
\put(0,-.9){\makebox(0,0)[t]{\Large$\bullet$}}
\put(.4,-1.3){\makebox(0,0)[tl]{\scriptsize $X_1$}}
}
\end{picture}\hskip -1.8cm \raisebox{.51cm}{.}
\]
Clearly $\Tr(\widehat G) = G$ which proves that~(\ref{E1}) is surjective.
So, we have two spectral sequences, $(E^{p,q}_*,d_*)$ and
$(F^{p,q}_*,d_*)$, such that
\[
(E^{p,q}_0,d_0) =  (\Grtr^{p,q}(d),\delta),\
(F^{p,q}_0,d_0) =  (\Gr_{\bullet\nabla\wc}^{p,q}(d),\delta),
\] 
and the map $\Tr_* : (E^{p,q}_*,d_*) \to (F^{p,q}_*,d_*)$ induced by
the trace map $\Tr : \Grtr^* (d) \to \Gr^*_{\bullet\nabla\wc}(d)$. The
map $\Tr_1 : (E^{p,q}_1,d_1) \to (F^{p,q}_1,d_1)$ of the first levels
of the spectral sequences is~(\ref{24}) and we identified this map
with epimorphism~(\ref{E1}). It is also clear that the first terms of
both spectral sequences are supported by the diagonal $p+q=0$, so
these spectral sequences degenerate at this level. A standard argument
then implies that the map $H^0(\Tr) : H^0(\Grtr^* (d),\delta) \to
H^0(\Gr^*_{\bullet\nabla\wc}(d),\delta)$ in Proposition~\ref{P12} is
an epimorphism.
\end{proof}


\def\cprime{$'$}

\end{document}